\documentclass[preprint,12pt,authoryear]{elsarticle}

\usepackage{amssymb}
\usepackage{amsthm}
\usepackage{amsmath}
\usepackage{xcolor}
\usepackage{booktabs}
\usepackage{xspace}

\usepackage[a4paper,top=2.45cm,bottom=2.45cm,left=2.7cm,right=2.7cm,marginparwidth=2cm]{geometry}

\usepackage{tikz}
\usetikzlibrary{positioning}
\usetikzlibrary{arrows,decorations.markings,calc,patterns}

\usepackage[aboveskip=0pt,belowskip=-5pt]{subcaption}

\usepackage{enumitem}

\journal{enter the reviewing process}

\newcommand{\solvedS}{{\sc Solved}$^{\text{S}}$\xspace}
\newcommand{\foundref}{{\sc Found}$^{\text{R}}$\xspace}
\newcommand{\nearref}{{\sc Near}$^{\text{R}}$\xspace}
\newcommand{\alternat}{{\sc Altern}\xspace}
\newcommand{\timeF}{{\sc Time}$^{\text{BFR}}$\xspace}
\newcommand{\success}{{\sc Success}\xspace}

\newcommand{\harmonic}{\textup{\texttt{harmonic}}\xspace}
\newcommand{\daisy}{\textup{\texttt{daisy mamil}}\xspace}
\newcommand{\daisyF}{\textup{\texttt{daisy mamil3$^\text{F}$}}\xspace}
\newcommand{\daisyP}{\textup{\texttt{daisy mamil3$^\text{P}$}}\xspace}
\newcommand{\hiv}{\textup{\texttt{hiv}}\xspace}
\newcommand{\hivF}{\textup{\texttt{hiv$^\text{F}$}}\xspace}
\newcommand{\hivP}{\textup{\texttt{hiv$^\text{P}$}}\xspace}
\newcommand{\lv}{\textup{\texttt{Lotka Volterra}}\xspace}
\newcommand{\lvF}{\textup{\texttt{Lotka Volterra$^\text{F}$}}\xspace}
\newcommand{\lvP}{\textup{\texttt{Lotka Volterra$^\text{P}$}}\xspace}
\newcommand{\FHN}{\textup{\texttt{FHN}}\xspace}
\newcommand{\alphapinene}{\textup{\texttt{alpha pinene}}\xspace}
\newcommand{\crauste}{\textup{\texttt{Crauste}}\xspace}
\newcommand{\crausteF}{\textup{\texttt{Crauste$^\text{F}$}}\xspace}
\newcommand{\crausteP}{\textup{\texttt{Crauste$^\text{P}$}}\xspace}
\newcommand{\BBG}{\textup{\texttt{BBG}}\xspace}

\newcommand{\BARON}{\textup{\texttt{BARON$^{\text{\sc g}}$}}\xspace}
\newcommand{\Couenne}{\textup{\texttt{Couenne$^{\text{\sc g}}$}}\xspace}
\newcommand{\Octeract}{\textup{\texttt{Octeract$^{\text{\sc g}}$}}\xspace}
\newcommand{\BONMIN}{\textup{\texttt{BONMIN$^{\text{\sc l}}$}}\xspace}
\newcommand{\Ipopt}{\textup{\texttt{Ipopt$^{\text{\sc l}}$}}\xspace}
\newcommand{\Knitro}{\textup{\texttt{Knitro$^{\text{\sc l}}$}}\xspace}
\newcommand{\CONOPT}{\textup{\texttt{CONOPT$^{\text{\sc l}}$}}\xspace}
\newcommand{\SNOPT}{\textup{\texttt{SNOPT$^{\text{\sc l}}$}}\xspace}

\newcommand{\trapecio}{\textup{\texttt{Trapezoid}}\xspace}
\newcommand{\adams}{\textup{\texttt{Adams-Moulton}}\xspace}
\newcommand{\runge}{\textup{\texttt{Runge-Kutta}}\xspace}
\newcommand{\simpson}{\textup{\texttt{Simpson}}\xspace}
\newcommand{\euler}{\textup{\texttt{Euler}}\xspace}

\newcommand{\modeA}{\textup{\texttt{Baseline}}\xspace}
\newcommand{\modeB}{\textup{\texttt{ExtraTol}}\xspace}
\newcommand{\modeC}{\textup{\texttt{SoftCons}}\xspace}

\newcommand{\xb}{\pmb{x}}
\newcommand{\yb}{\pmb{y}}

\newcommand{\pb}{\pmb{p}}
\newcommand{\kb}{\pmb{k}}

\usepackage[final,bookmarksnumbered,breaklinks,colorlinks]{hyperref}
\hypersetup{urlcolor=blue,linkcolor=blue,citecolor=blue}

\begin{document}

\begin{frontmatter}



\title{Parameter estimation in ODEs: assessing the potential of local and global solvers}


\author[inst1]{M. Fernández de Dios}
\author[inst1,inst2]{Ángel M. González-Rueda}
\author[inst3]{Julio R. Banga}
\author[inst1,inst2]{Julio González-Díaz}
\author[inst3]{David R. Penas}

\affiliation[inst1]{organization=Department of Statistics, Mathematical Analysis and Optimization and MODESTYA Research Group. University of Santiago de Compostela. Santiago de Compostela, Spain}

\affiliation[inst2]{organization=CITMAga (Galician Center for Mathematical Research and Technology), 15782 Santiago de Compostela, Spain}

\affiliation[inst3]{organization=Computational Biology Lab, MBG-CSIC (Spanish National Research Council), 36143 Pontevedra, Spain}

\begin{abstract}
We consider the problem of parameter estimation in dynamic systems described by ordinary differential equations. A review of the existing literature emphasizes the need for deterministic global optimization methods due to the nonconvex nature of these problems. Recent works have focused on expanding the capabilities of specialized deterministic global optimization algorithms to handle more complex problems. Despite advancements, current deterministic methods are limited to problems with a maximum of around five state and five decision variables, prompting ongoing efforts to enhance their applicability to practical problems. 
Our study seeks to assess the effectiveness of state-of-the-art general-purpose global and local solvers in handling realistic-sized problems efficiently, and evaluating their capabilities to cope with the nonconvex nature of the underlying estimation problems.
\end{abstract}


\begin{highlights}
\item We show that current mathematical programming solvers can be effective for not-so-small problems.
\item Our approach provides insights on the challenges posed by local optimality.
\item We assess the impact of different discretization techniques and mathematical programming models.
\item Our study can be used as a benchmark to assess the performance of specialized algorithms.

\end{highlights}


\end{frontmatter}



Published in the journal Optimization and Engineering (2025)

\url{https://doi.org/10.1007/s11081-025-09978-9}

\section{Introduction}
\label{sec:intro}

Mathematical models play a crucial role in describing and analyzing real-world phenomena across engineering, natural sciences, and medicine. Here we consider the problem of parameter estimation in models described by nonlinear ordinary differential equations (ODEs) \citep{Schittkowski2002}. That is, given a parametric dynamic model and experimental data, we seek to identify the parameter values that minimize the mismatch between model predictions and data. These problems play a key role in many areas, including systems biology \citep{mendes1998non,Ashyraliyev2009,Banga2008,chou2009recent} and chemical kinetics \citep{Esposito2000,Singer2005}. 

From the mathematical point of view, parameter estimation in ODE-based models belong to the class of dynamic optimization problems. Existing solution approaches can be broadly categorized into direct and indirect methods \citep{biegler2010nonlinear}. Indirect methods concentrate on deriving a solution that adheres to the classical necessary conditions for optimality, manifested as a two-point boundary value problem \citep{BrysonHo69}. Direct methods apply discretization schemes to transform the original infinite-dimensional problem into a finite-dimensional nonlinear programming problem. Within direct methods, sequential approaches proceed by discretizing the control variables only, transforming the problem into an outer nonlinear programming (NLP) problem with an inner initial value problem embedded \citep{goh1988control,vassiliadis1994solution}. Simultaneous approaches \citep{biegler2007overview} discretize both control and state profiles, transforming the ODE system into algebraic equations for optimization. While the simultaneous approach offers advantages with respect to sequential approaches (e.g. easier handling of unstable dynamic systems and path constraints, and facilitating automatic differentiation), it also introduces challenges by potentially leading to large scale nonlinear programming problems that are difficult to solve.

Even though a number of local optimization methods have been developed to estimate parameters in dynamic models \citep{Schittkowski2002,Edsberg1995}, they typically converge to local solutions due to the nonconvex nature of these problems \citep{Esposito2000,Chen2010,Ljung2013}. 
A common strategy to surmount this drawback is to use a multi-start strategy with local methods \citep{Raue2013,frohlich2017scalable}, i.e., to repeat local optimizations from random initial points inside the parameter's bounds. Other studies have advocated for the use of global optimization approaches, i.e., methods with specialized numerical strategies seeking the globally optimal solution \citep{Esposito2000,Moles2003}.

During the last two decades, several stochastic and hybrid global optimization methods have been presented, especially in the area of systems biology \citep{Moles2003,Balsa-Canto2008-ck,Jia2012-rs,egea2010evolutionary,sun2011parameter}. Even though some of these methods have shown good performance even with large-scale problems \citep{Villaverde2015,Penas2017,Villaverde2018-iq}, an intrinsic limitation persists: their inability to guarantee global optimality. Consequently, when these techniques do not yield a satisfactory fit, it remains uncertain whether this discrepancy arises from the model's inadequacy to explain the data or from the method's inability to locate the global solution.

Given that the problem of parameter estimation in dynamic systems, and more generally, the problem of dynamic optimization, constitutes a nonconvex optimization problem, certifying global optimality requires deterministic global optimization methods. Pioneering research in this area commenced in the early 2000s \citep{Esposito2000,Papamichail2002,Papamichail2004,Singer2005,Chachuat2006,Lin2006,Lin2007,PerezGalvan2017}. In the area of systems biology, examples of relevant works include \cite{polisetty2006identification,panning2008deterministic,Miro2012,pitt2018critical}. Recent and comprehensive overviews of the methods and challenges involved in rigorously solving these dynamic optimization problems are given by \cite{wilhelm2019global} and \cite{Song2022}.

However, as noted by \cite{Song2022}, despite all these efforts, specialized deterministic methods for global dynamic optimization are currently limited to handling problems with a maximum of approximately five state variables and five decision variables. Consequently, there is an ongoing pursuit to broaden the applicability of deterministic global optimization to dynamic systems, making them suitable for solving problems of practical significance. For example, \citep{Sass2024} have recently presented a spatial branch-and-bound algorithm that exploits the structure of these problems when large datasets are considered.

Given the aforementioned limitations, the main goal of this paper is to study the potential of mathematical programming modeling and state-of-the-art solvers to handle a set of challenging parameter estimation problems in ODEs. Specifically, we focus on assessing the effectiveness and reliability of contemporary global and local optimization solvers. Our thorough computational experiments allow to assess to what extent global solvers can handle realistic-sized problems within acceptable computation times and, also, if local solvers can effectively tackle the nonconvex characteristics of the underlying parameter estimation problems without getting stuck at local optima.

Importantly, our numerical results show that state-of-the-art mathematical programming solvers can be effective for not-so-small problems in relatively short time (executions are limited to ten minutes): not only we find that both the local and global solvers are capable of consistently finding the optimal solutions, but we also show that the global solvers can certify global optimality. It is worth noting that these results were obtained with standard discretization schemes and out-of-the-box configurations of the solvers which, moreover, were not allowed to use any kind of parallelization capabilities. This hints at the potential of tackling even larger problems by relying on finely tuned solver configurations and high-performance computing. Last, but not least, we believe our study can serve as a benchmark to assess the potential of present and future specialized algorithms.

The rest of the paper is structured as follows. In Section~\ref{sec:framework} we present a detailed description of the framework for the analysis. In Section~\ref{sec:setup} we describe the specifics of the numerical study. Section~\ref{sec:performance} is devoted to a general overview of the computational results and, then, we devote Section~\ref{sec:local} to discuss the implications of the results for relevant aspects in parameter estimation such as local optimality, identifiability and flatness of the objective function. Finally, we conclude in Section~\ref{sec:conclusion}.

\section{Framework for the analysis}\label{sec:framework}
In this work we focus on parameter estimation in dynamic models that can be described by a nonlinear system of differential algebraic equations, DAEs. The system depends on some parameters that we want to estimate. To this aim, some experimental data from the process is usually available (observations of the dynamic states over time), so the goal is to obtain the value of the parameters that delivers the best possible fit of the experimental data. A general approach to describe the dynamics of the model is the following:
\begin{eqnarray*}
\frac{d\xb(t)}{dt}& = & f(t,\xb(t),\pb),\  t\in \left[t_0,t_f\right] \\
\yb(t)            & = & g(\xb(t),\pb),\  t\in \left[t_0,t_f\right] \\
\xb(t_0)          & = & \bar{\xb}_{0}, 
\end{eqnarray*}
where $\xb \in \mathbb{R}^{n_{s}}$ is the vector of state variables, $\bar{\xb}_{0}$ represents their initial conditions and $\pb \in \mathbb{R}^{n_{p}}$ is the vector of unknown parameters. Furthermore, $f$ is a nonlinear function describing the dynamics of the problem and $g$ is the observation function that gives the vector of observed states $\yb\in \mathbb{R}^{n_{y}}$ predicted by the model.

Let $(\tau_1,\bar{\yb}_{1})$,$\ldots$,$(\tau_n,\bar{\yb}_n)$ be the input experimental data, where $\bar{\yb}_{i}$ is the measured value of variable $\yb$ at some time $\tau_{i}\in \left[t_0,t_f\right]$. Thus, we want to solve the optimization problem of finding the value of the $\pb$ parameters that minimizes the error between the model predictions $\yb$ and the observable measurements $\bar{\yb}$. This leads to the formulation of the following Non Linear Programming (NLP) problem with DAEs:\footnote{Different error functions can be used to measure this error. We follow the standard in the field: the method of least squares.}
\begin{equation}
\begin{array}{rl}
\min        & \displaystyle \sum_{i=1}^{n}||\yb(\tau_{i})-\bar{\yb}_{i}||^{2}\\[18pt]
\text{s.t.} & \displaystyle \frac{d\xb(t)}{dt} = f(t,\xb(t),\pb), \ t\in \left[t_0,t_f\right]\\[10pt]
            & \yb(t)             = g(\xb(t),\pb), \ t\in \left[t_0,t_f\right]\\
            & \xb(t_0)           = \bar{\xb}_{0},  \\
						& \underline{\pb}\leq \pb \leq  \bar{\pb},
\end{array}
\label{eq:proborig}
\end{equation}
where $\underline{\pb}$ and  $\bar{\pb}$ represent known lower and upper bounds for the parameters.

\subsection{Direct transcription approach: \modeA \text{formulation}}
In order to solve this problem, a direct transcription approach is employed (see~\cite{betts2020} for details), which consists of discretizing problem \eqref{eq:proborig}. More precisely, the time domain is discretized into $M=\frac{t_{f}-t_{0}}{h}$ 
time steps, being $h$ the step size and $M$ the size of the mesh. Thus, a finite dimensional approximation of the original problem is created by discretizing all the variables and constraints at discrete time points $t_{m}$ with $m\in \{0,\ldots,M\}$, where $t_{M}=t_{f}$. Denote by
\[
\pmb{\xi} = (\xb_0,\yb_{0},\xb_1,\yb_{1},\ldots,\xb_M,\yb_{M},\pb)
\] 
the optimization variables arising from the discretization, where $\xb_m = \xb(t_{m})$ and $\yb_m = \yb(t_{m})$ for all $m\in \{0,\ldots,M\}$. Further, the ODEs are approximated using a differential numerical discretization scheme that effectively reformulates them into a system of algebraic equations. Therefore, the original dynamic problem~\eqref{eq:proborig} can be reformulated as the following standard NLP problem:
\begin{equation}
\begin{array}{rlll}
\multicolumn{4}{c}{\modeA \text{ formulation}}\\
\min & \displaystyle \sum_{i=1}^{n}||\yb_{m(i)}-\bar{\yb}_{i}||^{2} & & \\[15pt]
\text{s.t.} & \displaystyle H(\pmb{\xi})= \pmb{0} \\
            & \yb_{m} = g(\xb_{m},\pb),\ \forall m\in \{0,\ldots,M\}\\
            & \xb_{0} = \bar{\xb}_{0}  \\
						& \underline{\pb}\leq \pb \leq  \bar{\pb},
\end{array}
\label{eq:prob_discr1}
\end{equation}
where $H$ is the system of algebraic equations obtained after applying the chosen discretization scheme and $m(i)$ gives the correspondence between the times at which the observable measures are taken, $\tau_{i}$, and the discretization points $t_{m}$, i.e.,
for each $i\in \{1\ldots,n\}$, $m(i)$ gives the value such that $t_{m(i)}=\tau_{i}$.\footnote{Thus, we are implicitly assuming that the discrete time points of the observable measures $\tau_{i}$ are a subset of the discretization points $t_{m}$. If this were not the case, an interpolation approximation could be applied to recover this property.} Naturally, the larger M is, the better~\eqref{eq:prob_discr1} approximates problem~\eqref{eq:proborig}. For our numerical study we implemented five classic  discretization schemes (see~\cite{burden2001numerical}). For the sake of completeness, we present below the resulting approximation of the ODEs:
\begin{itemize}
\item The \euler method. For each $m\in \{0,\ldots,M-1\}$ we have
\begin{equation}
\xb_{m+1} = \xb_{m}+h f\left(t_{m},\xb_{m},\pb\right). 
\label{eq:euler_scheme}
\end{equation}

\item The \trapecio method. For each $m\in \{0,\ldots,M-1\}$ we have
\begin{equation}
\xb_{m+1} = \xb_{m}+\dfrac{h}{2} \big(f\left(t_{m},\xb_{m},\pb\right)+f(t_{m+1},\xb_{m+1},\pb)\big).
\label{eq:trapezoid_scheme}
\end{equation}

\item The \adams method with step 3. For each $m\in \{0,\ldots,M-3\}$ we have\footnote{Initial conditions for $\xb_{1}$ and $\xb_{2}$ are taken from the \trapecio method.}
\begin{equation}
\begin{split}
\xb_{m+3} = \xb_{m+2} + \dfrac{h}{24}\Big(& 9 f(t_{m+3},\xb_{m+3},\pb)+19f(t_{m+2},\xb_{m+2},\pb) \\
                                          & -5 f(t_{m+1},\xb_{m+1},\pb)+f(t_{m},\xb_{m},\pb)\Big).
\label{eq:adams_scheme}
\end{split}
\end{equation}

\item The \simpson method. For each $m\in \{1,\ldots,M-1\}$ we have\footnote{Initial conditions for $\xb_{1}$ are taken from the \trapecio method.}
\begin{equation}
\xb_{m+1} = \xb_{m-1}+\dfrac{h}{3} \big(f\left(t_{m+1},\xb_{m+1},\pb\right)+4f(t_{m},\xb_{m},\pb)+f(t_{m-1},\xb_{m-1},\pb)\big).
\label{eq:simpson_scheme}
\end{equation}

\item The \runge method. For each $m\in \{1,\ldots,M-1\}$ we have
\begin{equation}
\xb_{m+1} = \xb_{m} + \dfrac{h}{6}(\kb^{1}_{m}+2\kb^{2}_{m}+2\kb^{3}_{m}+\kb^{4}_{m}),\ \text{where}\\
\label{eq:runge_scheme}
\end{equation}
\begin{eqnarray*}
\kb^{1}_{m}   & = & f(t_{m},\xb_{m},\pb) \\
\kb^{2}_{m}   & = & f(t_{m}+\dfrac{h}{2},\xb_{m}+\dfrac{h}{2}\cdot\kb_{1},\pb) \\
\kb^{3}_{m}   & = & f(t_{m}+\dfrac{h}{2},\xb_{m}+\dfrac{h}{2}\cdot\kb_{2},\pb) \\
\kb^{4}_{m}   & = & f(t_{m+1},\xb_{m}+h\cdot\kb_{3},\pb). \\
\end{eqnarray*}
\end{itemize}

\subsection{Variations of the mathematical programming model}

Given the difficulty of the optimization problems to be solved, we also consider two standard variations of the \modeA \text{formulation}, with the goal of assessing the impact of these alternative formulations in the performance of the state-of-the-art solvers. 

First, we define a relaxation of \modeA by introducing a fixed feasibility tolerance in the algebraic constraints. Specifically, given $\epsilon>0$, the \modeB formulation is given by
\begin{equation}
\begin{array}{rl}
\multicolumn{2}{c}{\modeB \text{ formulation}}\\
\min & \displaystyle \sum_{i=1}^{n}||\yb_{m(i)}-\bar{\yb}_{i}||^{2} \\[15pt]
\text{s.t.} & \displaystyle -\epsilon\leq H(\pmb{\xi})\leq \epsilon \\
            & \yb_{m} = g(\xb_{m},\pb),\ \forall m\in \{0,\ldots,M\}\\
            & \xb_{0} = \bar{\xb}_{0}  \\
						& \underline{\pb}\leq \pb \leq  \bar{\pb}. \\
\end{array}
\label{eq:prob_discr3}
\end{equation}
Note that feasibility tolerances can usually be controlled with solver-specific options. Yet, \modeB is solver independent, which enables a direct comparison of the impact of these tolerances across solvers without having to deal with the particularities of their respective configuration options.

A second variation of \modeA \text{formulation} is to treat all the algebraic constraints of the model as soft constraints. This is achieved by adding, to each algebraic constraint, a couple of slack variables, which are then heavily penalized in the objective function. In this case, we obtain \modeC formulation given by
\begin{equation}
\begin{array}{rl}
\multicolumn{2}{c}{\modeC \text{ formulation}}\\
\min & \displaystyle \sum_{i=1}^{n}||\yb_{m(i)}-\bar{\yb}_{i}||^{2} + P\cdot \sum_{j=1}^{n_1}s_{j} \\[18pt]
\text{s.t.} & \displaystyle -\pmb{s}\leq H(\pmb{\xi})\leq \pmb{s} \\
            & \yb_{m} = g(\xb_{m},\pb),\  \forall m\in \{0,\ldots,M\}\\
            & \xb_{0} = \bar{\xb}_{0}  \\
						& \underline{\pb}\leq \pb \leq  \bar{\pb} \\
						& \pmb{s}\geq \pmb{0},
\end{array}
\label{eq:prob_discr2}
\end{equation}
where $\pmb{s}\in \mathbb{R}^{n_1}$ are the slack variables added to the algebraic equations of the ODEs system. Further, $P \geq 0$ is the penalization in the objective function for using the slack variables, i.e., the penalization imposed on the violations of the constraints.




\section{Setting up the numerical study}\label{sec:setup}
The proposed framework has been implemented in AMPL~\citep{ampl}, using a PERL script \citep{perl} to automatically generate, for each problem instance, the AMPL model files for the different NLP formulations and the different discretization schemes. All numerical experiments have been performed on the supercomputer Finisterrae III, provided by Galicia Supercomputing Centre (CESGA), using computational nodes powered with two 32-core Intel Xeon Ice Lake 8352Y CPUs with 256GB of RAM and 1TB SSD.\footnote{It is worth noting that NEOS Server \citep{NEOS:1998} was also intensively used during the initial stages of the project to test and validate different model formulations.}

\subsection{Pool of parameter estimation problems}

We have applied our framework to a series of well known parameter estimation problems from the literature, which we outline below (refer to~\ref{app:description} for a deeper description):
	\begin{itemize}
		\item \alphapinene: it represents the mechanism for thermal isomerization of $\alpha$-pinene \citep{alphapinene}. It has 5 parameters and 5 states (fully observed).
		\item \BBG (Biomass batch growth): it describes the microbial growth in a stirred fed-batch bioreactor \citep{BBG}. It has 4 parameters and 2 states (fully observed).
		\item \FHN: the Fitzhugh-Nagumo model \citep{FHN1,FHN2}. It has 3 parameters and 2 states with only one observed state (partially observed). 
		\item \harmonic: it represents a harmonic oscillator (pendulum), as considered in \cite{Go2023}. It has 2 parameters and 2 states (fully observed).
		\item \lv: it is a two species Lotka-Volterra predatory-prey model, with 3 parameters and 2 states, as considered by \cite{Go2023}. Two versions of this problem are considered depending on the number of observed variables considered when solving it. We denote by \lvF the case with 2 observed states (fully observed) and by \lvP the problem with only one observed state (partially observed).
		\item \daisy: it is based on the examples used in the analysis of DAISY identifiability software \citep{daisy}, and it represents a 3-compartment model. It has 5 parameters and 3 states. Two versions of this problem are considered: \daisyF (fully observed) and \daisyP (with partially observable data: 2 observed states). 
		\item \hiv: it is is based on a model of HIV infection dynamics coming from \citep{hiv}. It has 10 parameters and 5 states. Two versions of this problem are considered: \hivF (fully observed) and \hivP (with partially observable data: 3 observed states and observed data for the sum of the other 2 states).
		\item \crauste: it is based on the system described in \citep{crauste}, and it has 16 parameters and 5 state variables. Two versions of this problem are considered: \crausteF (fully observed) and \crausteP (with partially observable data: 3 observed states and observed data for the sum of the other 2 states). 
	\end{itemize}

The procedure to generate the noiseless measurement data for the observed variables of each problem is as follows.\footnote{The only exception is \alphapinene which, given its simplicity, was fed with real (and noisy) measurements.} We take nominal values for the parameters of the different systems of ODES from past literature and solve the ODE system for those values of the parameters with the R package \texttt{pracma} \citep{pracma} (using Runge-Kutta). Then, from the results given by \texttt{pracma}, we select a sample of points of the state variables. As we describe in \ref{app:description}, for most of the problems these points are uniformly distributed on a given time of the form $[0,T]$, with $T$ ranging from $1$ to $20$. 
 Table~\ref{tab:observed_states} shows, for each problem, the number of i)~state variables, ii)~observed variables, iii)~parameters\footnote{When appropriate, we indicate in parenthesis the number of initial conditions to be estimated.} and iv)~sample time points.
\begin{table}[!htbp]
\centering
\begin{tabular}{l|ccccc}
\toprule
        & State Var. & Observed Var. & Parameters & Observed times \\[-0.2cm]
Problem             & $n_{s}$    & $n_{y}$       &  $n_{p}$   &    $n$               \\[-0.05cm]
\midrule
\alphapinene & 5 & 5  & 5  & 8   \\
\BBG         & 2 & 2  & 4  & 7  \\
\FHN         & 2 & 1  & 3  & 6  \\
\harmonic    & 2 & 2  & 2 (2)  & 10  \\
\lvF         & 2 & 2  & 3 (2) & 20  \\
\lvP         & 2 & 1  & 3 (2) & 20  \\
\daisyF      & 3 & 3  &  5 (3) & 20  \\
\daisyP      & 3 & 2  &  5 (3) & 20  \\
\hivF        & 5 & 5  & 10 (5) & 20  \\
\hivP        & 5 & 4  & 10 (5) & 20  \\
\crausteF    & 5 & 5  & 13 (5) & 20  \\
\crausteP    & 5 & 4  & 13 (5) & 20  \\
\bottomrule
\end{tabular}
\caption{Number of state variables, observed variables, parameters and time points for measurements.}
\label{tab:observed_states}
\end{table}

It is worth mentioning that, in the literature, most of these problems are solved under fixed initial conditions. However, for many of our problems (such as \harmonic and the four problems for which we have both fully and partially versions), we assume these conditions are unknown and need to be estimated as well. This increases the difficulty for solving them.

\subsection{Configurations of the elements of the computational framework}

Each parameter estimation problems was tackled with different configurations of the elements involved in our framework. We now detail all such configurations: 

\begin{description}
	\item[I. Mathematical formulations:] we run the three mathematical programming formulations described in the previous section: \modeA, \modeB and \modeC. The configurable parameters for  \modeB and \modeC are chosen as follows:
	\begin{itemize}
			\item For \modeB, two values for the feasibility tolerance $\epsilon$ are considered: $10^{-4}$ and $10^{-6}$.\footnote{In order to avoid some numerical issues observed in some preliminary runs, slightly different options are considered for \lvP ($10^{-5}$, $10^{-6}$ and $10^{-8}$) and \crauste ($10^{-5}$, $10^{-7}$ and $10^{-9}$).}
		\item For \modeC, two values for $P$ are considered: $10^3$ and $10^5$.
	\end{itemize}
	
	\item[II. Discretization numerical schemes:] the 5 discretization schemes introduced in the previous section are considered: \euler, \trapecio, \adams, \simpson and \runge. Further, each problem is run with two different mesh sizes ($M$), as specified in Table~\ref{tab:mesh_sizes}.\footnote{For \lvP problem, on top of the mesh sizes of 100 and 1000, a bigger mesh of size 10000 was also tested.}

\begin{table}[h!]
\centering
{\footnotesize
\begin{tabular}{cccccccc}
\toprule
\lv & \crauste & \FHN & \harmonic & \alphapinene & \BBG & \daisy & \hiv\\
\midrule
100          & 100  & 200  & 230  & 1230 & 120  & 100  & 100 \\
1000         & 1000 & 2000 & 2300 & 3690 & 1200 & 1000 & 1000 \\
\midrule
\bottomrule
\end{tabular}
}
\caption{Mesh sizes for each problem.}
\label{tab:mesh_sizes}
\end{table}

	\item[III. Mathematical programming solvers:] given the interest in analyzing whether the solutions obtained with our framework are global optima, it is pertinent to distinguish between local and global solvers. Specifically, three global solvers and five local solvers are employed in the numerical experiments:
\begin{description}
	\item[Global Solvers: ] \Couenne 0.5.7 \citep{couenne}, \BARON 21.1.13 (2021.01.13) \citep{Baron} and \Octeract Engine v4.4.1\citep{octeract}.
	\item[Local Solvers: ] \BONMIN 1.8.7 \citep{bonmin}, \Knitro 13.2.0 \citep{knitro}, \Ipopt 3.12.13 \citep{ipopt}, \CONOPT 3.17A \citep{conopt} and \SNOPT 7.5-1.2 \citep{snopt}.
\end{description}
\end{description}

Each solver is allowed to use only one thread and a time limit of 10 minutes is imposed. All other options of the solvers are set to their default values. All the combinations of the above configurations result in a total of 5280 executions. In particular, we tested at least 400 configurations on each parameter estimation problem: 5 mathematical programming formulations, 5 discretization schemes with 2 mesh sizes each and 8 state-of-the-art solvers.

\subsection{Performance metrics}

In order to evaluate the quality of the solutions given by our framework, two metrics are employed. On the one hand, it is important to check if the reference solution is found. Given a problem, let $\pb^{\text{\sc ref}}$ be the reference solution of the parameters and let $\pb$ be the estimation obtained with our approach. Then, the \textit{maximum of the relative error} is computed as follows:
\[
\text{MaxRE} = \max_{i=1,\ldots,n_{p}} \left|\frac{p_{i}-p^{\text{\sc ref}}_{i}}{p^{\text{\sc ref}}_{i}}\right|.
\]
Note that, the smaller MaxRE is, the closer of the solution to the reference one. In the numerical study, the criterion to consider that the solution is acceptably close to the reference solution is that $\text{MaxRE}<0.1$. The error showed in the figures of the computational results corresponds to MaxRE.

On the other hand, since some of the problems can present a potential lack of identifiability, it could happen that the parameter estimation obtained is not really close to the reference solution, but the adjustment with respect to the observed state variables is good. In order to analyze this, we use the following \textit{normalized mean squared error} measure:
\[
\text{NRMSE} = \frac{\sqrt{\frac{\sum_{i=1}^{n}||\yb_{m(i)}-\bar{\yb}_{i}||^{2}}{n_{y}\cdot n}}}{\bar{y}_{max} - \bar{y}_{min}},
\]
where $\bar{y}_{max}=\max\{\max\{\bar{\yb}_{1}\},\ldots,\max\{\bar{\yb}_{n}\}\}$ and $\bar{y}_{min}=\min\{\min\{\bar{\yb}_{1}\},\ldots,\min\{\bar{\yb}_{n}\}\}$.

Recall that, with the exception of \alphapinene, all problems are fed with noiseless synthetic data, whose generation ensures that the NRMSE at the reference solution should be zero or very close to it (depending on the precision of the discretization).

%
%
%
%
%
%
%
%

\section{Numerical results: general overview}\label{sec:performance}

This section focuses on the study of the overall performance of the mathematical programming solvers and the different modeling choices. The results are illustrated with a series of tables, each of them containing the following columns:
\begin{description}
\item[\solvedS.] Proportion of runs in which the solver returned ``solved'' as its final status.
\item[\foundref.] Proportion of runs in which the solver actually found the reference solution ($\text{MaxRE}\leq 0.1$). 
\item[\nearref.] Proportion of runs in which the solver did not find the reference solution but was close to it ($0.1<\text{MaxRE}\leq 0.5$).
\item[\alternat.] Proportion of runs in which the solution was not close to the reference solution ($\text{MaxRE}>0.5$) but, still, $\text{NRMSE}<0.0001$ (potential lack of identifiability).
\item[\timeF.] Smallest solve time, in seconds, among the runs in which the solver found the reference solution.
\item[\success.] Proportion of runs in which the solver execution was successful, in the sense of not experimenting numerical issues, having abrupt terminations or failing to return a feasible solution.
\end{description}

A supplementary Online Appendix (\url{https://bit.ly/ParamEstimODEs}) contains a more comprehensive set of results, including more complete tables and additional figures.

\subsection{Results by problem}

Table~\ref{tab:byproblem} shows the results disaggregated by problem and sorted by column \foundref. A first finding that stands out is that, for each and every problem, the reference solution was found at least once. Moreover, in more than half of the problems, the reference solution was found by more than one third of the configurations, which shows the robustness of the mathematical programming formulations, since no fine tuning of the underlying configurations is needed to found the reference solutions and even certifying its global optimality. Column \timeF shows that running times are also promising since, with the exception of \alphapinene, for each problem there is at least one configuration that has found the reference solution in less than one second and, very often, in less that a tenth of a second. These running times are particularly encouraging, given that they include problems such as \crauste and \hiv, which have more than 10 parameters to estimate. Yet, having quick running times should not come as a surprise, since many of the configurations involve relatively coarse discretizations. What is truly interesting is that, despite this coarseness, the reference solution is often found.

\begin{table}[!htbp]
\centering
\begin{tabular}{lcccccc}
\toprule
Problem & \solvedS & \foundref & \nearref & \alternat & \timeF & \success\\
\midrule
\daisyF & 0.735 & 0.885 & 0.022 & 0.002 & 0.029 & 0.965 \\
\lvF & 0.743 & 0.835 & 0.022 & 0\phantom{.000} & 0.026 & 0.960 \\
\harmonic & 0.715 & 0.733 & 0.002 & 0\phantom{.000} & 0.032 & 0.927 \\
\hivF & 0.820 & 0.688 & 0.122 & 0.082 & 0.033 & 0.975 \\
\BBG & 0.595 & 0.438 & 0.052 & 0\phantom{.000} & 0.084 & 0.802 \\
\alphapinene & 0.250 & 0.365 & 0.005 & 0\phantom{.000} & 1.786 & 0.667 \\
\daisyP & 0.665 & 0.343 & 0.092 & 0.030 & 0.028 & 0.958 \\
\crausteF & 0.796 & 0.298 & 0.173 & 0.362 & 0.054 & 0.958 \\
\hivP & 0.792 & 0.195 & 0.082 & 0.490 & 0.103 & 0.965 \\
\crausteP & 0.781 & 0.185 & 0.144 & 0.375 & 0.069 & 0.969 \\
\FHN & 0.415 & 0.117 & 0.007 & 0.030 & 0.514 & 0.802 \\
\lvP & 0.694 & 0.024 & 0.376 & 0.393 & 0.033 & 0.867 \\
\bottomrule
\end{tabular}
\caption{Summary results by problem.}
\label{tab:byproblem}
\end{table}

The results in Table~\ref{tab:byproblem} also hint at two other aspects particularly relevant in parameter estimation problems and that we explore more deeply in Section~\ref{sec:local}. First, in some problems the \solvedS value is higher than the value in \foundref, suggesting that there are multiple local optima in the problem at hand. At the same time, a value in \foundref higher than the one in \solvedS indicates that global solvers have found the optimal solution, but failed to certify it. Second, values different from zero in \alternat come from the existence of high quality solutions different from the reference one, which may be the result of a lack of identifiability in the problem at hand; indeed, with the exception of \crausteF, the largest values in \alternat are in partially specified problems.

\subsection{Results by solver}

We move now to Table~\ref{tab:bysolver}, which shows the results by solver, again sorted by column \foundref. The first aspect that stands out is that two of the global optimization solvers involved in the study come out on top, which is particularly valuable given that they are often capable of providing global optimality certificates for the discretization at hand. As expected, local solvers are faster, although \BARON is very competitive.


\begin{table}[!htbp]
\centering
\begin{tabular}{lcccccc}
\toprule
Solver & \solvedS & \foundref & \nearref & \alternat & \timeF & \success \\
\midrule
\BARON & 0.667 & 0.518 & 0.135 & 0.168 & 0.076 & 0.992 \\
\Couenne & 0.691 & 0.482 & 0.130 & 0.089 & 0.138 & 0.841 \\
\Knitro & 0.961 & 0.430 & 0.130 & 0.156 & 0.026 & 0.994 \\
\Octeract & 0.524 & 0.408 & 0.115 & 0.182 & 1.286 & 0.867 \\
\BONMIN & 0.897 & 0.403 & 0.111 & 0.117 & 0.038 & 0.897 \\
\Ipopt & 0.886 & 0.395 & 0.124 & 0.211 & 0.038 & 0.933 \\
\CONOPT & 0\phantom{.000} & 0.306 & 0.098 & 0.215 & 0.044 & 0.921 \\
\SNOPT & 0.752 & 0.221 & 0.045 & 0.212 & 0.040 & 0.764 \\
\bottomrule
\end{tabular}
\caption{Summary results by solver.}
\label{tab:bysolver}
\end{table}

\subsection{Results by discretization scheme}

In Table~\ref{tab:bydiscretization} we present the results by discretization scheme and, although there are no major differences, \trapecio seems to perform slightly better than the rest.


\begin{table}[!htbp]
\centering
\begin{tabular}{lcccccc}
\toprule
Scheme & \solvedS & \foundref & \nearref & \alternat & \timeF & \success \\
\midrule
\trapecio & 0.723 & 0.462 & 0.104 & 0.177 & 0.026 & 0.951 \\
\adams & 0.693 & 0.399 & 0.104 & 0.160 & 0.029 & 0.891 \\
\simpson & 0.691 & 0.396 & 0.126 & 0.149 & 0.028 & 0.914 \\
\runge & 0.604 & 0.365 & 0.097 & 0.166 & 0.093 & 0.813 \\
\euler & 0.650 & 0.356 & 0.125 & 0.192 & 0.029 & 0.937 \\
\bottomrule
\end{tabular}
\caption{Summary results by discretization scheme.}
\label{tab:bydiscretization}
\end{table}

\subsection{Results by mathematical programming formulation}

Finally, Table~\ref{tab:byformulation} studies the impact of the variations of the mathematical programming models to be solved. Again, the performance seems to be similar across the three approaches and, if anything, the use of \modeC may be slightly worse.


\begin{table}[!htbp]
\centering
\begin{tabular}{lcccccc}
\toprule
Form. & \solvedS & \foundref & \nearref & \alternat & \timeF & \success \\
\midrule
\modeA & 0.643 & 0.441 & 0.081 & 0.146 & 0.032 & 0.872 \\
\modeB & 0.716 & 0.410 & 0.141 & 0.203 & 0.028 & 0.911 \\
\modeC & 0.637 & 0.356 & 0.092 & 0.141 & 0.026 & 0.904 \\
\bottomrule
\end{tabular}
\caption{Summary results by mathematical programming formulation.}
\label{tab:byformulation}
\end{table}

\subsection{Additional disaggregations of the results}

\ref{app:performanceplus} studies different disaggregations of the above results, with the goal of getting further insights, not only of the individual elements involved in the different configurations, but also of potential interactions. They show that \daisyF, \harmonic and \lvF are the easiest problems to solve and, moreover, that global optimization solvers seem to be better at finding the reference solution. Interestingly, Table~\ref{tab:byproblemformulation} shows that the choice of the mathematical programming formulation can have a significant impact. \modeA performs extremely well in \hivF and \daisyF problems, finding the reference solutions in more than 90\% of the configurations, whereas the next best performing formulation for \hivF is under 65\%. On the other hand, \modeB is clearly the best formulation for \harmonic, whereas \modeC comes out narrowly on top for \lvF. Regarding the disaggregation by discretization scheme and mathematical programming formulation, we see in Table~\ref{tab:discretizationformulation} that the configurations with \trapecio tend to come out on top, while \euler and \runge are at the bottom. Moreover, for all discretization schemes, the configurations with \modeA seem to dominate those with \modeB, which themselves seem to dominate those configurations with \modeC.

\subsection{Running times by problem}

Section~3 in the Online Appendix also provides detailed information about running times. In particular, it presents, for each problem, the top 10 configurations in terms of running times among those that have found the reference solution (top 10 among all solvers and also top 10 among global solvers). 

The results show that \Knitro is the fastest solver in most of the problems. Some exceptions are \FHN, where \Ipopt, \CONOPT and \BONMIN are also quite fast, and \lvP, \daisyP, \crausteP and \crausteF, where again \Ipopt and \BONMIN are relatively close to \Knitro. Regarding running times, local solvers are very fast, and \alphapinene is the only problem not solved by any configuration in less than a second (the top 10 fastest configurations go from 1.8 seconds to 8.9 seconds). Moreover, more than half of the problems are solved in less than 0.1 seconds by the fastest configuration.

If we restrict attention to global optimization solvers, the results are also quite promising. Setting again aside \alphapinene, for which no global solver was able to certify global optimality for the discretization at hand after 10 minutes, all other problems where solved to global optimality, by at least one configuration, in less than 10 seconds. \BARON is the fastest global optimization solver in most problems and, for some of them such as \BBG and \hivP, it is even competitive with the fastest local solver, \Knitro. On the other hand, \Octeract is the fastest in \lvF and \Couenne is competitive with \BARON in \hivP, \hivF and \crausteP.

It may seem surprising that \alphapinene, despite coming from a relatively simple ODE system, has turned out to be relatively challenging, particularly for the certification of global optimality. Yet, recall that, differently from the other problems in the study, \alphapinene was fed with real (and noisy) measurements instead of synthetic ones. We have observed in some additional experiments with this problem that, with noiseless synthetic data, global certificates are consistently obtained within the time limits. The use of real measurements in this particular problem may have led to some practical lack of identifiability driven by the flatness of the objective function. These issues are discussed in the following section and, for the particular case of \alphapinene, they can be assessed in the figures in~\ref{sec:alpha}.

\section{Numerical results: local optimality, flatness, and identifiability}
\label{sec:local}
The approach followed in this paper, in which each problem is tackled with a wide range of configurations, including different families of local and global optimization algorithms, can be a useful tool to pinpoint structural characteristics of the underlying parameter estimation problems. In this section we illustrate the potential of our approach to accomplish this. We do so through a series of examples in which one can recognize patterns for multiplicity of local optima, flatness of the objective function, and lack of identifiability (i.e., multiplicity of global optima). \ref{app:plots} contains detailed results for all problems, similar to the selected ones that we discuss below.

\subsection{Multiplicity of local optima}

Overall, the existence of multiple local optima with significantly different objective functions (sum of squared errors in the resulting fit) does not seem to be an issue in the problems under study. In particular, the solutions provided by local optimizers are predominantly as close to the reference solution as the ones provided (and certified) by the global optimization solvers. 

\begin{figure}[!htbp]
\centering
	\begin{subfigure}{0.495\textwidth}
		\begin{tikzpicture}
			\node(s1) at (0,0){\includegraphics[width=0.9\textwidth]{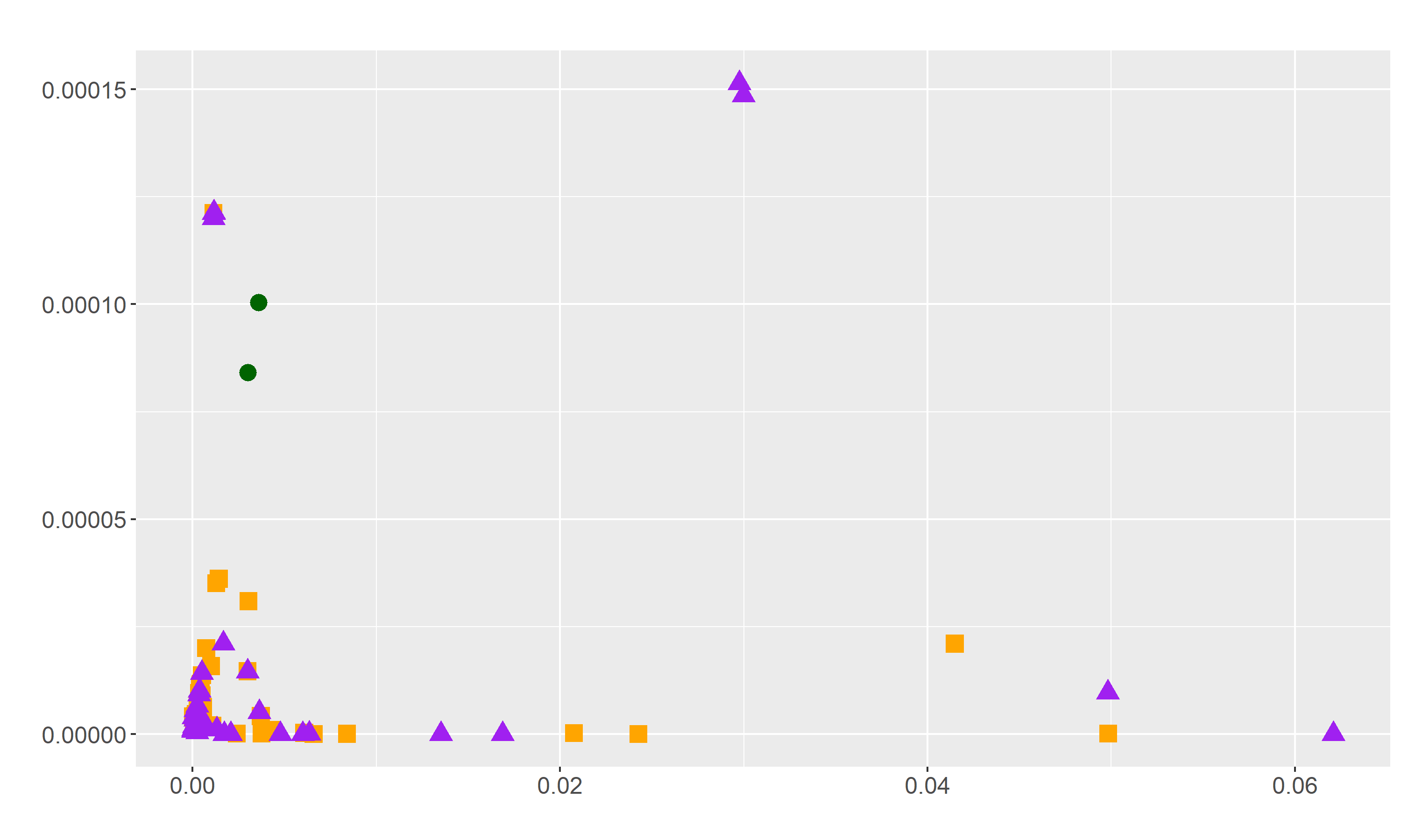}};
			\node(s2) at (2,0){\includegraphics[width=1cm]{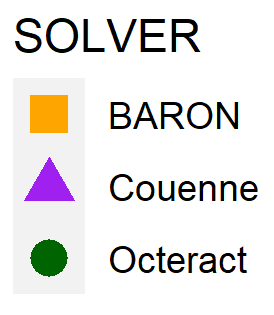}};
			\node[rotate=90](y) at (-3.45,0){{\color{darkgray}\tiny{$\log(\text{NRMSE}+1)$}}};
			\node(x) at (0,-2){{\color{darkgray}\tiny{MaxRE}}};
		\end{tikzpicture}
		\caption{Global solvers: Distance vs error when ``solved''.}
	\end{subfigure}		
	\begin{subfigure}{0.495\textwidth}
		\begin{tikzpicture}
			\node(s1) at (0,0){\includegraphics[width=0.9\textwidth]{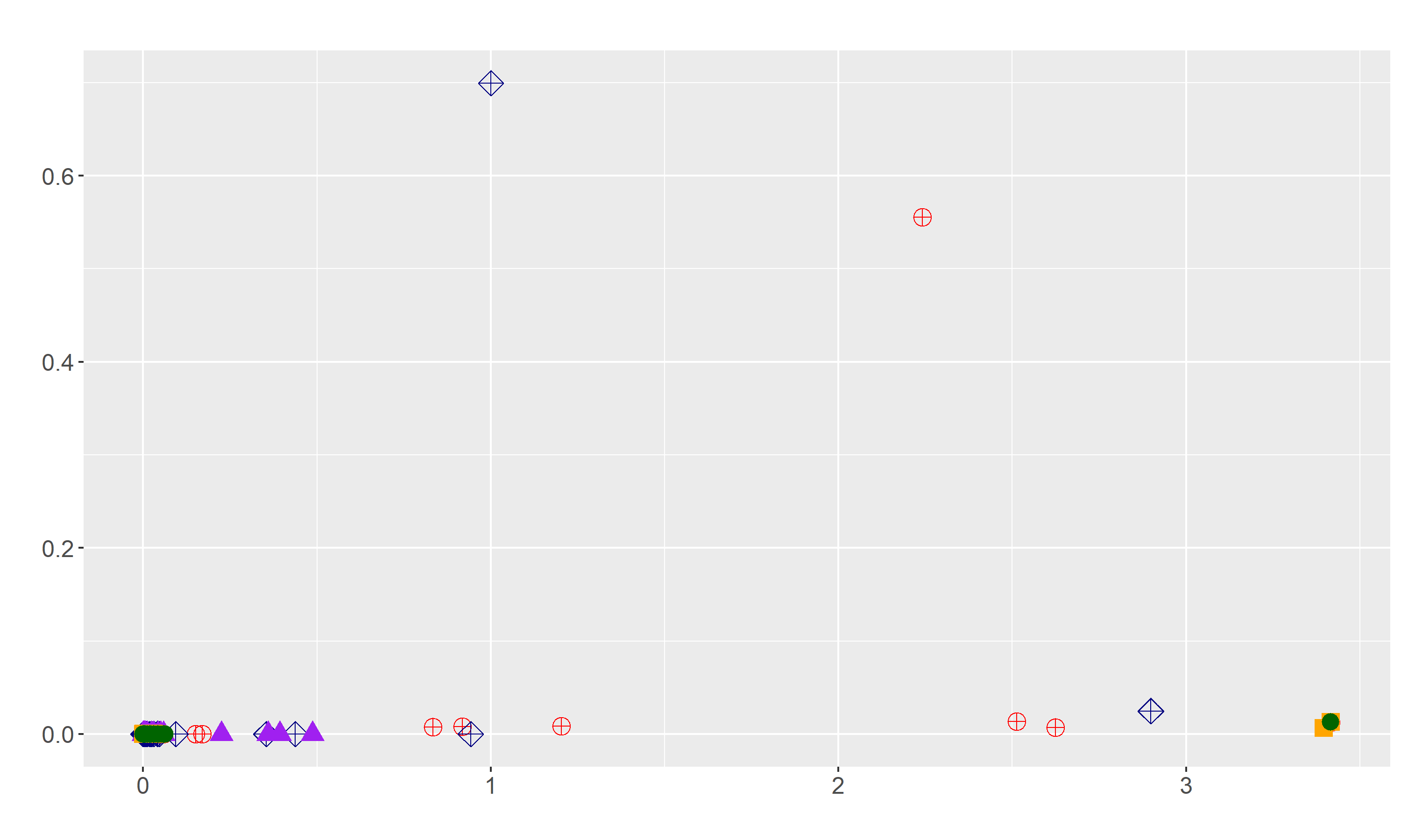}};
			\node(s2) at (2,0){\includegraphics[width=1cm]{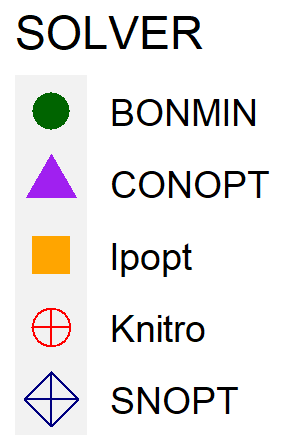}};
			\node[rotate=90](y) at (-3.45,0){{\color{darkgray}\tiny{$\log(\text{NRMSE}+1)$}}};
			\node(x) at (0,-2){{\color{darkgray}\tiny{MaxRE}}};
		\end{tikzpicture}
		\caption{Local solvers: Distance to ref. solution vs error.}
	\end{subfigure}
\caption{Studying local optimality in \daisyF.}
\label{fig:local}
\end{figure}

Figure~\ref{fig:local} focuses on problem \daisyF. The two plots represent the distance to the reference solution on the $x$-axis and the error on the $y$-axis, for the different configurations run for this problem. On the left, we represent the solutions provided by global solvers in configurations where they were able to certify global optimality on the discretized problem at hand. On the right we represent all the solutions provided by local solvers (just ensuring some local optimality condition).

We can see that, although the solutions provided by the global solvers tend to be closer to the reference solution, the ones provided by the local solvers are not far away and the associated errors are typically very small as well. In this case, the fact that the solutions of the local solvers are slightly more spread than the ones provided by the global solvers may be an indication of having a relatively flat objective function near the global optimum (which we discuss below), but there is no problem of getting stuck at local optima of very poor quality in terms of the error in the resulting fit.

\subsection{Flatness of the objective function}

\begin{figure}[!htbp]
\centering
\begin{subfigure}{0.495\textwidth}
		\centering
		\begin{tikzpicture}
			\node(s1) at (0,0){\includegraphics[width=0.9\textwidth]{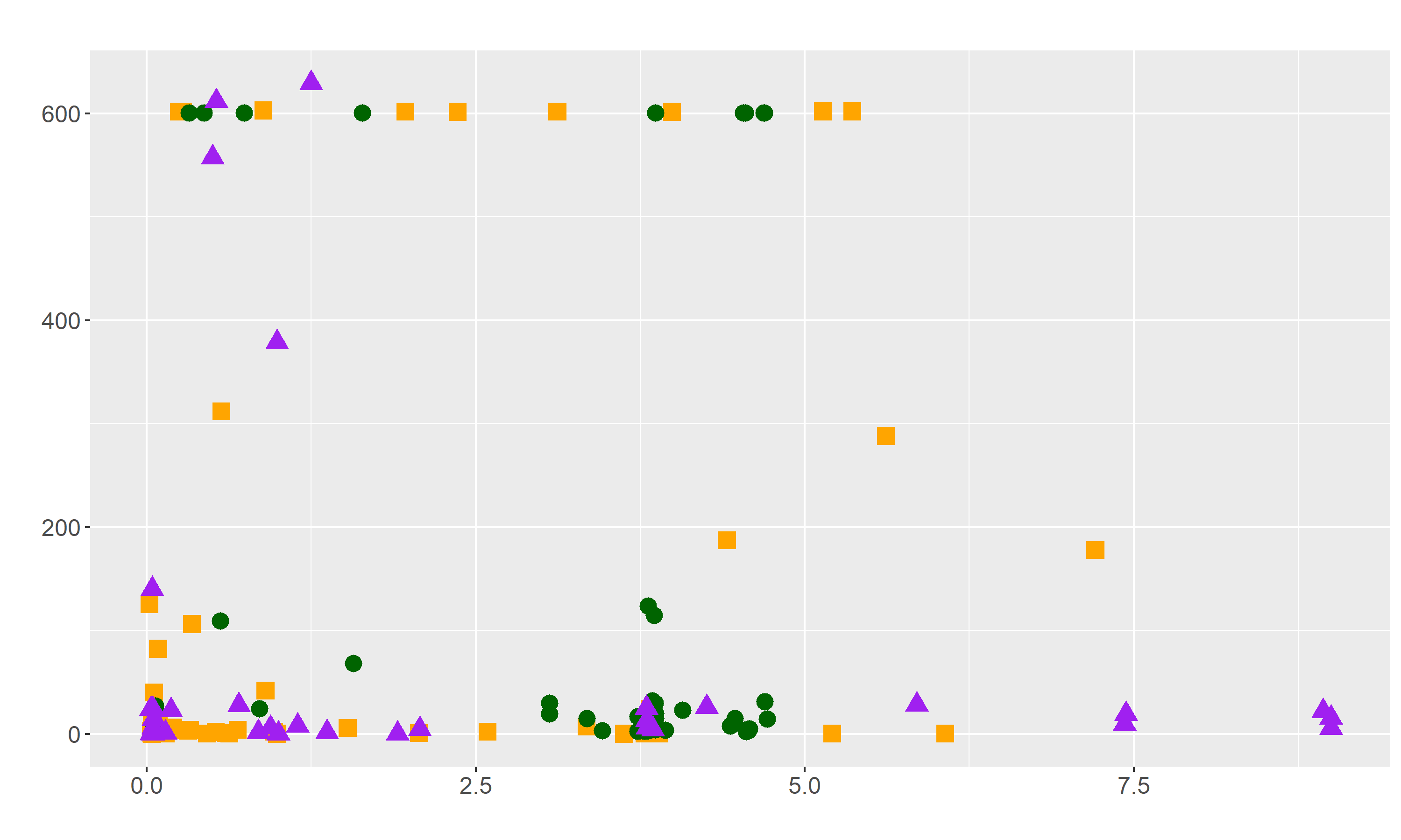}};
			\node[rotate=90](y) at (-3.45,0){{\color{darkgray}\tiny{time}}};
			\node(x) at (0,-2){{\color{darkgray}\tiny{MaxRE}}};
		\end{tikzpicture}
		\caption{Distance to ref. solution vs time.}
\end{subfigure}
\begin{subfigure}{0.495\textwidth}
		\centering
		\begin{tikzpicture}
			\node(s1) at (0,0){\includegraphics[width=0.9\textwidth]{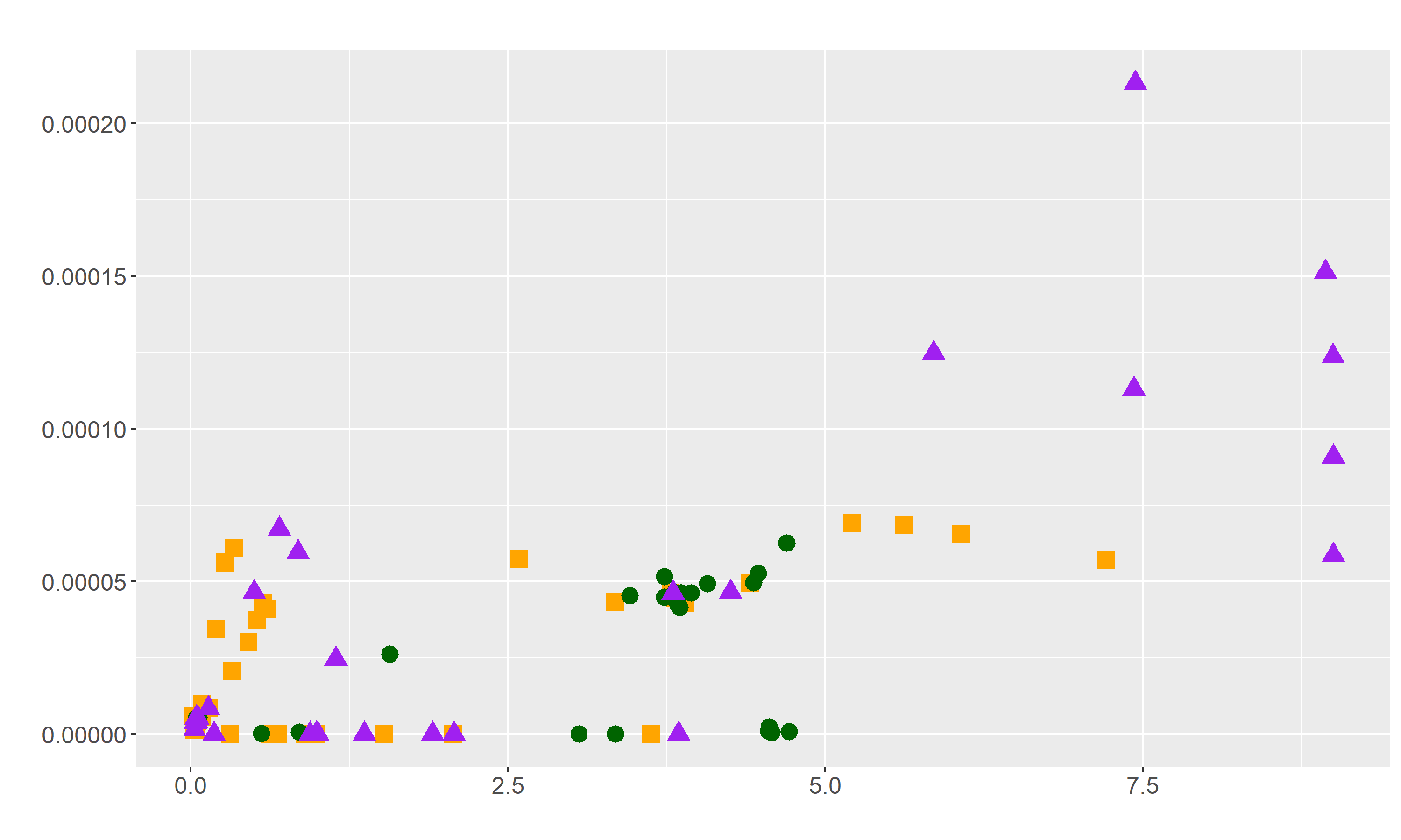}};
			\node(s2) at (2,0){\includegraphics[width=1cm]{figs/leyend_global.png}};
			\node[rotate=90](y) at (-3.45,0){{\color{darkgray}\tiny{$\log(\text{NRMSE}+1)$}}};
			\node(x) at (0,-2){{\color{darkgray}\tiny{MaxRE}}};
		\end{tikzpicture}
		\caption{Distance to ref. solution vs error when ``solved''.}
\end{subfigure}
\caption{Studying flatness in problem \hivP.}
\label{fig:flat}
\end{figure}

If the objective function is very flat near the reference solution, convergence becomes harder and one might need to specify more demanding stopping criteria in order to get better solutions. One example of this can be seen in Figure~\ref{fig:flat} for problem \hivP.

The two plots represent, for global optimization solvers, the distance to the reference solution on the $x$-axis. On the left we have running time on the $y$-axis and, on the right, the $y$-axis represents the error for the instances in which global optimality was ``certified''. The plot of the left shows that most of the configurations successfully terminated within the time limit certifying global optimality of the final solution. Differently from the situation for \daisyF in Figure~\ref{fig:local}, many of the obtained solutions are still relatively far from the reference solution. Yet, the plot on the right shows that all of these solutions still provide very small errors, and so imposing more strict termination criteria for the optimizers might mitigate this behavior for such a problem.

\subsection{Lack of identifiability}

Lack of identifiability can be an important challenge for parameter estimation problems. It arises when different parameter configurations deliver fits with the same quality. In other words, lack of identifiability corresponds to situations in which the underlying optimization problem has multiple global optima. It seems that lack of identifiability is not present in any of the problems we have studied since, as can be seen in all subfigures in \ref{app:plots}, it is never the case that global solvers provide very different solutions for a given problem beyond what we have classified above as flatness in the discussion of \hivP in Figure~\ref{fig:flat}. This result is in agreement with the findings of \citep{Go2023}, where several of the problems considered here were shown to be globally identifiable.

\section{Conclusions}
\label{sec:conclusion}
The main goal of this paper was to study the potential of mathematical programming modeling and state-of-the-art solvers to handle well known parameter estimation problems. The results are very promising, significantly surpassing the capacity to solve these problems reported in past literature. Our study has also helped to assess different modeling aspects such as the chosen discretization scheme or different variants on the formulation of the different optimization problems.

A natural direction for future research is to explore the robustness and scalability of the analysis, systematically introducing noisy measurements and studying significantly larger problems. In this respect, it is worth emphasizing that in our numerical analysis we have not relied on any sort of parallelization, which could definitely be exploited by solvers with these kind of capabilities such as \Octeract and \Knitro. Further, it is also worth noting that all the analysis developed in this paper has been done using the default settings of each solver and, thus, fine tuning these settings for each pair solver-problem might enable the solution of significantly larger problems. A related approach would be to develop a methodology for identifying the most promising configurations for solving a given problem, such as the choice of discretization scheme, mathematical programming formulation, or solver.

The analysis we have developed here can also serve as a benchmark for evaluating the performance of new specialized algorithms, specially global optimization ones, since they can now be compared to solvers such as \Couenne, \BARON, and \Octeract.

\section*{Declarations}

\subsection*{Ethics approval and consent to participate} 
Not applicable.

\subsection*{Consent for publication} 
Not applicable.

\subsection*{Funding} 
This work is part of the R\&D projects  PID2021-124030NB-C31 and PID2021-124030NB-C32 funded by MICIU/AEI/10.13039/501100011033/ and by ERDF/EU. This research was also funded by Grupos de Referencia Competitiva ED431C-2021/24 from the Consellería de Cultura, Educación e Universidades, Xunta de Galicia. JRB acknowledges support from grant PID2020-117271RB-C22 (BIODYNAMICS) funded by MCIN/AEI/10.13039/501100011033, 
from grant PID2023-146275NB-C22 (DYNAMO-bio) funded by MICIU/AEI/ 10.13039/501100011033 and ERDF/EU, and from grant CSIC PIE 202470E108 (LARGO). The authors acknowledge CESGA (Centro de Supercomputación de Galicia) for providing access to its FinisTerrae III supercomputer.

\subsection*{Availability of data and materials} 
Upon acceptance, all AMPL models used in the computational study will be made publicly available according to the journal standards.

\subsection*{Competing interests} 
The authors declare that they have no competing interests.

\subsection*{Authors' contributions}
M.F-D implemented the different solution methodologies as well as their adjustment to the different problems under study. M.F-D also carried out the computational study and collected the numerical results.

All five authors participated in the selection of the target problems, the design of the mathematical programming formulations, the design of the computational study, and the analysis, interpretation and discussion of the numerical results.

\subsection*{Acknowledgements}
The authors thank Eliseo Pita Vilariño for his contribution to the automation of the filtering and the processing of the computational results.

 \bibliographystyle{elsarticle-harv} 
 \bibliography{OptimEDO-refs}

\appendix

\section{Detailed description of the problems}
\label{app:description}
This appendix contains the detailed description of the parameter estimation problems used in the numerical experiments.

\section*{\textbf{Problem:} \alphapinene}

\begin{description}[labelindent=1cm]
	\item[ODE system:] \mbox{} \vspace{-1.5\abovedisplayskip}
	\begin{eqnarray*}
	\dfrac{dx_{1}}{dt} & = & -(p_1+p_2)x_{1} \\
	\dfrac{dx_{2}}{dt} & = & p_1x_{1} \\
	\dfrac{dx_{3}}{dt} & = & p_2x_{1}-(p_3+p_4)x_3+p_5x_5 \\
	\dfrac{dx_{4}}{dt} & = & p_3x_{3} \\
	\dfrac{dx_{5}}{dt} & = & -p_4x_{3}+p_{5}x_5
	\end{eqnarray*}
	\item[Initial conditions:] $x_1(0)= 100, x_2(0)=x_3(0)=x_4(0)=x_5(0)=0$.
	\item[Observed states:] $y_i(t)=x_i(t), \, \,i=1,\ldots,5$.
	\item[Measurements for each observed variable:] The time interval considered is $[0, 36900]$, with 8 measurements at time points 1230, 3060, 4920, 7800, 10680, 15030, 22620, and 36420.
	\item[Parameter bounds:] $0\leq p_i \leq 1, \,i=1,\ldots,5.$.
\end{description}


\section*{\textbf{Problem:} \BBG}

\begin{description}[labelindent=1cm]
	\item[ODE system:] \mbox{} \vspace{-1.5\abovedisplayskip}
	\begin{eqnarray*}
	\dfrac{dC_{b}}{dt} & = & \mu_{max}\dfrac{C_{s}C_{b}}{K_{s}+C_{s}}-k_{d}C_{b} \\
	\dfrac{dC_{s}}{dt} & = & -\dfrac{\mu_{max}}{yield}\dfrac{C_{s}C_{b}}{K_{s}+C_{s}} 
	\end{eqnarray*}
	\item[Initial conditions:] $C_{b}(0) = 2$ and $C_{s}(0) = 30$.
	\item[Observed states:] $y_1(t) = C_{b}(t)$ and $y_2(t) = C_{s}(t)$.
	\item[Measurements for each observed variable:] 7 measurements uniformly distributed in the interval $[0,12]$.
	\item[Parameter bounds:] $0.0001\leq p_i \leq 100, \,i=1,\ldots,4$, where $(p_1,p_2,p_3,p_4)=(\mu_{max},K_{s},k_{d},yield)$.
\end{description}


\section*{\textbf{Problem:} \FHN}

\begin{description}[labelindent=1cm]
	\item[ODE system:] \mbox{} \vspace{-1.5\abovedisplayskip}
	\begin{eqnarray*}
	\dfrac{dV}{dt} & = & g\left(V-\dfrac{V^{3}}{3}+R\right) \\
	\dfrac{dR}{dt} & = & \dfrac{1}{g}(V-a+bR) \\
	\end{eqnarray*}
	\item[Initial conditions:] $V(0) = -1$ and $R(0) = 1$.
	\item[Observed states:] $y_1(t) = V(t)$.
	\item[Measurements for each observed variable:] 6 measurements uniformly distributed in the interval $[0,20]$.
	\item[Parameter bounds:] $10^{-5}\leq p_i \leq 10^{5}, \,i=1,\ldots,3$, where $(p_1,p_2,p_3)=(g,a,b)$.
\end{description}


\section*{{\textbf{Problem:} \harmonic}}

\begin{description}[labelindent=1cm]
	\item[ODE system:] \mbox{} \vspace{-1.5\abovedisplayskip}
	\begin{eqnarray*}
	\dfrac{dx_{1}}{dt} & = & -p_1x_{2} \\
	\dfrac{dx_{2}}{dt} & = & \dfrac{1}{p_2}x_{1}
	\end{eqnarray*}
	\item[Initial conditions:] $0 \leq x_{i}(0)\leq 1.5, \,i=1,2$.
	\item[Observed states:] $y_i(t)=x_i(t), \, \,i=1,2$.
	\item[Measurements for each observed variable:] 10 measurements uniformly distributed in the interval $[0,2.3]$.
	\item[Parameter bounds:] $0.0001\leq p_i \leq 10, \,i=1,2$.
\end{description}


\section*{\textbf{Problem:} \lv}

\begin{description}[labelindent=1cm]
	\item[ODE system:] \mbox{} \vspace{-1.5\abovedisplayskip}
	\begin{eqnarray*}
	\dfrac{dr}{dt} & = & k_1r-k_2rw \\
	\dfrac{dw}{dt} & = & k_2rw-k_3w \\
	\end{eqnarray*}
	\item[Initial conditions:] $90 \leq r(0) \leq 110$ and $90 \leq w(0) \leq 110$.
	\item[Observed states:] \mbox{} \vspace{-1.5\abovedisplayskip}
	\[
	\begin{array}{l|r}
	\text{\lvF}   & \text{\lvP} \\
	y_1(t) = r(t) & y_1(t) = r(t) \\
	y_2(t) = w(t) &              
	\end{array}
	\]
	\item[Measurements for each observed variable:] 20 measurements uniformly distributed in the interval $[0,1]$.
	\item[Parameter bounds:] $0.0001 \leq p_i \leq 1, \,i=1,\ldots,3$, where $(p_1,p_2,p_3)=(k_1,k_2,k_3)$.
\end{description}


\section*{\textbf{Problem:} \daisy}
\begin{description}[labelindent=1cm]
	\item[ODE system:] \mbox{} \vspace{-1.5\abovedisplayskip}
	\begin{eqnarray*}
	\dfrac{dx_{1}}{dt} & = & -(a_{21}+a_{31}+a_{01})x_1 + a_{12}x_2+a_{13}x_3 \\
	\dfrac{dx_{2}}{dt} & = & a_{21}x_1 - a_{12}x_2 \\
	\dfrac{dx_{3}}{dt} & = & a_{31}x_1 - a_{13}x_3 
	\end{eqnarray*}
	\item[Initial conditions:] \mbox{} \vspace{-1.5\abovedisplayskip}
	\begin{eqnarray*}
	-1 \leq & x_1(0) & \leq 2  \\
	-1 \leq & x_2(0) & \leq 2  \\
	-1 \leq & x_3(0) & \leq 2 
	\end{eqnarray*}
	\item[Observed states:] \mbox{} \vspace{-1.5\abovedisplayskip}
	\begin{equation*}
	\begin{array}{l|r}
	\text{\daisyF}   & \text{\daisyP} \\
	y_1(t) = x_1(t) & y_1(t) = x_1(t) \\
	y_2(t) = x_2(t) & y_2(t) = x_2(t) \\
	y_3(t) = x_3(t) &               
	\end{array}
	\end{equation*}
	\item[Measurements for each observed variable:] 20 measurements uniformly distributed in the interval $[0,1]$.
	\item[Parameter bounds:] $-1 \leq p_i \leq 2, \,i=1,\ldots,5$,
where $(p_1,p_2,p_3,p_4,p_5)=(a_{21},a_{31},a_{01},a_{12},a_{13})$.
\end{description}


\section*{\textbf{Problem:} \hiv}
\begin{description}[labelindent=1cm]
	\item[ODE system:] \mbox{} \vspace{-1.5\abovedisplayskip}
	\begin{eqnarray*}
	\dfrac{dx}{dt} & = & \lambda-dx-\beta xv\\
	\dfrac{dy}{dt} & = & \beta xv - ay \\
	\dfrac{dv}{dt} & = & ky-uv \\
	\dfrac{dw}{dt} & = & cxyw-cqyw-bw \\ 
	\dfrac{dz}{dt} & = & cqyw-hz \\
	\end{eqnarray*}
	\item[Initial conditions:] \mbox{} \vspace{-1.5\abovedisplayskip}
	\begin{eqnarray*}
	0.001 \leq & x(0) & \leq 2 \\
	0.001 \leq & y(0) & \leq 2 \\
	0.001 \leq & v(0) & \leq 2 \\
	0.001 \leq & w(0) & \leq 2 \\
	0.001 \leq & z(0) & \leq 2
	\end{eqnarray*}
	\item[Observed states:] \mbox{} \vspace{-1.5\abovedisplayskip}
	\begin{equation*}
	\begin{array}{l|l}
	\text{\hivF}   & \text{\hivP} \\
	y_1(t) = x(t)      & y_1(t) = x(t) \\
	y_2(t) = y(t)      & y_2(t) = y(t)+v(t) \\
	y_3(t) = v(t)      & y_3(t) = w(t) \\   
	y_4(t) = w(t)      & y_4(t) = z(t) \\    
	y_5(t) = z(t)      &         
	\end{array}
	\end{equation*}
	\item[Measurements for each observed variable:] 20 measurements uniformly distributed in the interval $[0,10]$.
	\item[Parameter bounds:] $0.0001 \leq p_i \leq 1, \,i=1,\ldots,10$, where $(p_{1},p_{2},\ldots,p_{10})=(\lambda,d,\ldots,h)$.
\end{description}


\section*{\textbf{Problem:} \crauste}

\begin{description}[labelindent=1cm]
	\item[ODE system:] \mbox{} \vspace{-1.5\abovedisplayskip}
	\begin{eqnarray*}
	\dfrac{dN}{dt} & = & \mu_{N}N-\delta_{NE}NP \\
	\dfrac{dE}{dt} & = & \delta_{NE}NP-\mu_{EE}E^{2}-\delta_{EL}E+\rho_{E}EP \\
	\dfrac{dS}{dt} & = & \delta_{EL}S-S\delta_{LM}-\mu_{LL}S^2-\mu_{LE}ES \\
	\dfrac{dM}{dt} & = & \delta_{LM}S-\mu_{M}M \\
	\dfrac{dP}{dt} & = & \rho_{P}P^2-\mu_{P}P-\mu_{PE}EP-\mu_{PL}SP \\
	\end{eqnarray*}

	\item[Initial conditions:] \mbox{} \vspace{-1.5\abovedisplayskip}
	\begin{eqnarray*}
	-1.1 \leq & N(0) & \leq 1.1  \\
	-1.1 \leq & E(0) & \leq 1.1  \\
	-1.1 \leq & S(0) & \leq 1.1  \\
	-1.1 \leq & M(0) & \leq 1.1  \\
	-1.1 \leq & P(0) & \leq 1.1  
	\end{eqnarray*}
	\item[Observed states:] \mbox{} \vspace{-1.5\abovedisplayskip}
	\begin{equation*}
	\begin{array}{l|l}
	\text{\crausteF}   & \text{\crausteP} \\
	y_1(t) = N(t) & y_1(t) = N(t) \\
	y_2(t) = E(t) & y_2(t) = E(t)  \\
	y_3(t) = S(t) & y_3(t) = S(t)+M(t) \\
	y_4(t) = M(t) & y_4(t) = P(t)   \\
	y_5(t) = P(t) &               
	\end{array}
	\end{equation*}
	\item[Measurements for each observed variable:] 20 measurements uniformly distributed in the interval $[0,1]$.
	\item[Parameter bounds:] $-2 \leq p_i \leq 2, \,i=1,\ldots,13$, where $(p_1,p_2,p_3,\ldots,p_{13})=(\mu_{N},\delta_{NE},\mu_{EE},\ldots,\mu_{PL})$.
\end{description}


\section{Numerical results: Digging deeper}
\label{app:performanceplus}

The goal of this section is to get additional insights from the computational analysis by by studying different disaggregations of the numerical results. The objective here is to get a better understanding, not only of the individual elements involved in the different configurations, but also of any potential interactions between them. For the sake of exposition, here we only present the top performing configurations, with respect to \foundref, for each considered disaggregation level.\footnote{The full tables are available in the Online Appendix.}

The results highlight once again that there is no need to fine tune a specific configuration to be able to successfully solve a given problem, since in most of them the reference solution is solved with a wide variety of solvers and configurations.

\begin{table}[!htbp]
\centering
\footnotesize{\begin{tabular}{lcccccc}
\toprule
Problem --- solver & \solvedS & \foundref & \nearref & \alternat & \timeF & \success \\
\midrule
\daisyF{} --- \BARON & 0.880 & 1.000 & 0\phantom{.000} & 0\phantom{.000} & 0.076 & 1.000 \\
\daisyF{} --- \Couenne & 1.000 & 1.000 & 0\phantom{.000} & 0\phantom{.000} & 0.220 & 1.000 \\
\harmonic{} --- \Octeract & 0.720 & 0.980 & 0.020 & 0\phantom{.000} & 1.286 & 1.000 \\
\harmonic{} --- \BARON & 0.760 & 0.960 & 0\phantom{.000} & 0\phantom{.000} & 0.143 & 1.000 \\
\daisyF{} --- \BONMIN & 1.000 & 0.960 & 0\phantom{.000} & 0\phantom{.000} & 0.046 & 1.000 \\
\daisyF{} --- \CONOPT & 0\phantom{.000} & 0.920 & 0.080 & 0\phantom{.000} & 0.072 & 1.000 \\
\daisyF{} --- \Ipopt & 1.000 & 0.920 & 0\phantom{.000} & 0\phantom{.000} & 0.056 & 1.000 \\
\lvF{} --- \BONMIN & 1.000 & 0.900 & 0\phantom{.000} & 0\phantom{.000} & 0.038 & 1.000 \\
\lvF{} --- \Ipopt & 0.980 & 0.880 & 0\phantom{.000} & 0\phantom{.000} & 0.038 & 0.980 \\
\lvF{} --- \BARON & 0.580 & 0.860 & 0\phantom{.000} & 0\phantom{.000} & 0.097 & 1.000 \\
\bottomrule
\end{tabular}}
\caption{Summary results by problem and solver.}
\label{tab:byproblemsolver}
\end{table}

Table~\ref{tab:byproblemsolver} confirms that \daisyF, \harmonic and \lvF are the easiest problems to solve and, moreover, that global optimization solvers seem to be better at finding the reference solution.

\begin{table}[!htbp]
\centering
\footnotesize{\begin{tabular}{lcccccc}
\toprule
Problem --- scheme & \solvedS & \foundref & \nearref & \alternat & \timeF & \success \\
\midrule
\lvF{} --- \simpson & 0.775 & 0.988 & 0.012 & 0\phantom{.000} & 0.028 & 1.000 \\
\lvF{} --- \trapecio & 0.775 & 0.975 & 0.025 & 0\phantom{.000} & 0.026 & 1.000 \\
\daisyF{} --- \euler & 0.700 & 0.925 & 0.025 & 0\phantom{.000} & 0.029 & 1.000 \\
\daisyF{} --- \trapecio & 0.750 & 0.925 & 0.025 & 0\phantom{.000} & 0.031 & 1.000 \\
\lvF{} --- \adams & 0.863 & 0.900 & 0.012 & 0\phantom{.000} & 0.029 & 0.912 \\
\daisyF{} --- \simpson & 0.725 & 0.887 & 0.012 & 0\phantom{.000} & 0.038 & 0.938 \\
\daisyF{} --- \runge & 0.750 & 0.875 & 0.025 & 0.012 & 0.093 & 0.938 \\
\lvF{} --- \runge & 0.725 & 0.825 & 0.050 & 0\phantom{.000} & 0.102 & 0.900 \\
\daisyF{} --- \adams & 0.750 & 0.812 & 0.025 & 0\phantom{.000} & 0.055 & 0.950 \\
\harmonic{} --- \trapecio & 0.787 & 0.800 & 0\phantom{.000} & 0\phantom{.000} & 0.039 & 0.950 \\
\bottomrule
\end{tabular}}
\caption{Summary results by problem and discretization scheme.}
\label{tab:byproblemdiscretization}
\end{table}

Table~\ref{tab:byproblemdiscretization} shows that, although \euler seems to be the worst performing discretization scheme overall, it can be particularly effective for specific problems. In particular, it is the best performing scheme in \daisyF, tied with \trapecio.

\begin{table}[!htbp]
\centering
\footnotesize{\begin{tabular}{lcccccc}
\toprule
Problem --- form. & \solvedS & \foundref & \nearref & \alternat & \timeF & \success \\
\midrule
\daisyF{} --- \modeA & 0.725 & 0.938 & 0\phantom{.000} & 0\phantom{.000} & 0.032 & 0.975 \\
\hivF{} --- \modeA & 0.838 & 0.925 & 0\phantom{.000} & 0\phantom{.000} & 0.056 & 0.938 \\
\daisyF{} --- \modeB & 0.738 & 0.875 & 0.056 & 0.006 & 0.029 & 0.956 \\
\daisyF{} --- \modeC & 0.738 & 0.869 & 0\phantom{.000} & 0\phantom{.000} & 0.040 & 0.969 \\
\harmonic{} --- \modeB & 0.812 & 0.863 & 0\phantom{.000} & 0\phantom{.000} & 0.034 & 0.981 \\
\lvF{} --- \modeC & 0.738 & 0.844 & 0.037 & 0\phantom{.000} & 0.026 & 0.975 \\
\lvF{} --- \modeB & 0.769 & 0.838 & 0.019 & 0\phantom{.000} & 0.030 & 0.963 \\
\lvF{} --- \modeA & 0.700 & 0.812 & 0\phantom{.000} & 0\phantom{.000} & 0.044 & 0.925 \\
\harmonic{} --- \modeC & 0.619 & 0.650 & 0\phantom{.000} & 0\phantom{.000} & 0.120 & 0.894 \\
\hivF{} --- \modeB & 0.856 & 0.644 & 0.150 & 0.206 & 0.033 & 1.000 \\
\bottomrule
\end{tabular}}
\caption{Summary results by problem and mathematical programming formulation.}
\label{tab:byproblemformulation}
\end{table}

Interestingly, Table~\ref{tab:byproblemformulation} shows that the choice of the mathematical programming formulation can have a significant impact. \modeA performs extremely well in \hivF and \daisyF problems, finding the reference solutions in more than 90\% of the configurations, whereas the next best performing formulation for \hivF is under 65\%. On the other hand, \modeB is clearly the best formulation for \harmonic, whereas \modeC comes out narrowly on top for \lvF.\footnote{The Online Appendix contains the full tables.} 

\begin{table}[!htbp]
\centering
\footnotesize{\begin{tabular}{lcccccc}
\toprule
Solver --- scheme & \solvedS & \foundref & \nearref & \alternat & \timeF & \success \\
\midrule
\Couenne{} --- \trapecio & 0.795 & 0.591 & 0.144 & 0.098 & 0.138 & 0.917 \\
\BARON{} --- \trapecio & 0.758 & 0.576 & 0.114 & 0.182 & 0.078 & 1.000 \\
\BARON{} --- \adams & 0.758 & 0.568 & 0.129 & 0.144 & 0.102 & 0.985 \\
\Knitro{} --- \trapecio & 0.977 & 0.530 & 0.098 & 0.174 & 0.026 & 1.000 \\
\Octeract{} --- \trapecio & 0.568 & 0.523 & 0.106 & 0.167 & 1.286 & 0.939 \\
\Couenne{} --- \simpson & 0.780 & 0.523 & 0.182 & 0.083 & 0.148 & 0.894 \\
\BARON{} --- \simpson & 0.689 & 0.500 & 0.159 & 0.167 & 0.095 & 1.000 \\
\BARON{} --- \runge & 0.667 & 0.485 & 0.121 & 0.136 & 0.306 & 0.977 \\
\Couenne{} --- \adams & 0.720 & 0.477 & 0.121 & 0.106 & 0.236 & 0.924 \\
\BARON{} --- \euler & 0.462 & 0.462 & 0.152 & 0.212 & 0.076 & 1.000 \\
\bottomrule
\end{tabular}}
\caption{Summary results by solver and discretization scheme.}
\label{tab:solverdiscretization}
\end{table}

Table~\ref{tab:solverdiscretization} shows that \trapecio and is the best performing configuration for global optimization solvers solvers.

\begin{table}[!htbp]
\centering
\footnotesize{\begin{tabular}{lcccccc}
\toprule
Solver --- form. & \solvedS & \foundref & \nearref & \alternat & \timeF & \success \\
\midrule
\Knitro{} --- \modeA & 0.888 & 0.576 & 0.096 & 0.128 & 0.032 & 0.984 \\
\BARON{} --- \modeB & 0.786 & 0.572 & 0.182 & 0.133 & 0.076 & 0.996 \\
\Couenne{} --- \modeA & 0.664 & 0.536 & 0.096 & 0.040 & 0.212 & 0.784 \\
\BONMIN{} --- \modeA & 0.856 & 0.504 & 0.120 & 0.024 & 0.046 & 0.856 \\
\Couenne{} --- \modeB & 0.705 & 0.491 & 0.158 & 0.109 & 0.138 & 0.849 \\
\BARON{} --- \modeA & 0.624 & 0.488 & 0.096 & 0.264 & 0.134 & 1.000 \\
\BARON{} --- \modeC & 0.552 & 0.472 & 0.100 & 0.160 & 0.097 & 0.984 \\
\Octeract{} --- \modeA & 0.616 & 0.456 & 0.056 & 0.176 & 1.563 & 0.856 \\
\Knitro{} --- \modeB & 0.975 & 0.449 & 0.168 & 0.186 & 0.028 & 0.993 \\
\Ipopt{} --- \modeB & 0.902 & 0.449 & 0.154 & 0.175 & 0.056 & 0.937 \\
\bottomrule
\end{tabular}}
\caption{Summary results by solver and mathematical programming formulation.}
\label{tab:solverformulation}
\end{table}

Table~\ref{tab:solverformulation} shows that \modeA is often the best performing formulation for the different optimization solvers, although \modeB comes out on top for \BARON and \Ipopt, for instance.

\begin{table}[!htbp]
\centering
\footnotesize{\begin{tabular}{lcccccc}
\toprule
Scheme -- form. & \solvedS & \foundref & \nearref & \alternat & \timeF & \success \\
\midrule
\trapecio{} --- \modeA & 0.715 & 0.530 & 0.065 & 0.145 & 0.032 & 0.950 \\
\trapecio{} --- \modeB & 0.750 & 0.474 & 0.125 & 0.206 & 0.028 & 0.950 \\
\adams{} --- \modeA & 0.700 & 0.445 & 0.075 & 0.145 & 0.093 & 0.870 \\
\simpson{} --- \modeA & 0.635 & 0.430 & 0.105 & 0.130 & 0.046 & 0.870 \\
\adams{} --- \modeB & 0.724 & 0.425 & 0.149 & 0.189 & 0.041 & 0.897 \\
\simpson{} --- \modeB & 0.741 & 0.419 & 0.140 & 0.193 & 0.037 & 0.925 \\
\trapecio{} --- \modeC & 0.695 & 0.415 & 0.100 & 0.160 & 0.026 & 0.953 \\
\runge{} --- \modeA & 0.545 & 0.405 & 0.055 & 0.140 & 0.106 & 0.760 \\
\euler{} --- \modeA & 0.620 & 0.395 & 0.105 & 0.170 & 0.032 & 0.910 \\
\runge{} --- \modeB & 0.664 & 0.371 & 0.138 & 0.189 & 0.093 & 0.840 \\
\bottomrule
\end{tabular}}
\caption{Summary results by discretization scheme and mathematical programming formulation.}
\label{tab:discretizationformulation}
\end{table}

Finally, regarding the disaggregation by discretization scheme and mathematical programming formulation, we see in Table~\ref{tab:discretizationformulation} that the configurations with \trapecio tend to come out on top and \euler and \runge are at the bottom. Moreover, for all discretization schemes, the configurations with \modeA seem to dominate those with \modeB, which themselves seem to dominate those configurations with \modeC.

\clearpage

\section{Graphical representations of the results}
\label{app:plots}

\subsection{\alphapinene}
\label{sec:alpha}

\begin{figure}[!htbp]
\centering
\begin{subfigure}{0.495\textwidth}
		\centering
		\begin{tikzpicture}
			\node(s1) at (0,0){\includegraphics[width=0.9\textwidth]{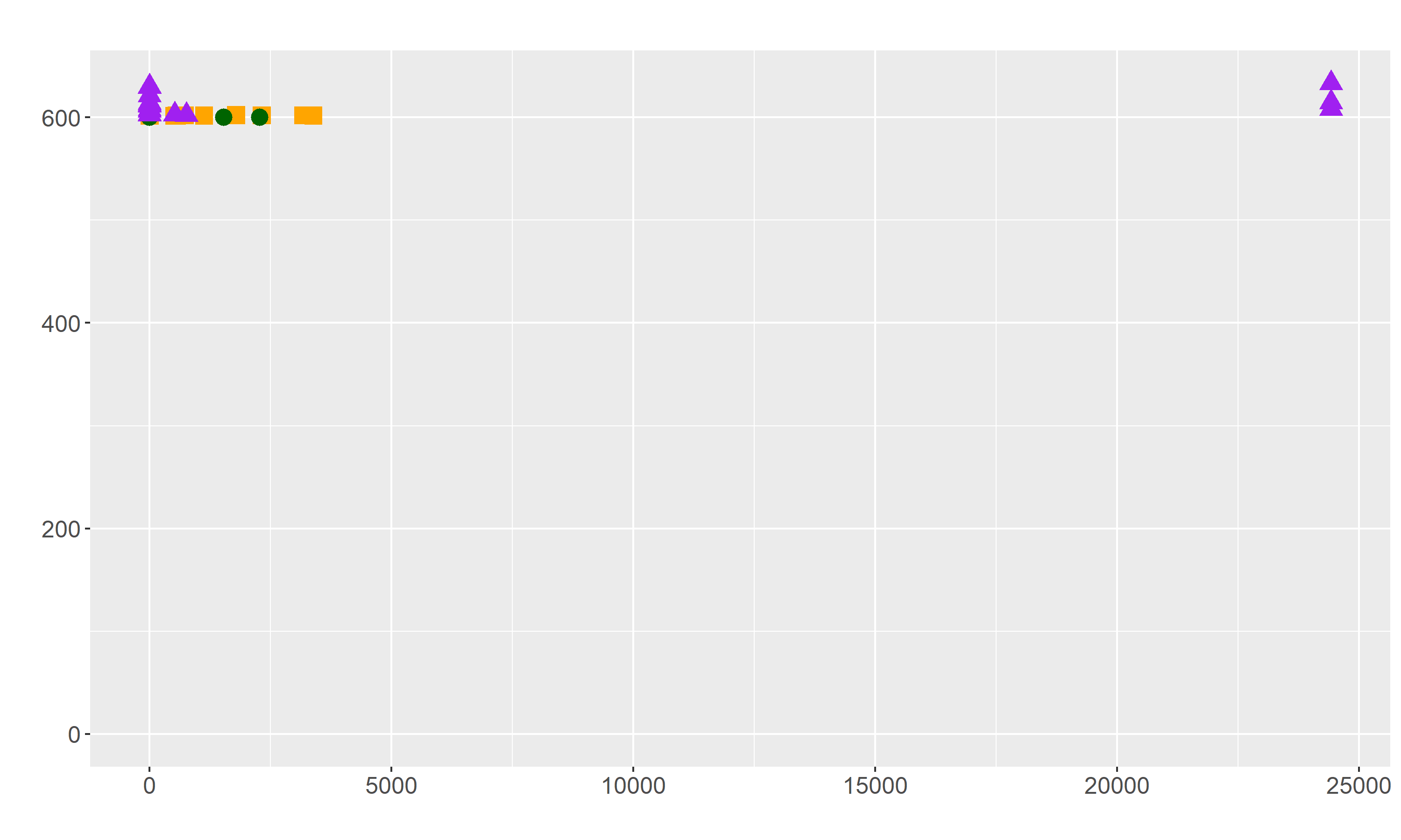}};
			\node[rotate=90](y) at (-3.35,0){{\color{darkgray}\tiny{time}}};
			\node(x) at (0,-2){{\color{darkgray}\tiny{distance}}};
		\end{tikzpicture}
		\caption{Global solvers: Distance to ref. solution vs time.}
\end{subfigure}
\begin{subfigure}{0.495\textwidth}
		\centering
		\begin{tikzpicture}
			\node(s1) at (0,0){\includegraphics[width=0.9\textwidth]{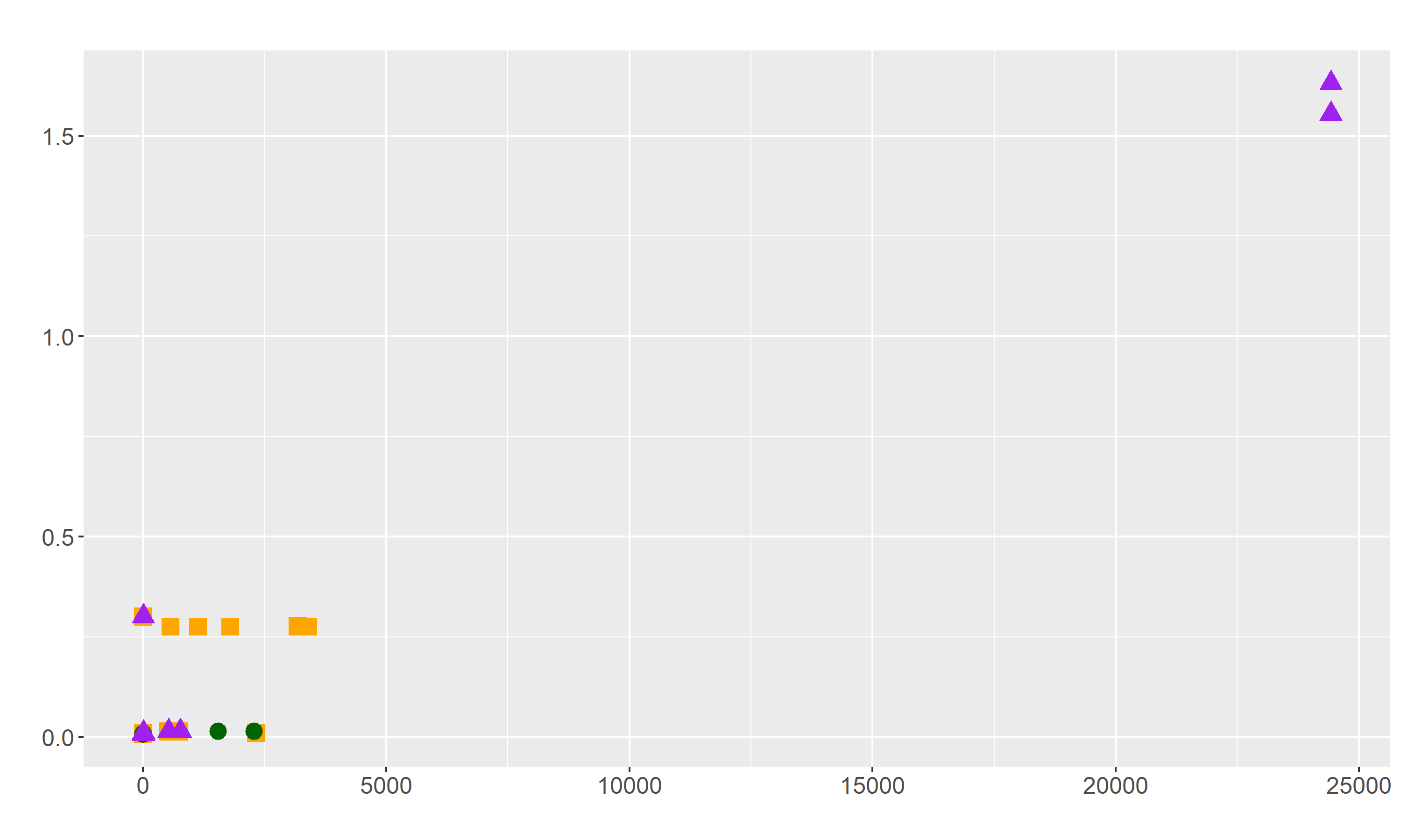}};
		\node(s2) at (2,0){\includegraphics[width=1cm]{figs/leyend_global.png}};
		\node[rotate=90](y) at (-3.4,0){{\color{darkgray}\tiny{$\log(\text{error}+1)$}}};
		\node(x) at (0,-2){{\color{darkgray}\tiny{distance}}};
		\end{tikzpicture}
		\caption{Global solvers: Distance to ref. solution vs error.}
\end{subfigure}

\vspace{0.2cm}

\begin{subfigure}{0.495\textwidth}
		\centering
		\begin{tikzpicture}
			\node(s1) at (0,0){\includegraphics[width=0.9\textwidth]{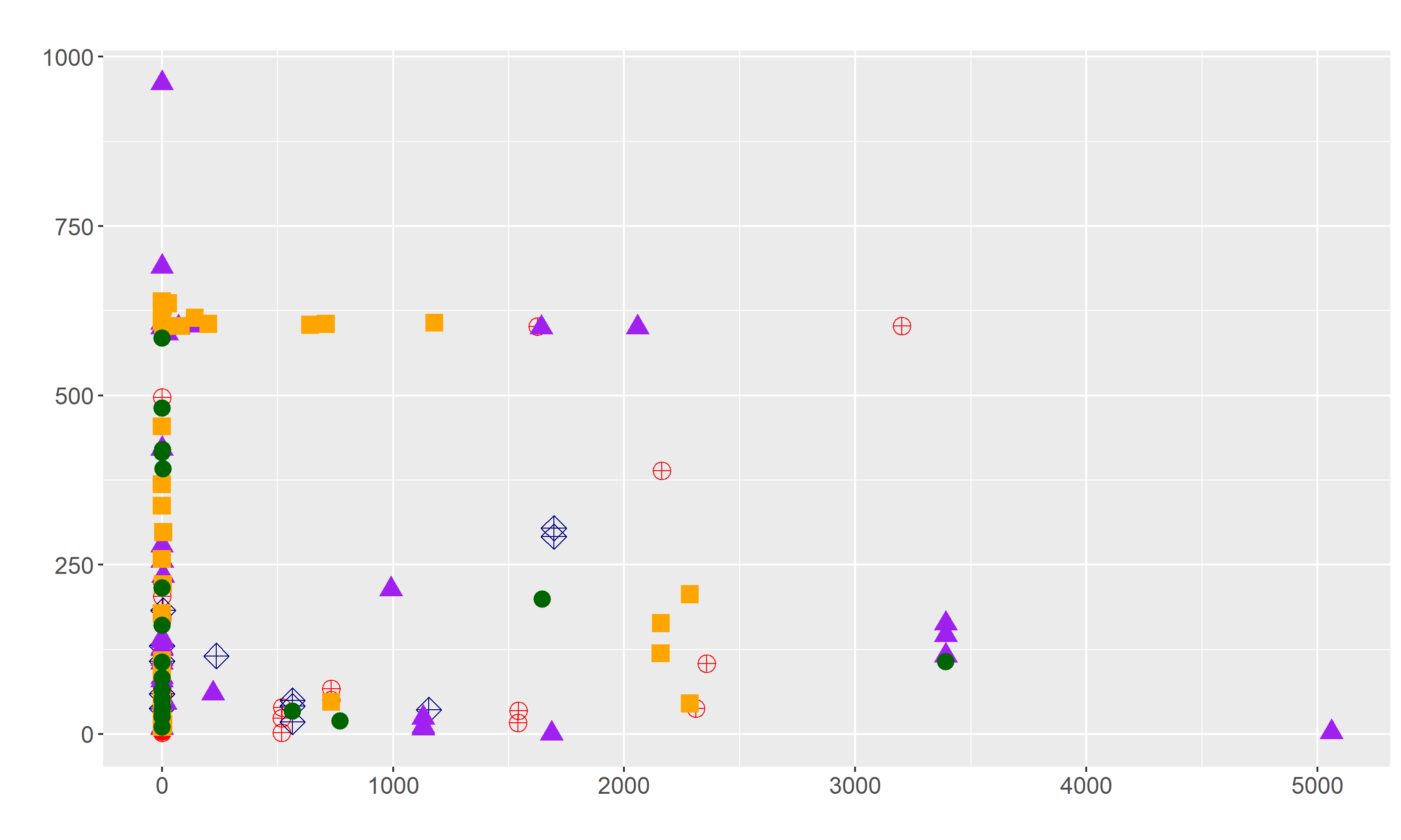}};
			\node[rotate=90](y) at (-3.35,0){{\color{darkgray}\tiny{time}}};
			\node(x) at (0,-2){{\color{darkgray}\tiny{distance}}};
		\end{tikzpicture}
		\caption{Local solvers: Distance to ref. solution vs time.}
\end{subfigure}
\begin{subfigure}{0.495\textwidth}
		\centering
		\begin{tikzpicture}
			\node(s1) at (0,0){\includegraphics[width=0.9\textwidth]{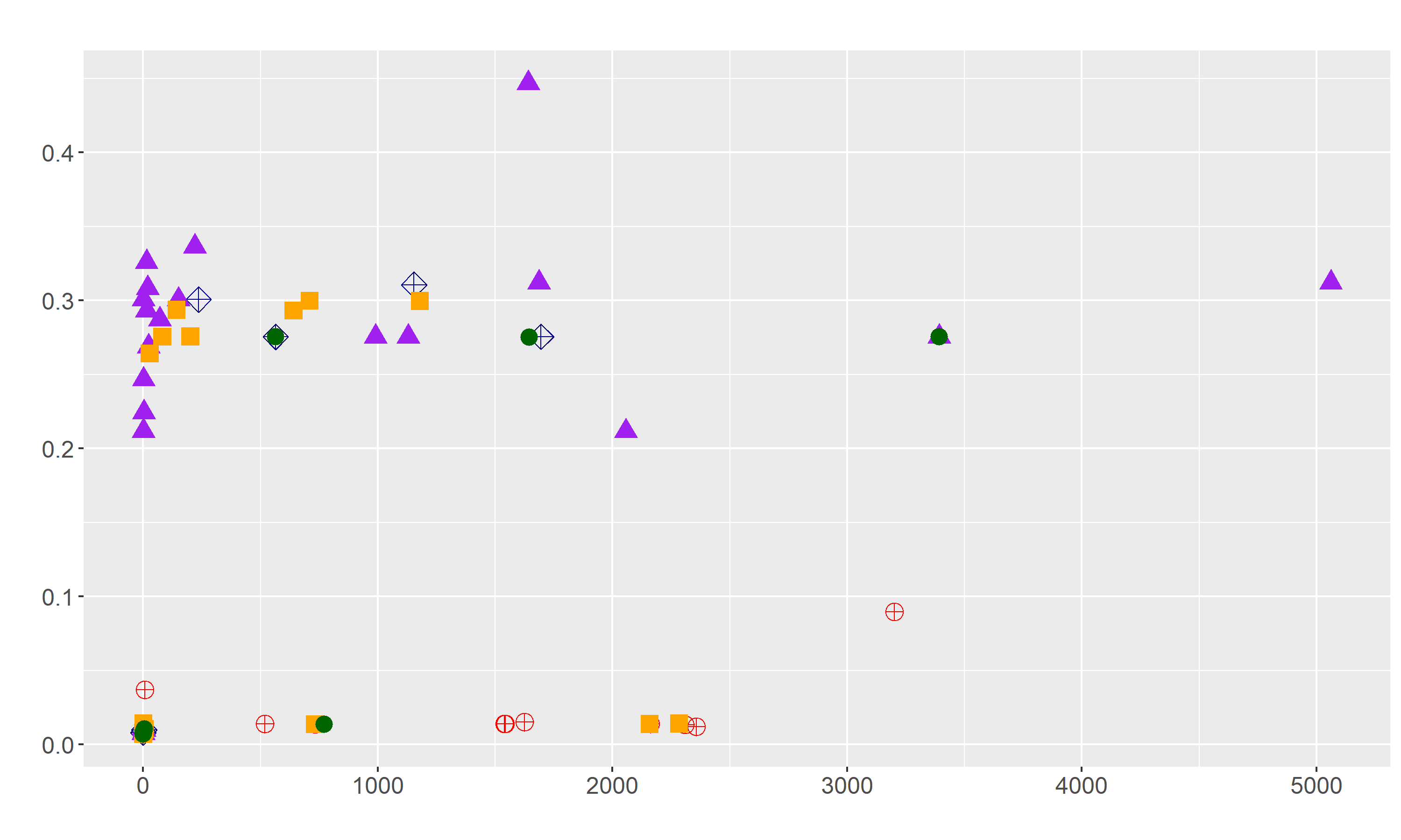}};
			\node(s2) at (2,0){\includegraphics[width=1cm]{figs/leyend_local.png}};
			\node[rotate=90](y) at (-3.4,0){{\color{darkgray}\tiny{$\log(\text{error}+1)$}}};
			\node(x) at (0,-2){{\color{darkgray}\tiny{distance}}};
		\end{tikzpicture}
		\caption{Local solvers: Distance to ref. solution vs error.}
\end{subfigure}

\caption{Results for problem \alphapinene.}
\label{fig:alphapinene}
\end{figure}

\newpage 
\subsection{\harmonic}
\begin{figure}[!htbp]
\centering
\begin{subfigure}{0.495\textwidth}
		\centering
		\begin{tikzpicture}
			\node(s1) at (0,0){\includegraphics[width=0.9\textwidth]{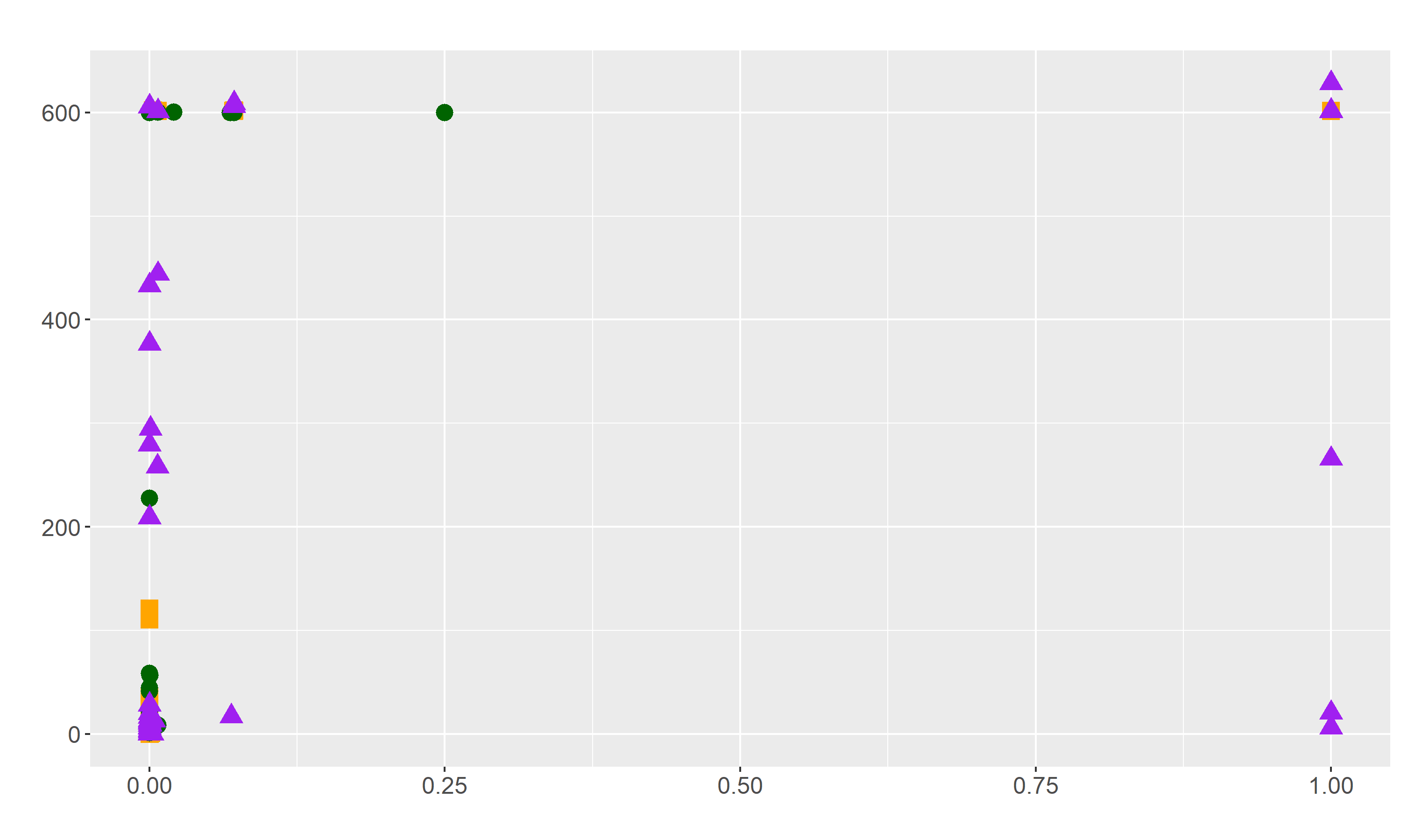}};
			\node[rotate=90](y) at (-3.35,0){{\color{darkgray}\tiny{time}}};
			\node(x) at (0,-2){{\color{darkgray}\tiny{distance}}};
		\end{tikzpicture}
		\caption{Global solvers: Distance to ref. solution vs time.}
\end{subfigure}
\begin{subfigure}{0.495\textwidth}
		\centering
		\begin{tikzpicture}
			\node(s1) at (0,0){\includegraphics[width=0.9\textwidth]{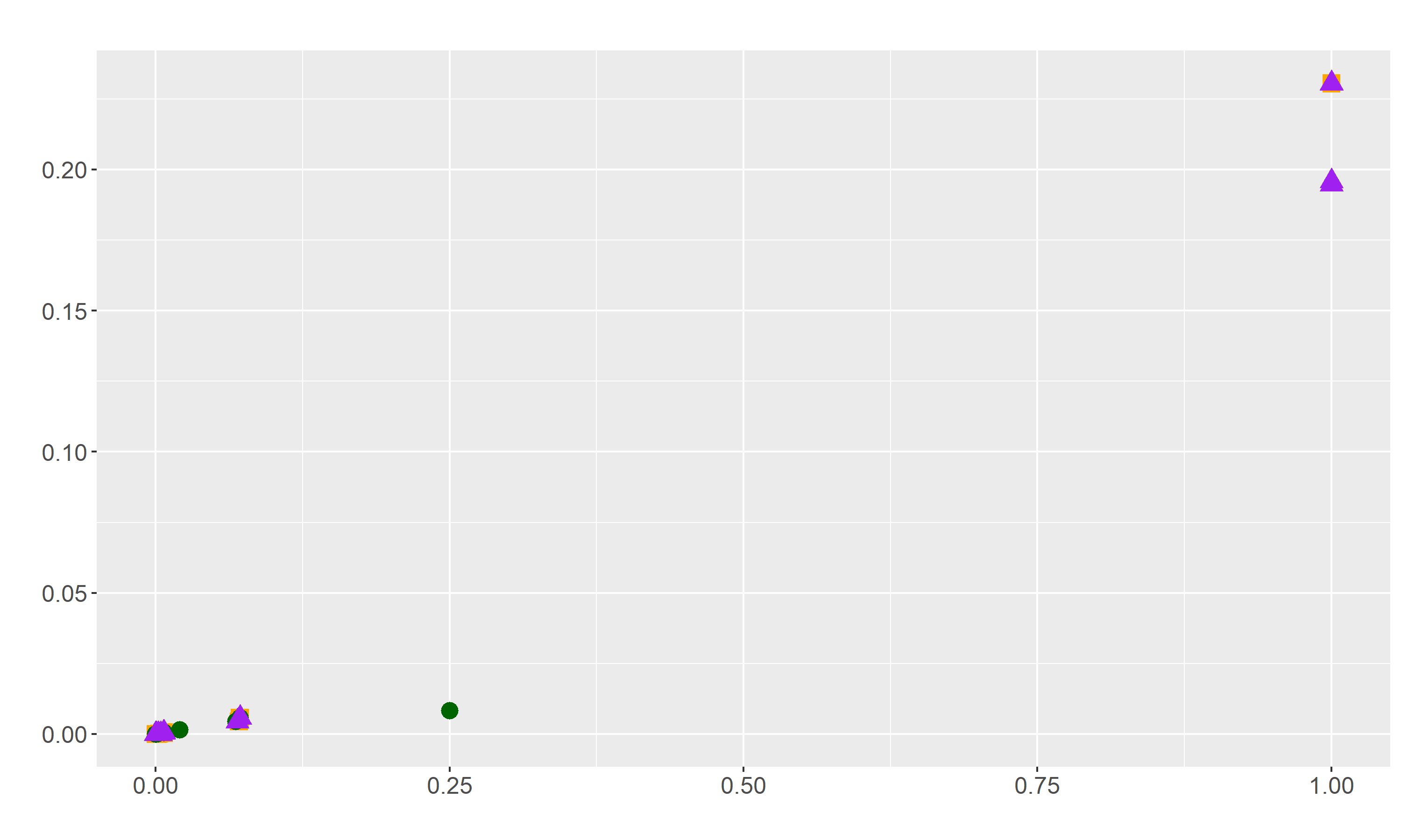}};
		\node(s2) at (2,0){\includegraphics[width=1cm]{figs/leyend_global.png}};
		\node[rotate=90](y) at (-3.4,0){{\color{darkgray}\tiny{$\log(\text{error}+1)$}}};
		\node(x) at (0,-2){{\color{darkgray}\tiny{distance}}};
		\end{tikzpicture}
		\caption{Global solvers: Distance to ref. solution vs error.}
\end{subfigure}

\vspace{0.2cm}

\begin{subfigure}{0.495\textwidth}
		\centering
		\begin{tikzpicture}
			\node(s1) at (0,0){\includegraphics[width=0.9\textwidth]{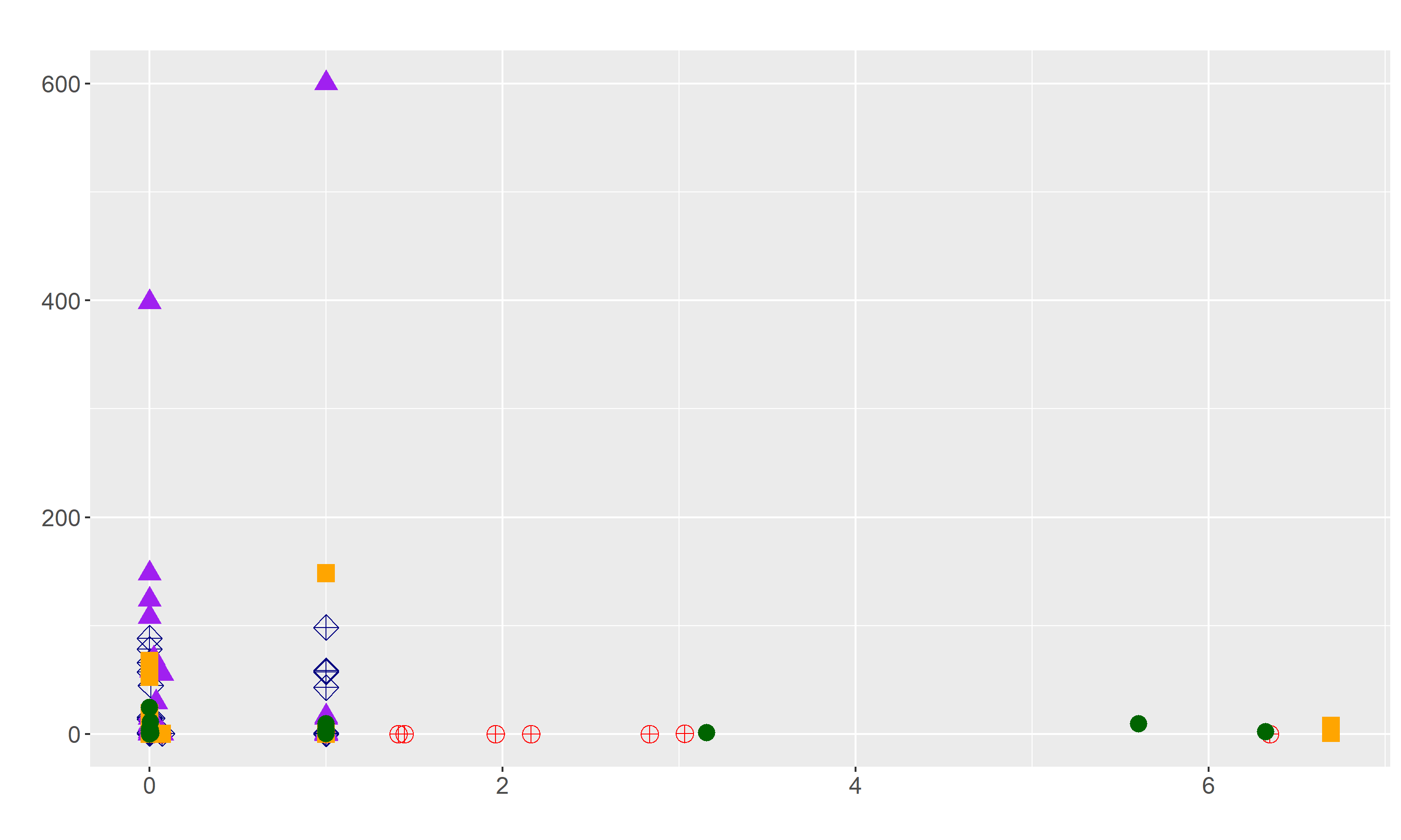}};
			\node[rotate=90](y) at (-3.35,0){{\color{darkgray}\tiny{time}}};
			\node(x) at (0,-2){{\color{darkgray}\tiny{distance}}};
		\end{tikzpicture}
		\caption{Local solvers: Distance to ref. solution vs time.}
\end{subfigure}
\begin{subfigure}{0.495\textwidth}
		\centering
		\begin{tikzpicture}
			\node(s1) at (0,0){\includegraphics[width=0.9\textwidth]{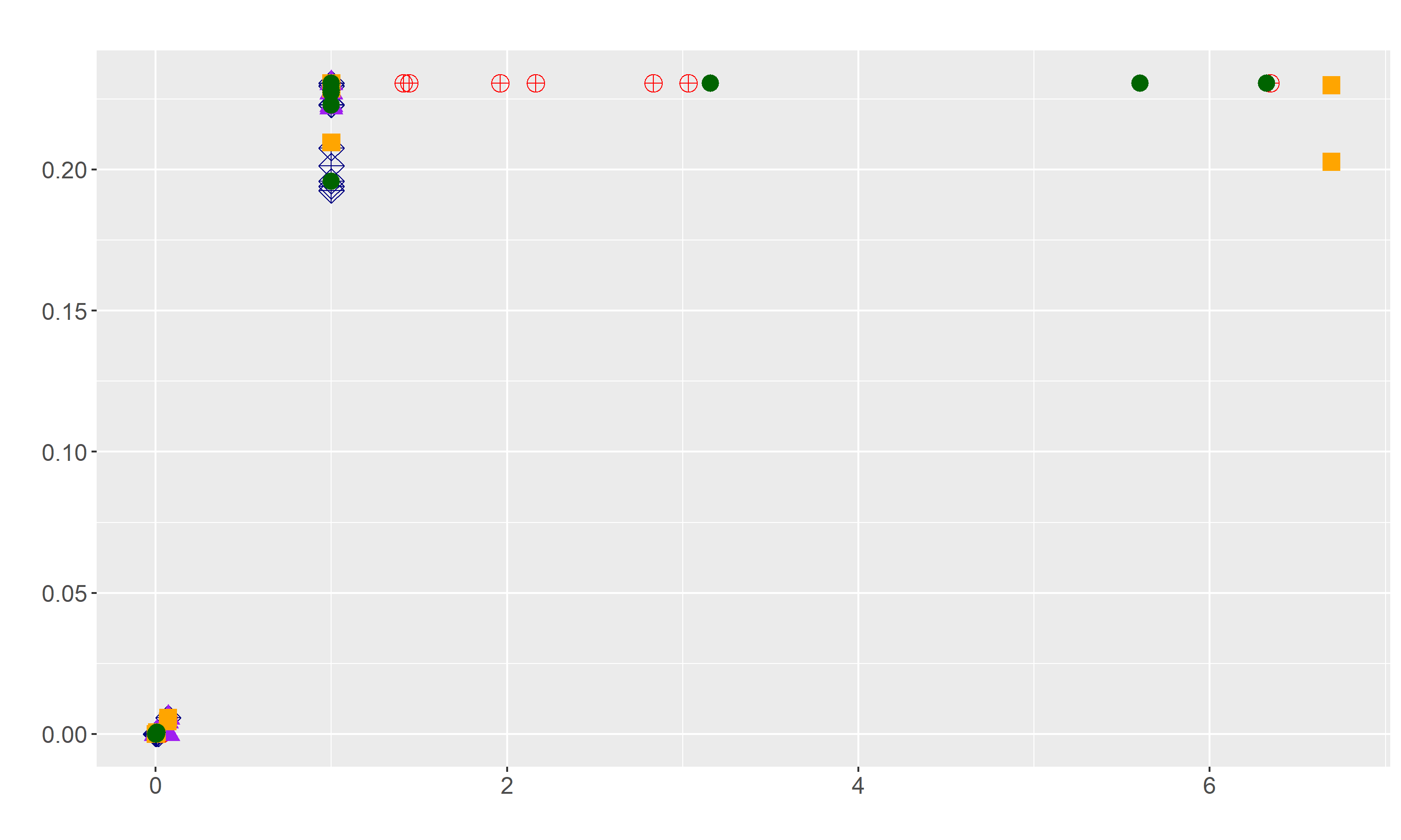}};
			\node(s2) at (2,0){\includegraphics[width=1cm]{figs/leyend_local.png}};
			\node[rotate=90](y) at (-3.4,0){{\color{darkgray}\tiny{$\log(\text{error}+1)$}}};
			\node(x) at (0,-2){{\color{darkgray}\tiny{distance}}};
		\end{tikzpicture}
		\caption{Local solvers: Distance to ref. solution vs error.}
\end{subfigure}

\begin{subfigure}{0.495\textwidth}
		\centering
		\begin{tikzpicture}
			\node(s1) at (0,0){\includegraphics[width=0.9\textwidth]{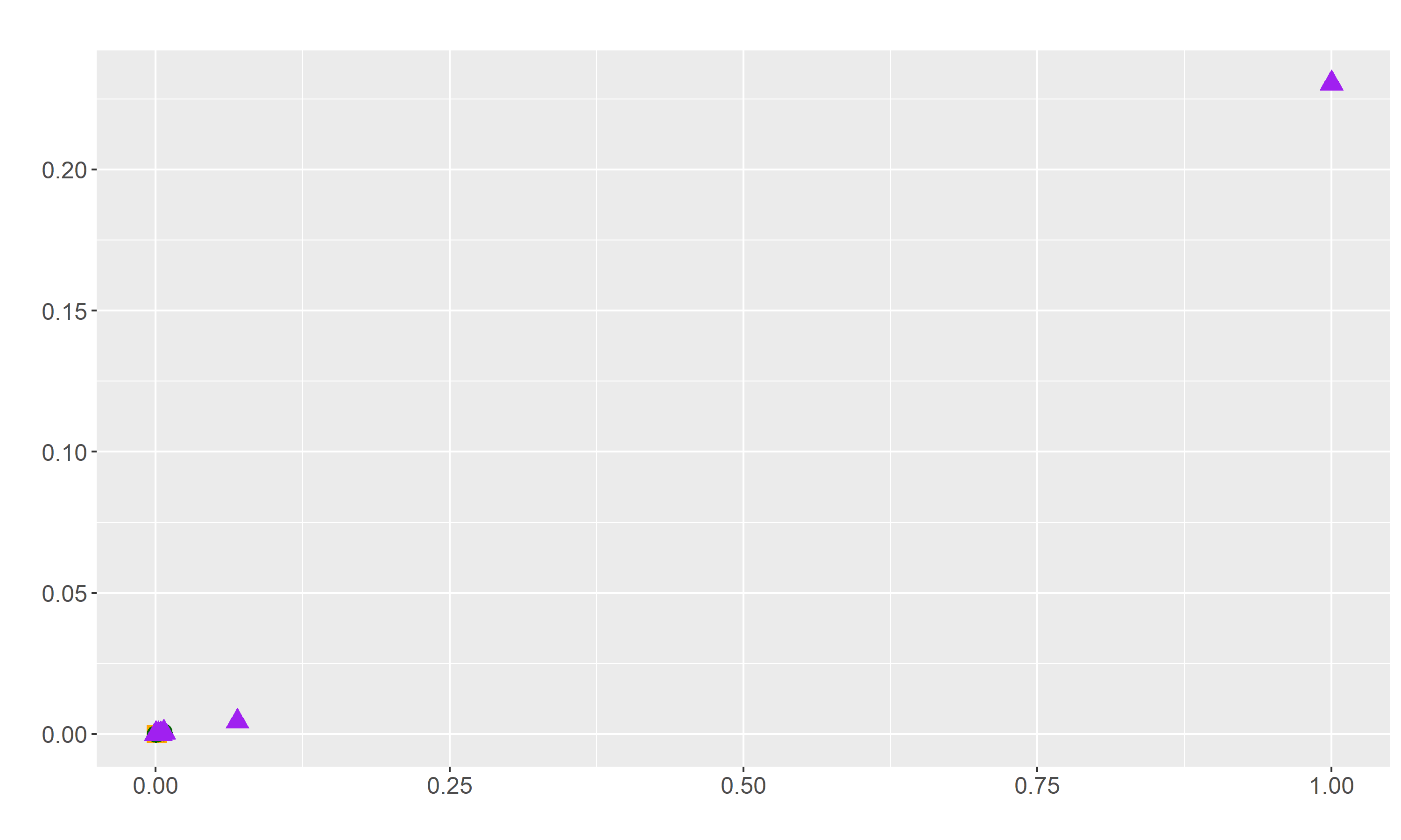}};
			\node(s2) at (2,0){\includegraphics[width=1cm]{figs/leyend_global.png}};
			\node[rotate=90](y) at (-3.4,0){{\color{darkgray}\tiny{$\log(\text{error}+1)$}}};
			\node(x) at (0,-2){{\color{darkgray}\tiny{distance}}};
		\end{tikzpicture}
		\caption{Global solvers: Distance to ref. solution vs error when ``solved''.}
\end{subfigure}
\caption{Results for problem \harmonic.}
\label{fig:harmonic}
\end{figure}

\newpage 

\subsection{\daisyF}
\begin{figure}[!htbp]
\centering
\begin{subfigure}{0.495\textwidth}
		\centering
		\begin{tikzpicture}
			\node(s1) at (0,0){\includegraphics[width=0.9\textwidth]{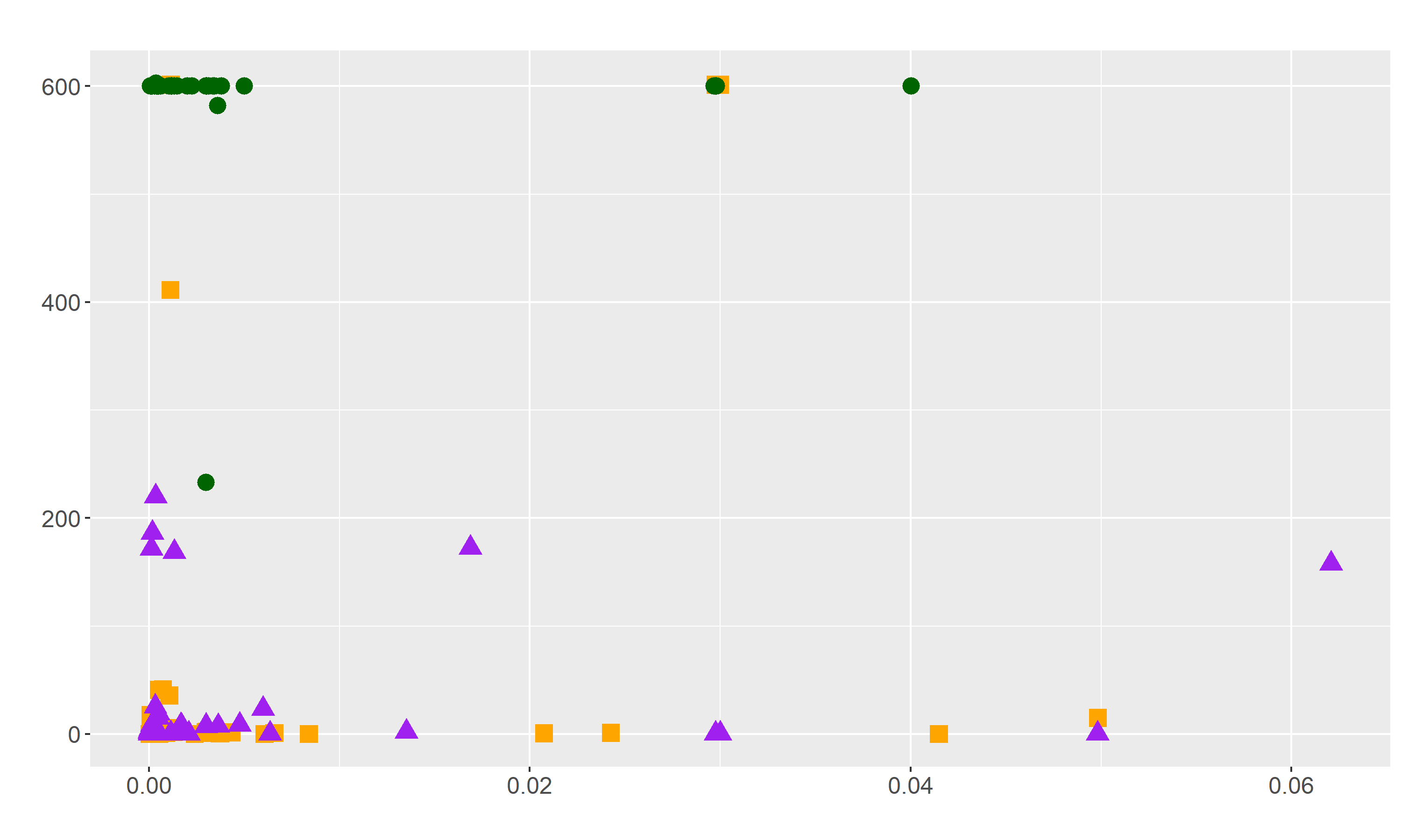}};
			\node[rotate=90](y) at (-3.35,0){{\color{darkgray}\tiny{time}}};
			\node(x) at (0,-2){{\color{darkgray}\tiny{distance}}};
		\end{tikzpicture}
		\caption{Global solvers: Distance to ref. solution vs time.}
\end{subfigure}
\begin{subfigure}{0.495\textwidth}
		\centering
		\begin{tikzpicture}
			\node(s1) at (0,0){\includegraphics[width=0.9\textwidth]{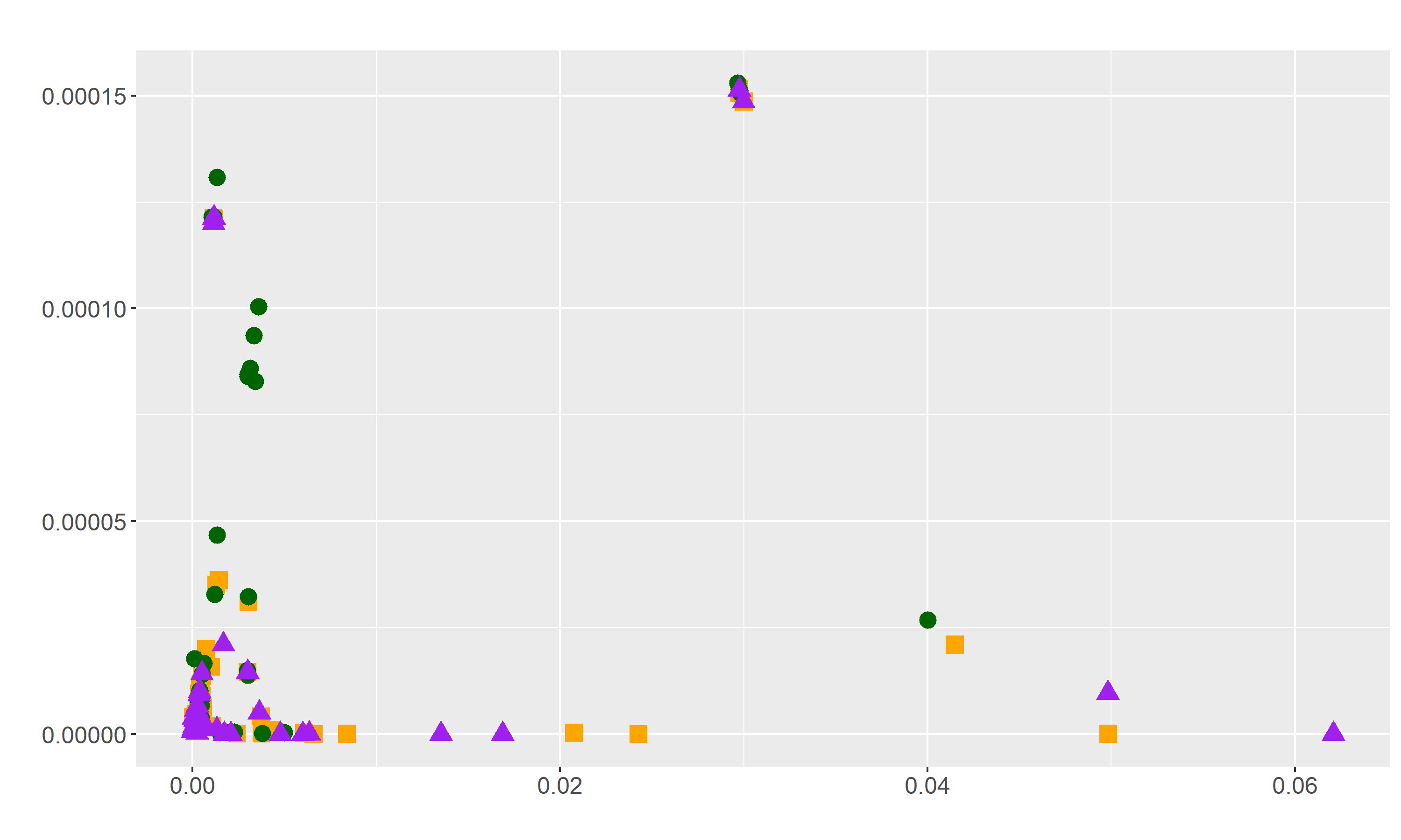}};
		\node(s2) at (2,0){\includegraphics[width=1cm]{figs/leyend_global.png}};
		\node[rotate=90](y) at (-3.4,0){{\color{darkgray}\tiny{$\log(\text{error}+1)$}}};
		\node(x) at (0,-2){{\color{darkgray}\tiny{distance}}};
		\end{tikzpicture}
		\caption{Global solvers: Distance to ref. solution vs error.}
\end{subfigure}

\vspace{0.2cm}

\begin{subfigure}{0.495\textwidth}
		\centering
		\begin{tikzpicture}
			\node(s1) at (0,0){\includegraphics[width=0.9\textwidth]{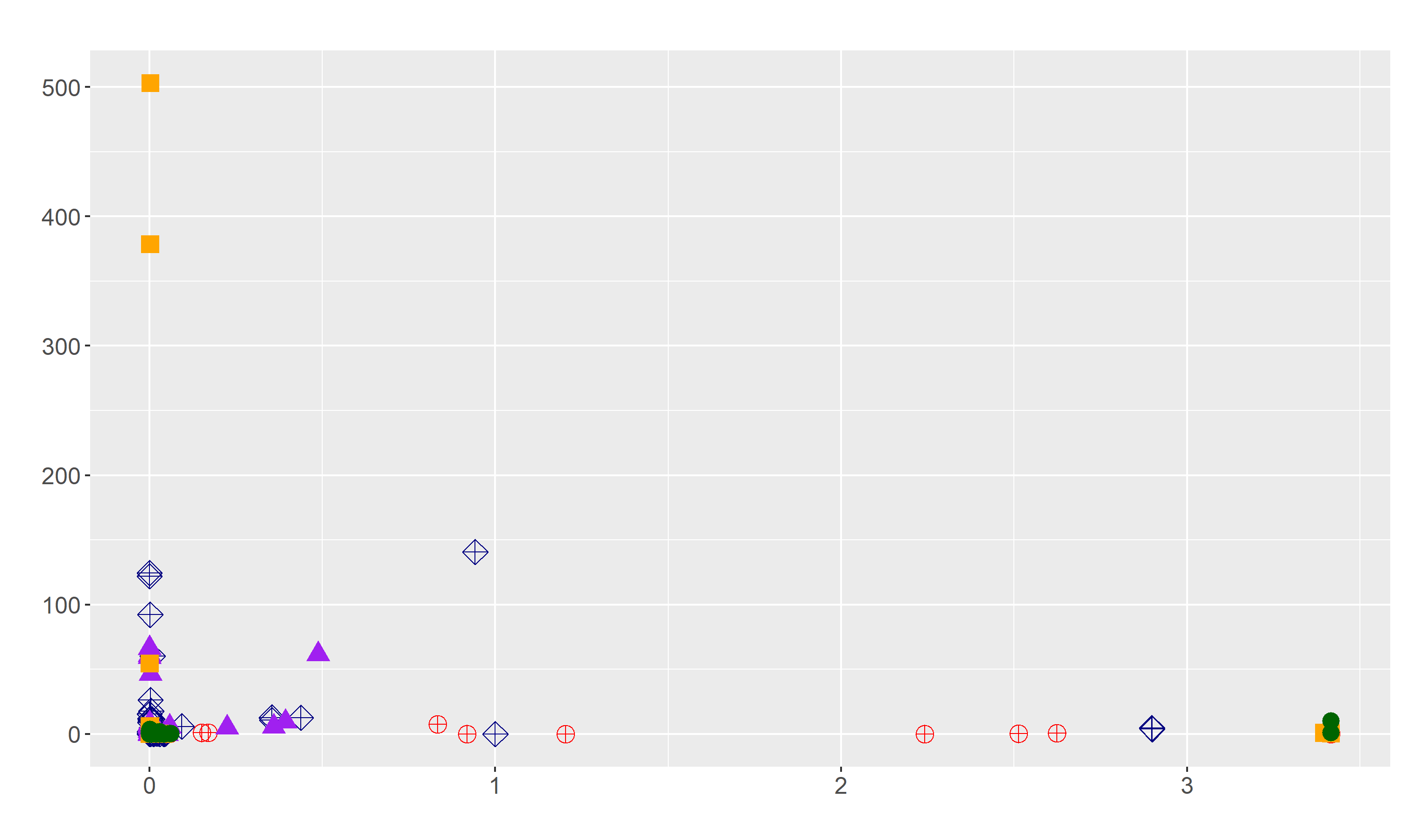}};
			\node[rotate=90](y) at (-3.35,0){{\color{darkgray}\tiny{time}}};
			\node(x) at (0,-2){{\color{darkgray}\tiny{distance}}};
		\end{tikzpicture}
		\caption{Local solvers: Distance to ref. solution vs time.}
\end{subfigure}
\begin{subfigure}{0.495\textwidth}
		\centering
		\begin{tikzpicture}
			\node(s1) at (0,0){\includegraphics[width=0.9\textwidth]{figs/daisy_mamil3F_distancia_erro_local.png}};
			\node(s2) at (2,0){\includegraphics[width=1cm]{figs/leyend_local.png}};
			\node[rotate=90](y) at (-3.4,0){{\color{darkgray}\tiny{$\log(\text{error}+1)$}}};
			\node(x) at (0,-2){{\color{darkgray}\tiny{distance}}};
		\end{tikzpicture}
		\caption{Local solvers: Distance to ref. solution vs error.}
\end{subfigure}

\begin{subfigure}{0.495\textwidth}
		\centering
		\begin{tikzpicture}
			\node(s1) at (0,0){\includegraphics[width=0.9\textwidth]{figs/daisy_mamil3F_distancia_erro_global_status0.png}};
			\node(s2) at (2,0){\includegraphics[width=1cm]{figs/leyend_global.png}};
			\node[rotate=90](y) at (-3.4,0){{\color{darkgray}\tiny{$\log(\text{error}+1)$}}};
			\node(x) at (0,-2){{\color{darkgray}\tiny{distance}}};
		\end{tikzpicture}
		\caption{Global solvers: Distance to ref. solution vs error when ``solved''.}
\end{subfigure}
\caption{Results for problem \daisyF.}
\label{fig:daisyF}
\end{figure}

\newpage 
\subsection{\daisyP}
\begin{figure}[!htbp]
\centering
\begin{subfigure}{0.495\textwidth}
		\centering
		\begin{tikzpicture}
			\node(s1) at (0,0){\includegraphics[width=0.9\textwidth]{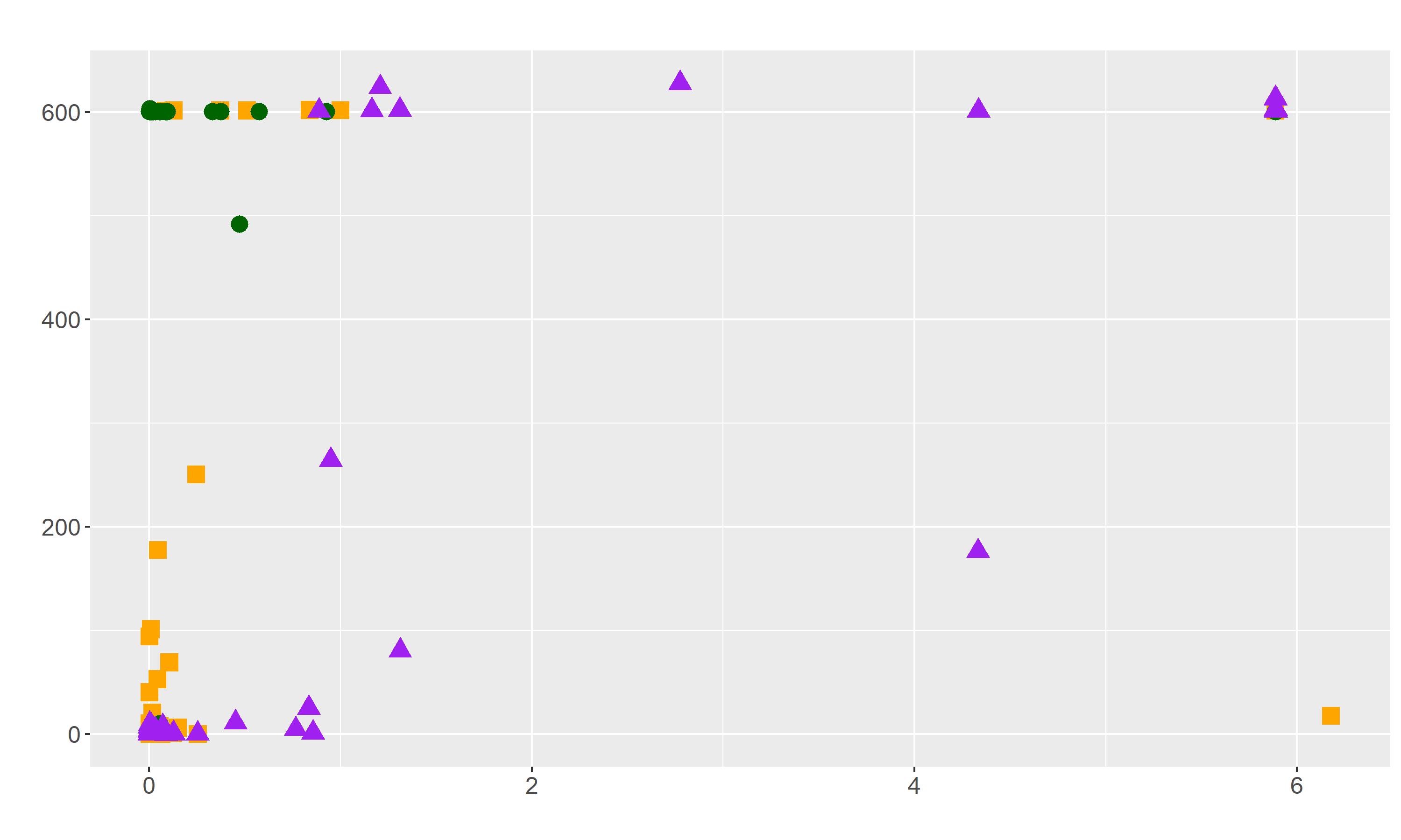}};
			\node[rotate=90](y) at (-3.35,0){{\color{darkgray}\tiny{time}}};
			\node(x) at (0,-2){{\color{darkgray}\tiny{distance}}};
		\end{tikzpicture}
		\caption{Global solvers: Distance to ref. solution vs time.}
\end{subfigure}
\begin{subfigure}{0.495\textwidth}
		\centering
		\begin{tikzpicture}
			\node(s1) at (0,0){\includegraphics[width=0.9\textwidth]{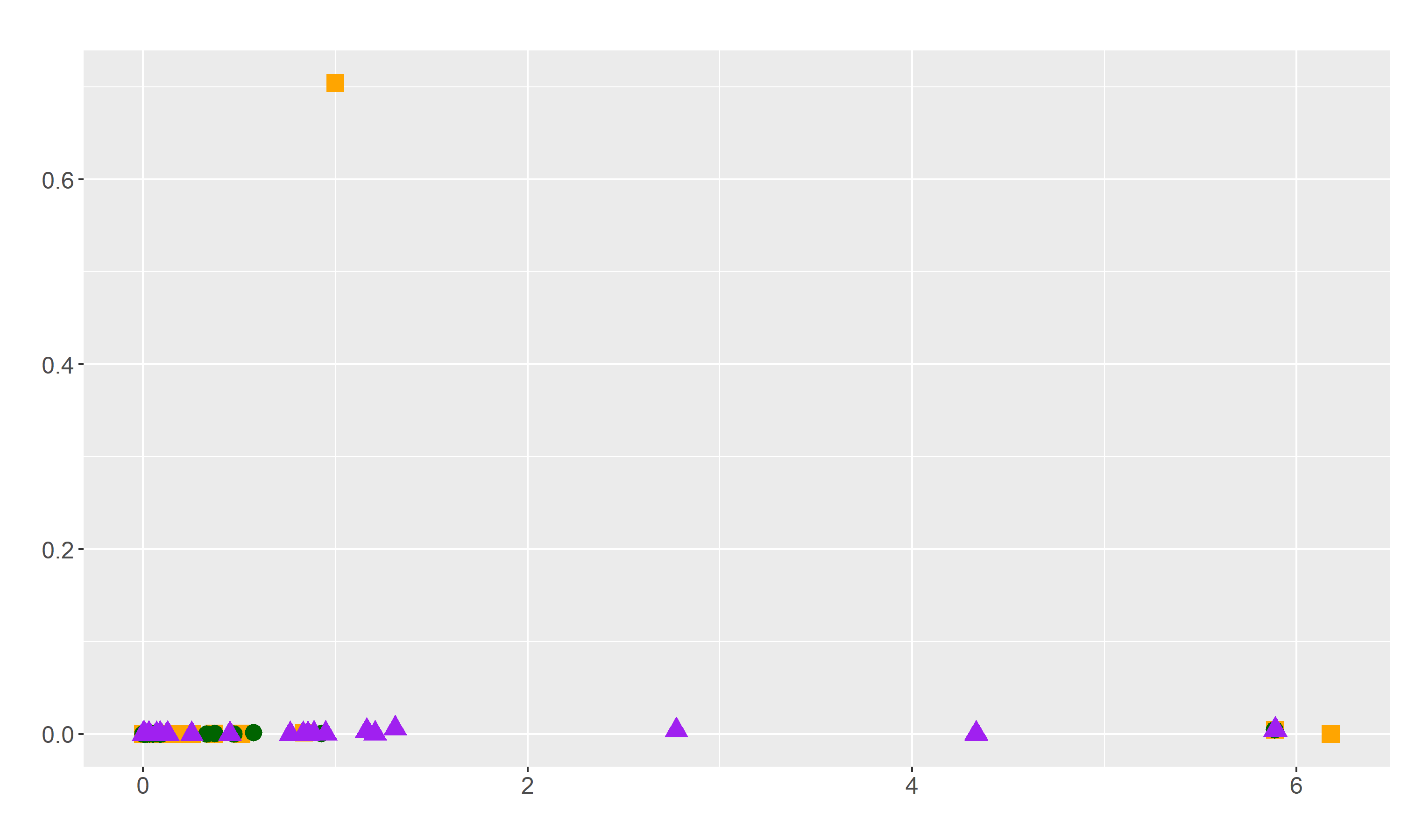}};
		\node(s2) at (2,0){\includegraphics[width=1cm]{figs/leyend_global.png}};
		\node[rotate=90](y) at (-3.4,0){{\color{darkgray}\tiny{$\log(\text{error}+1)$}}};
		\node(x) at (0,-2){{\color{darkgray}\tiny{distance}}};
		\end{tikzpicture}
		\caption{Global solvers: Distance to ref. solution vs error.}
\end{subfigure}

\vspace{0.2cm}

\begin{subfigure}{0.495\textwidth}
		\centering
		\begin{tikzpicture}
			\node(s1) at (0,0){\includegraphics[width=0.9\textwidth]{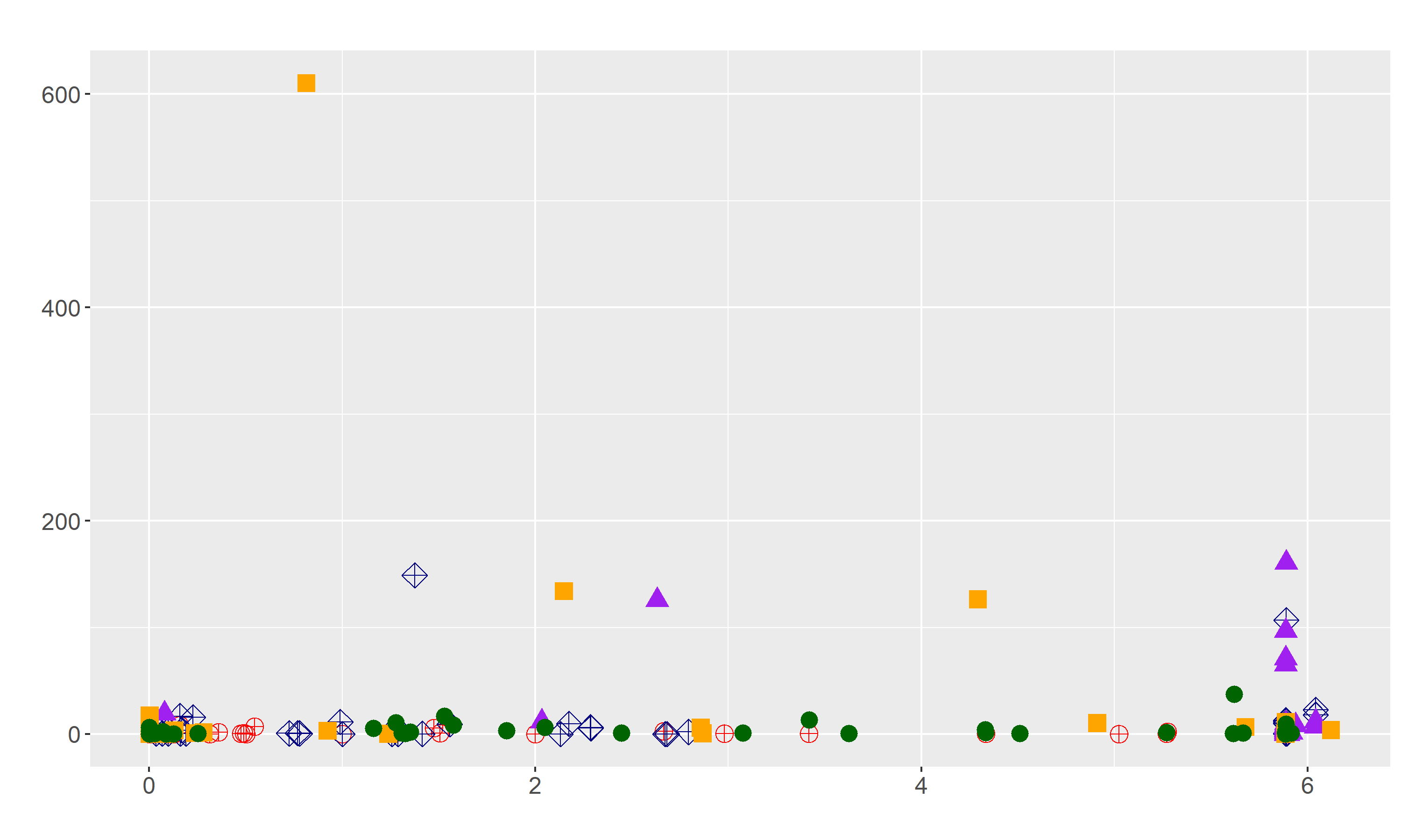}};
			\node[rotate=90](y) at (-3.35,0){{\color{darkgray}\tiny{time}}};
			\node(x) at (0,-2){{\color{darkgray}\tiny{distance}}};
		\end{tikzpicture}
		\caption{Local solvers: Distance to ref. solution vs time.}
\end{subfigure}
\begin{subfigure}{0.495\textwidth}
		\centering
		\begin{tikzpicture}
			\node(s1) at (0,0){\includegraphics[width=0.9\textwidth]{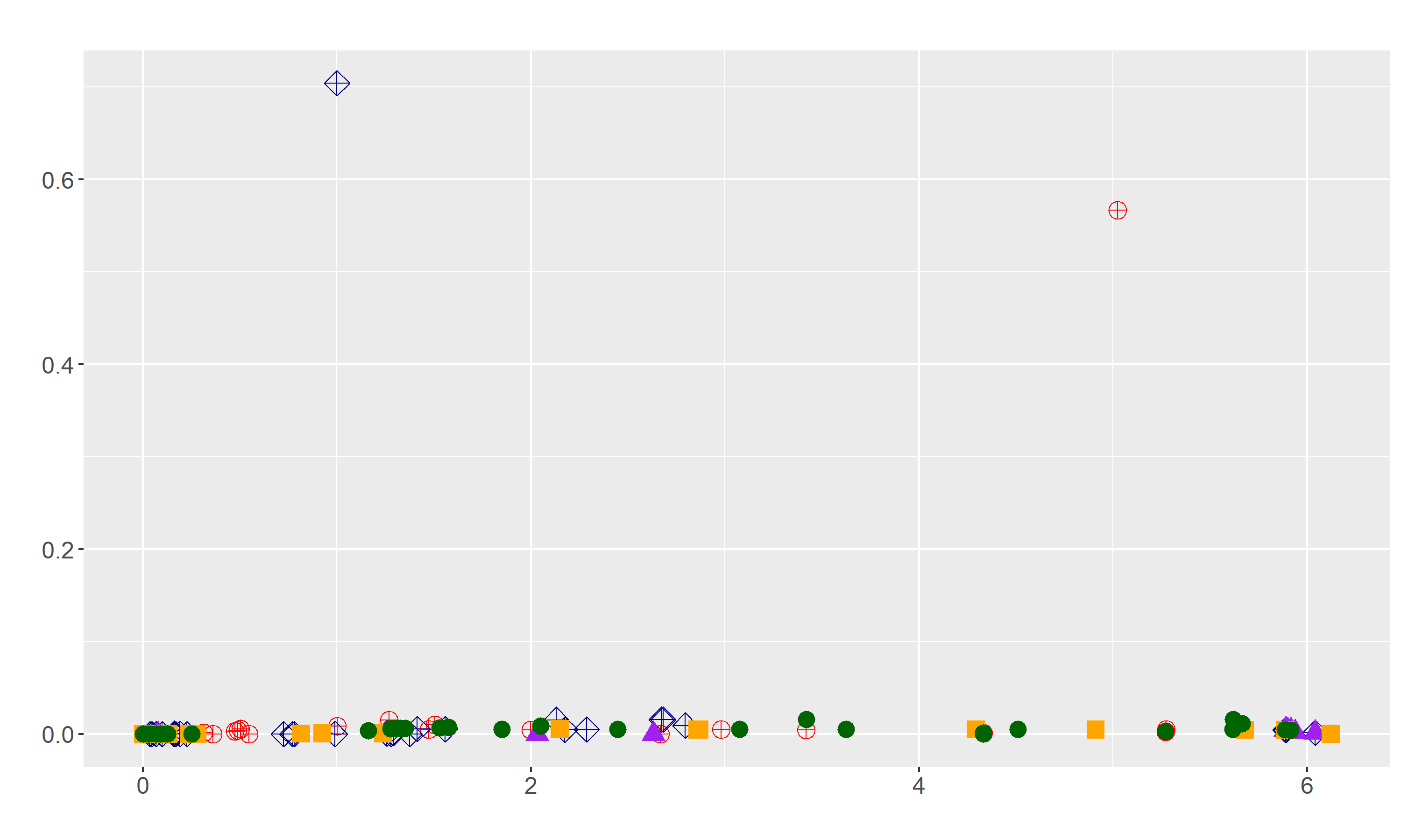}};
			\node(s2) at (2,0){\includegraphics[width=1cm]{figs/leyend_local.png}};
			\node[rotate=90](y) at (-3.4,0){{\color{darkgray}\tiny{$\log(\text{error}+1)$}}};
			\node(x) at (0,-2){{\color{darkgray}\tiny{distance}}};
		\end{tikzpicture}
		\caption{Local solvers: Distance to ref. solution vs error.}
\end{subfigure}

\begin{subfigure}{0.495\textwidth}
		\centering
		\begin{tikzpicture}
			\node(s1) at (0,0){\includegraphics[width=0.9\textwidth]{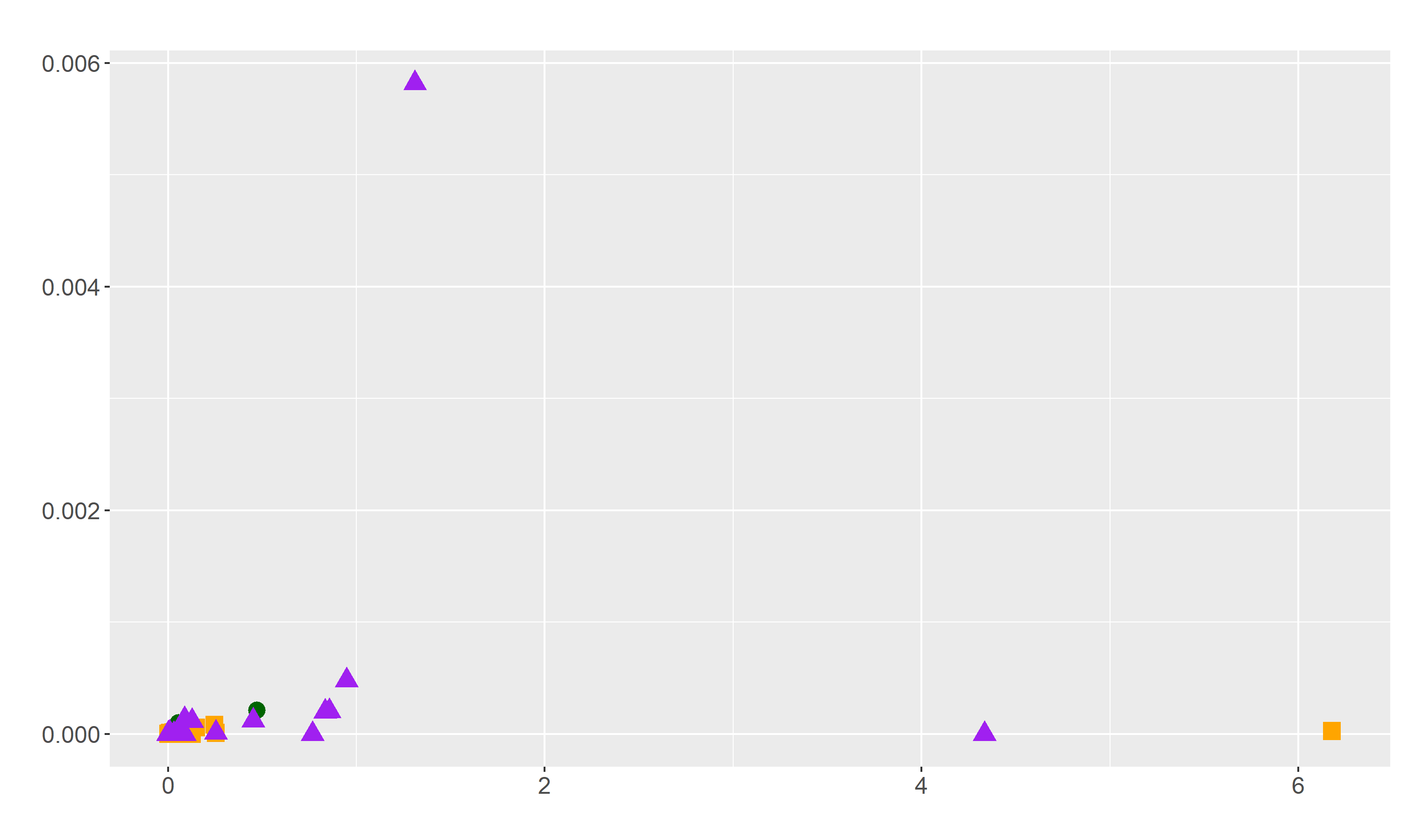}};
			\node(s2) at (2,0){\includegraphics[width=1cm]{figs/leyend_global.png}};
			\node[rotate=90](y) at (-3.4,0){{\color{darkgray}\tiny{$\log(\text{error}+1)$}}};
			\node(x) at (0,-2){{\color{darkgray}\tiny{distance}}};
		\end{tikzpicture}
		\caption{Global solvers: Distance to ref. solution vs error when ``solved''.}
\end{subfigure}
\caption{Results for problem \daisyP.}
\label{fig:daisyP}
\end{figure}

\newpage 
\subsection{\hivF}
\begin{figure}[!htbp]
\centering
\begin{subfigure}{0.495\textwidth}
		\centering
		\begin{tikzpicture}
			\node(s1) at (0,0){\includegraphics[width=0.9\textwidth]{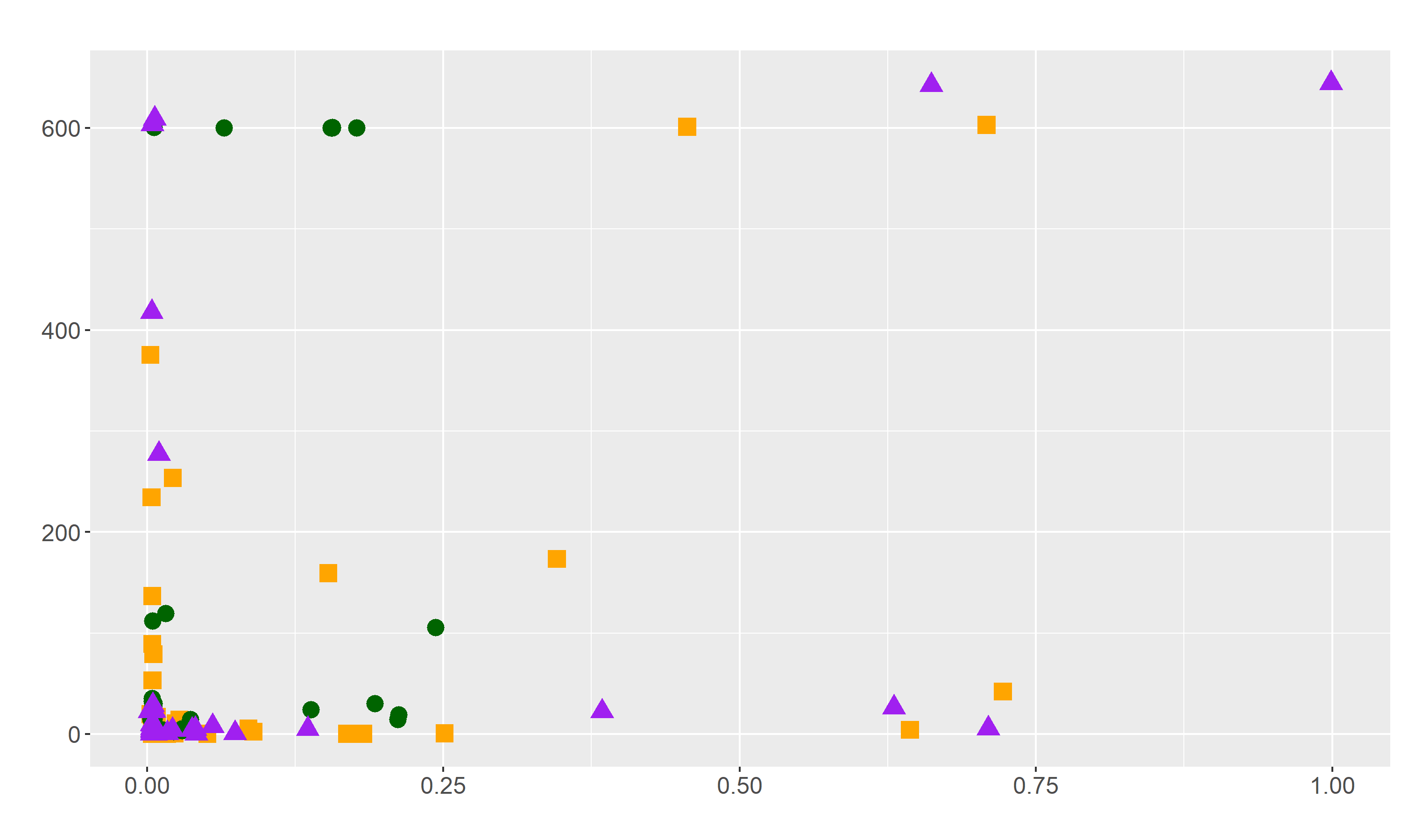}};
			\node[rotate=90](y) at (-3.35,0){{\color{darkgray}\tiny{time}}};
			\node(x) at (0,-2){{\color{darkgray}\tiny{distance}}};
		\end{tikzpicture}
		\caption{Global solvers: Distance to ref. solution vs time.}
\end{subfigure}
\begin{subfigure}{0.495\textwidth}
		\centering
		\begin{tikzpicture}
			\node(s1) at (0,0){\includegraphics[width=0.9\textwidth]{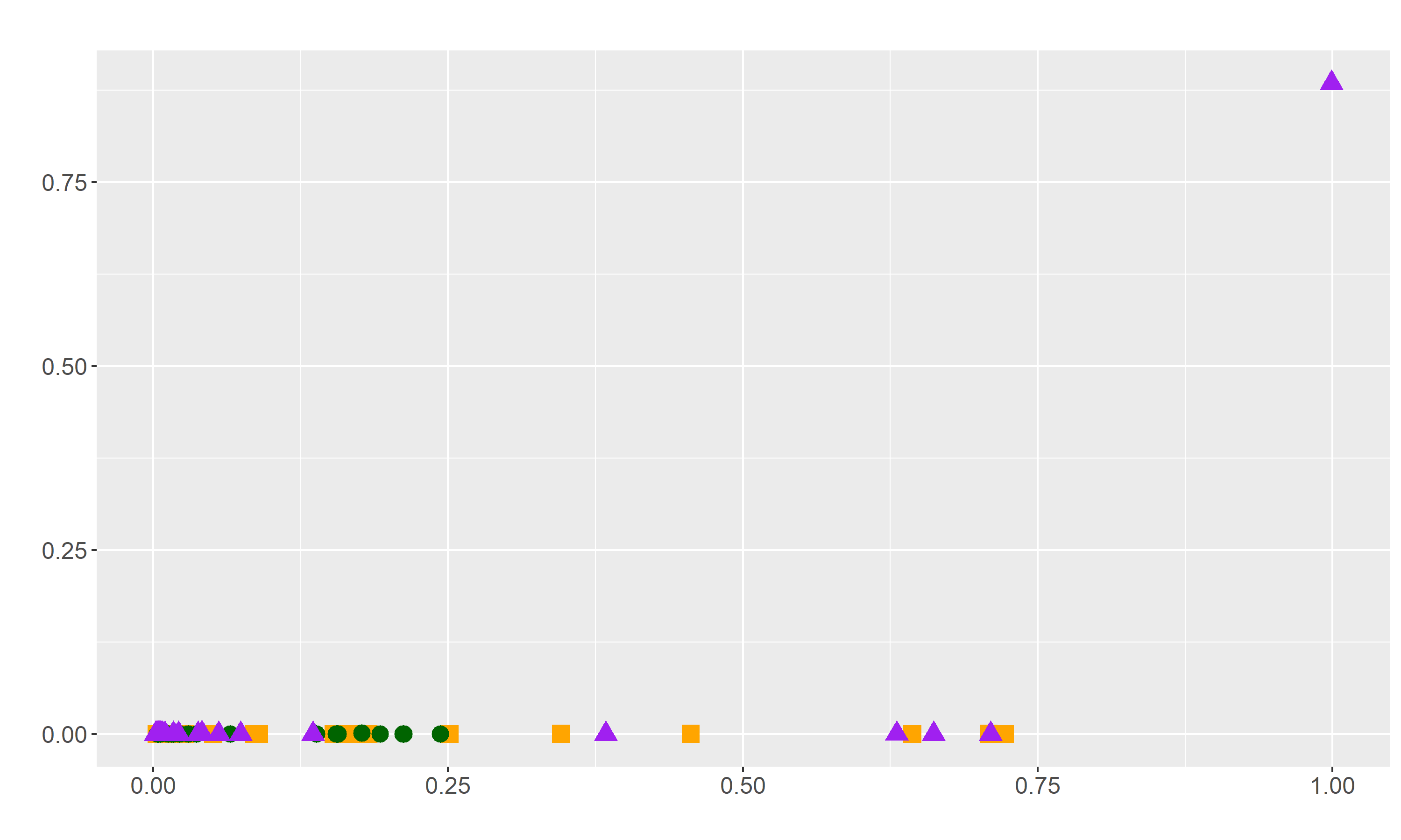}};
		\node(s2) at (2,0){\includegraphics[width=1cm]{figs/leyend_global.png}};
		\node[rotate=90](y) at (-3.4,0){{\color{darkgray}\tiny{$\log(\text{error}+1)$}}};
		\node(x) at (0,-2){{\color{darkgray}\tiny{distance}}};
		\end{tikzpicture}
		\caption{Global solvers: Distance to ref. solution vs error.}
\end{subfigure}

\vspace{0.2cm}

\begin{subfigure}{0.495\textwidth}
		\centering
		\begin{tikzpicture}
			\node(s1) at (0,0){\includegraphics[width=0.9\textwidth]{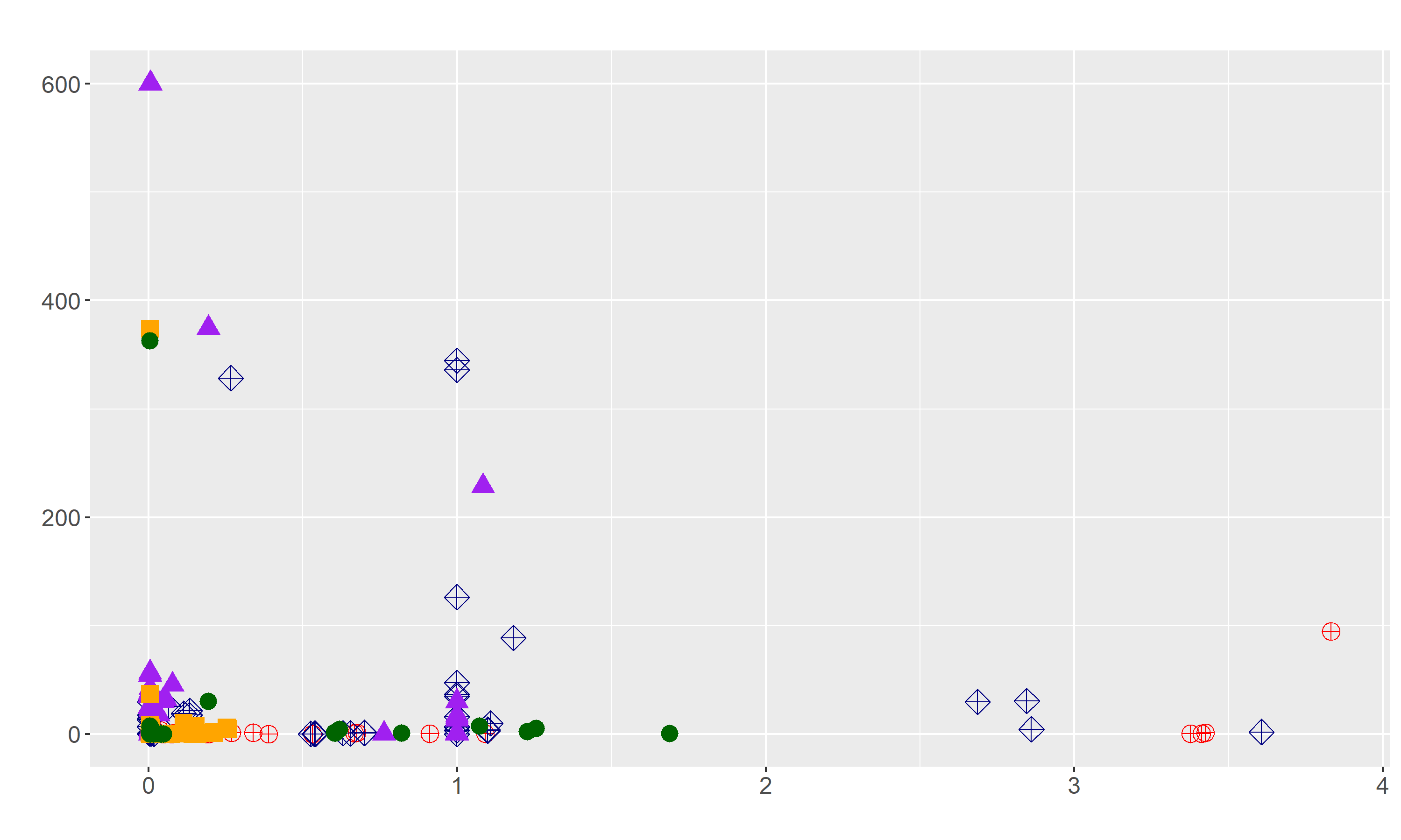}};
			\node[rotate=90](y) at (-3.35,0){{\color{darkgray}\tiny{time}}};
			\node(x) at (0,-2){{\color{darkgray}\tiny{distance}}};
		\end{tikzpicture}
		\caption{Local solvers: Distance to ref. solution vs time.}
\end{subfigure}
\begin{subfigure}{0.495\textwidth}
		\centering
		\begin{tikzpicture}
			\node(s1) at (0,0){\includegraphics[width=0.9\textwidth]{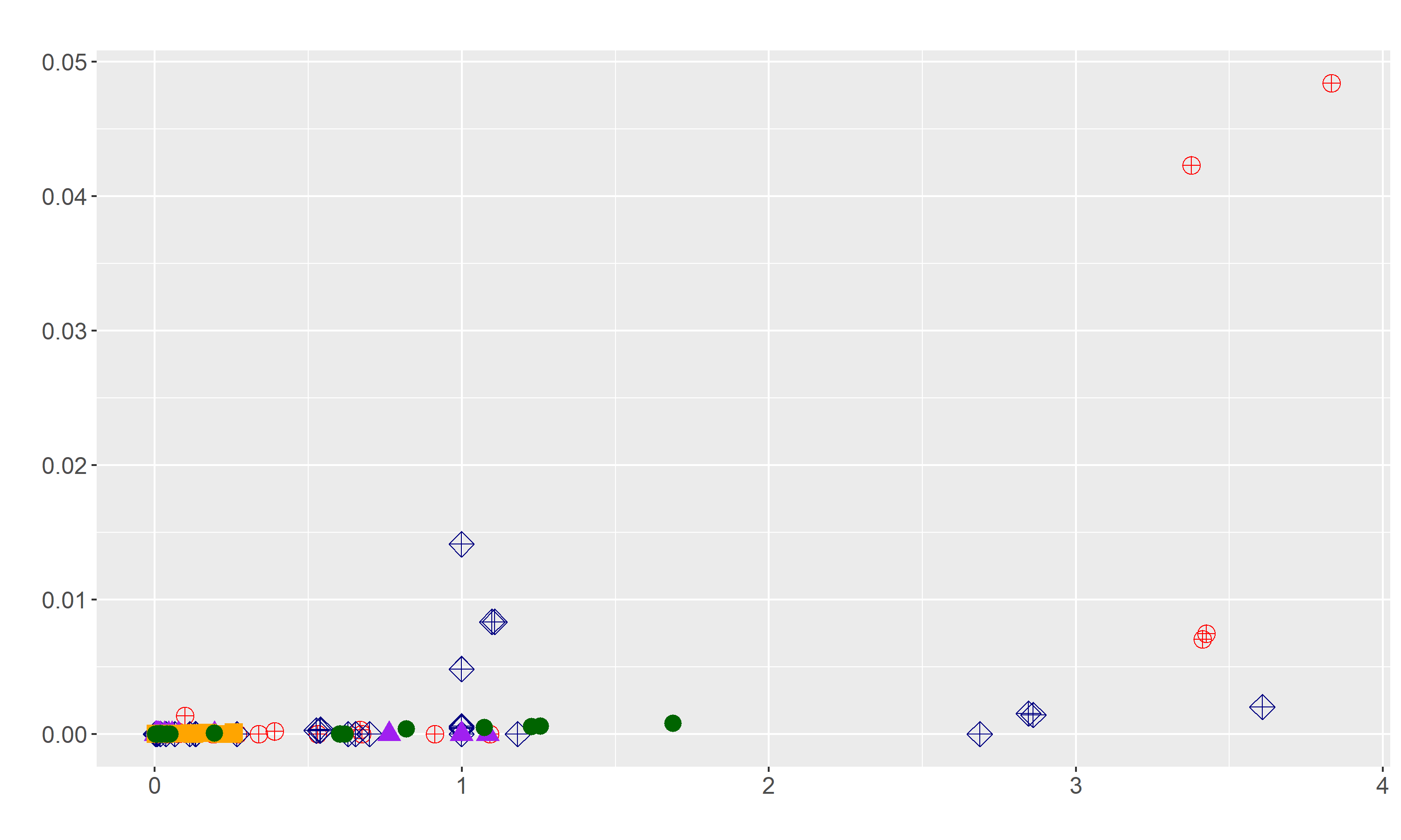}};
			\node(s2) at (2,0){\includegraphics[width=1cm]{figs/leyend_local.png}};
			\node[rotate=90](y) at (-3.4,0){{\color{darkgray}\tiny{$\log(\text{error}+1)$}}};
			\node(x) at (0,-2){{\color{darkgray}\tiny{distance}}};
		\end{tikzpicture}
		\caption{Local solvers: Distance to ref. solution vs error.}
\end{subfigure}

\begin{subfigure}{0.495\textwidth}
		\centering
		\begin{tikzpicture}
			\node(s1) at (0,0){\includegraphics[width=0.9\textwidth]{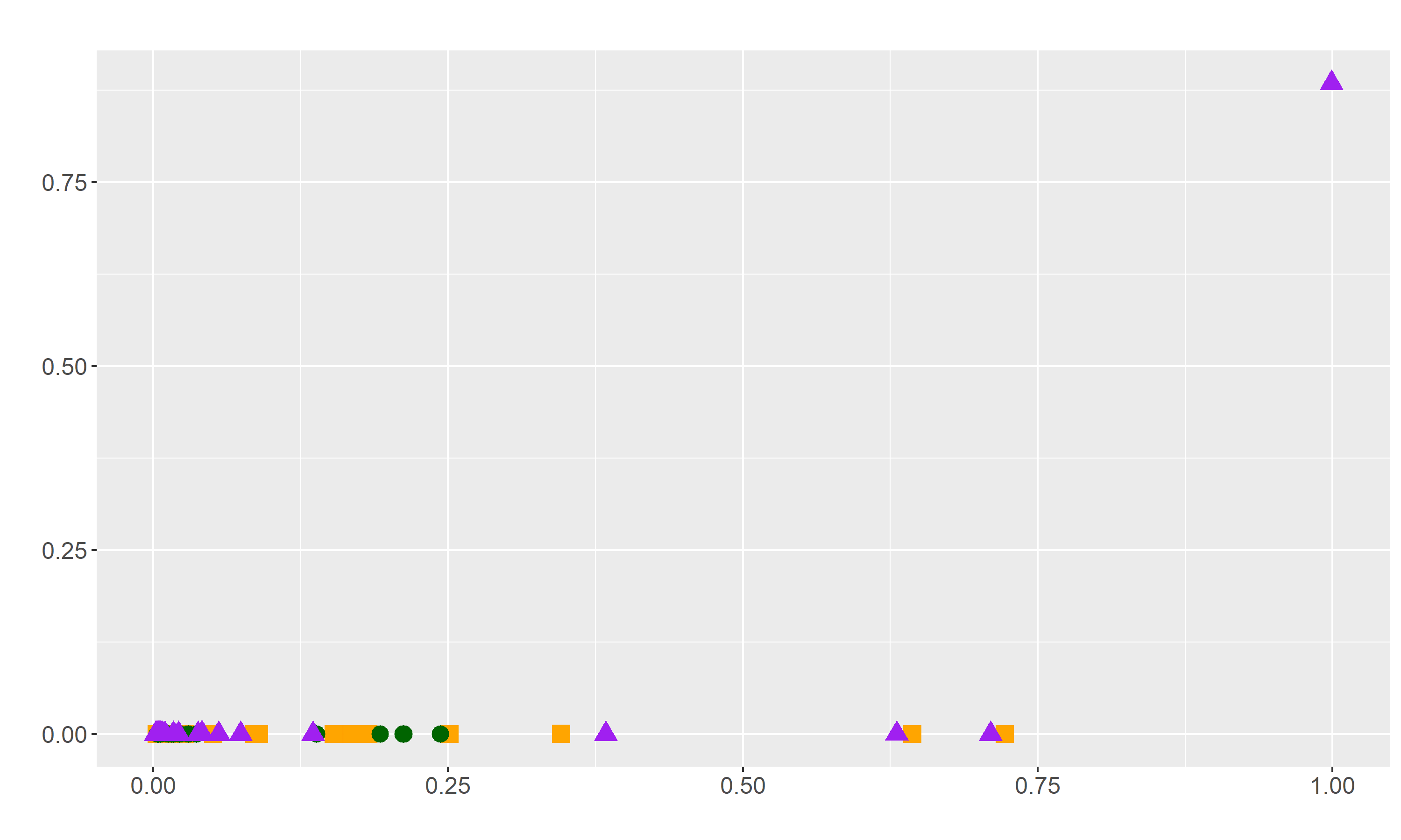}};
			\node(s2) at (2,0){\includegraphics[width=1cm]{figs/leyend_global.png}};
			\node[rotate=90](y) at (-3.4,0){{\color{darkgray}\tiny{$\log(\text{error}+1)$}}};
			\node(x) at (0,-2){{\color{darkgray}\tiny{distance}}};
		\end{tikzpicture}
		\caption{Global solvers: Distance to ref. solution vs error when ``solved''.}
\end{subfigure}

\caption{Results for problem \hivF.}
\label{fig:hivF}
\end{figure}

\newpage 
\subsection{\hivP}
\begin{figure}[!htbp]
\centering
\begin{subfigure}{0.495\textwidth}
		\centering
		\begin{tikzpicture}
			\node(s1) at (0,0){\includegraphics[width=0.9\textwidth]{figs/hivP_distancia_tempo_global.png}};
			\node[rotate=90](y) at (-3.35,0){{\color{darkgray}\tiny{time}}};
			\node(x) at (0,-2){{\color{darkgray}\tiny{distance}}};
		\end{tikzpicture}
		\caption{Global solvers: Distance to ref. solution vs time.}
\end{subfigure}
\begin{subfigure}{0.495\textwidth}
		\centering
		\begin{tikzpicture}
			\node(s1) at (0,0){\includegraphics[width=0.9\textwidth]{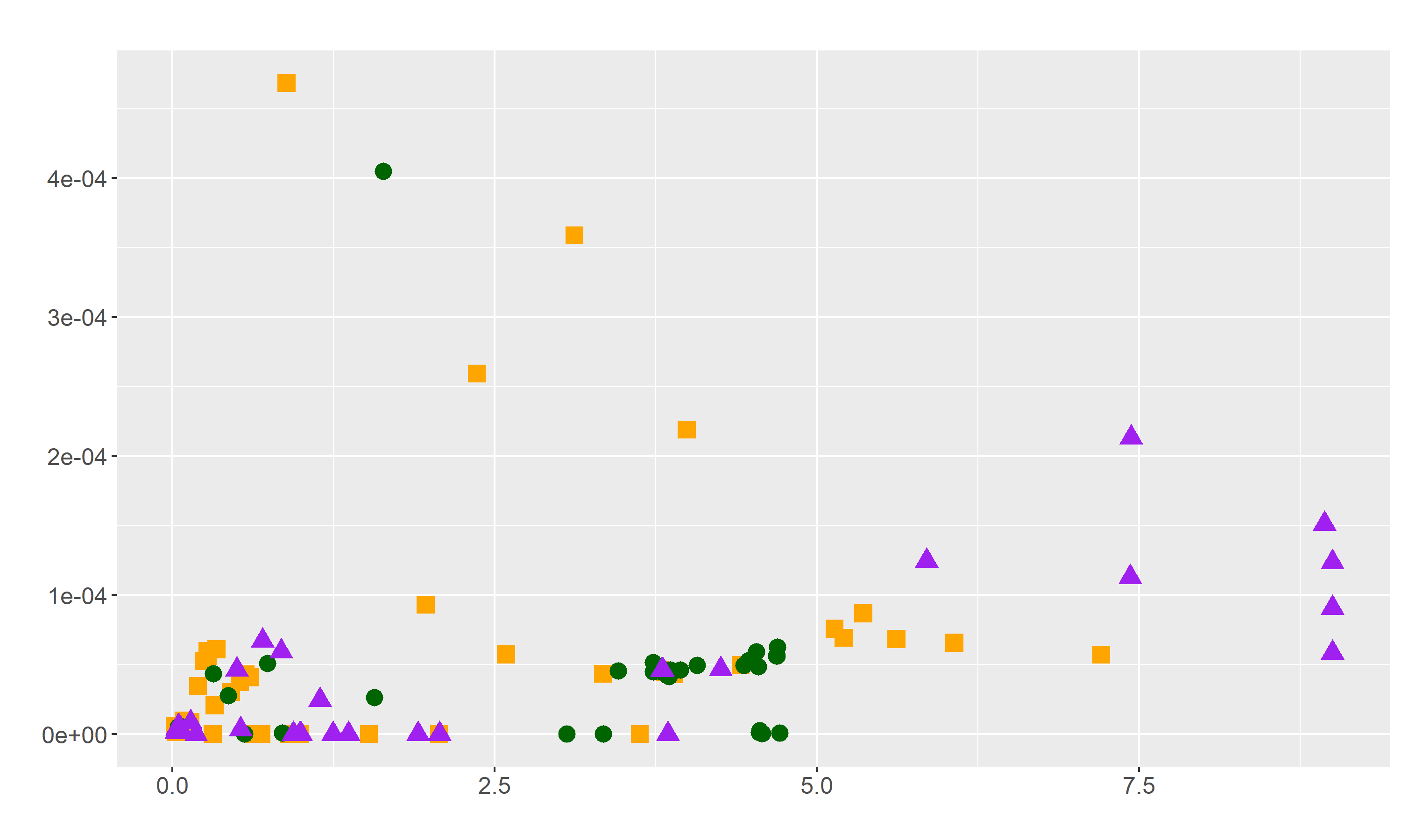}};
		\node(s2) at (2,0){\includegraphics[width=1cm]{figs/leyend_global.png}};
		\node[rotate=90](y) at (-3.4,0){{\color{darkgray}\tiny{$\log(\text{error}+1)$}}};
		\node(x) at (0,-2){{\color{darkgray}\tiny{distance}}};
		\end{tikzpicture}
		\caption{Global solvers: Distance to ref. solution vs error.}
\end{subfigure}

\vspace{0.2cm}

\begin{subfigure}{0.495\textwidth}
		\centering
		\begin{tikzpicture}
			\node(s1) at (0,0){\includegraphics[width=0.9\textwidth]{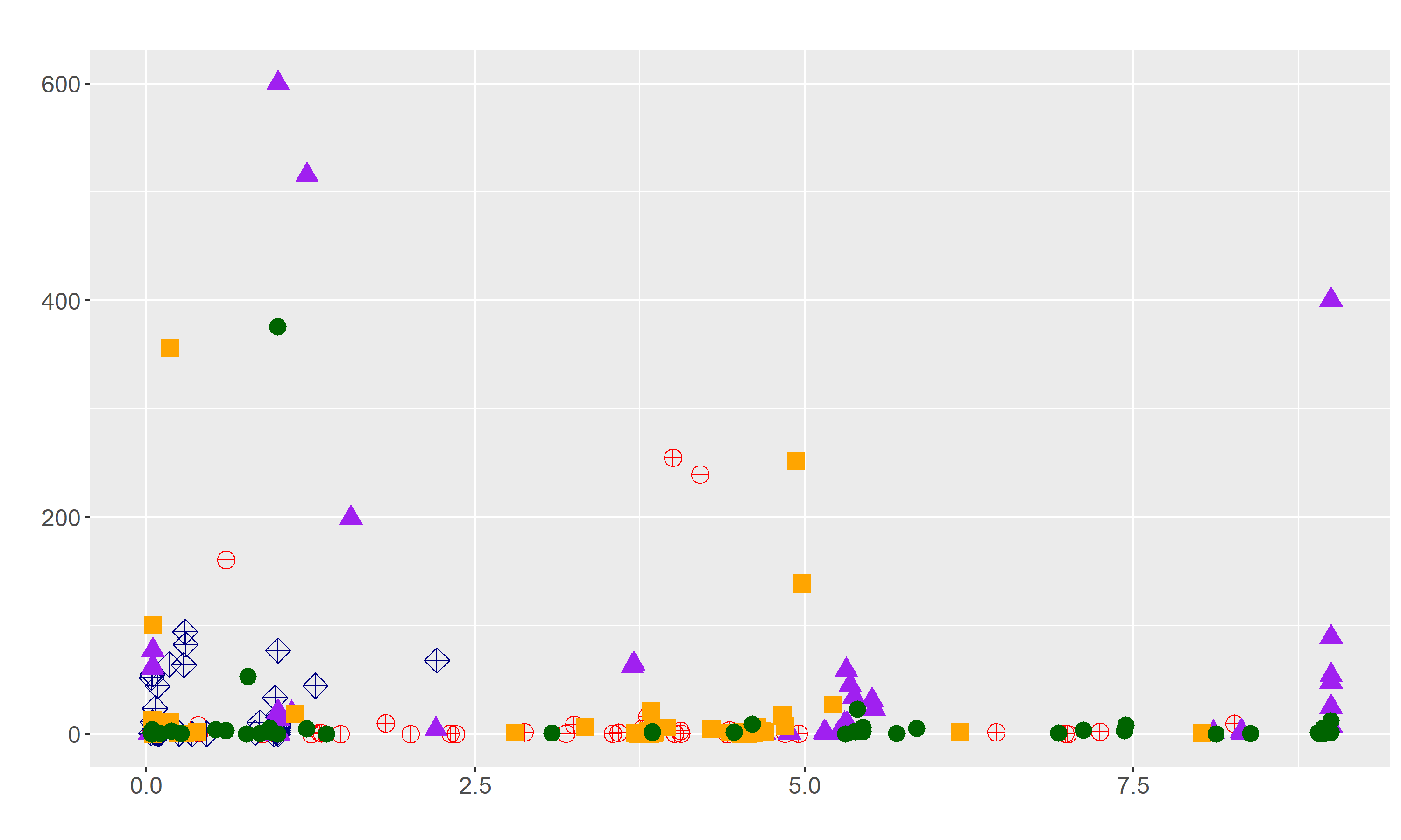}};
			\node[rotate=90](y) at (-3.35,0){{\color{darkgray}\tiny{time}}};
			\node(x) at (0,-2){{\color{darkgray}\tiny{distance}}};
		\end{tikzpicture}
		\caption{Local solvers: Distance to ref. solution vs time.}
\end{subfigure}
\begin{subfigure}{0.495\textwidth}
		\centering
		\begin{tikzpicture}
			\node(s1) at (0,0){\includegraphics[width=0.9\textwidth]{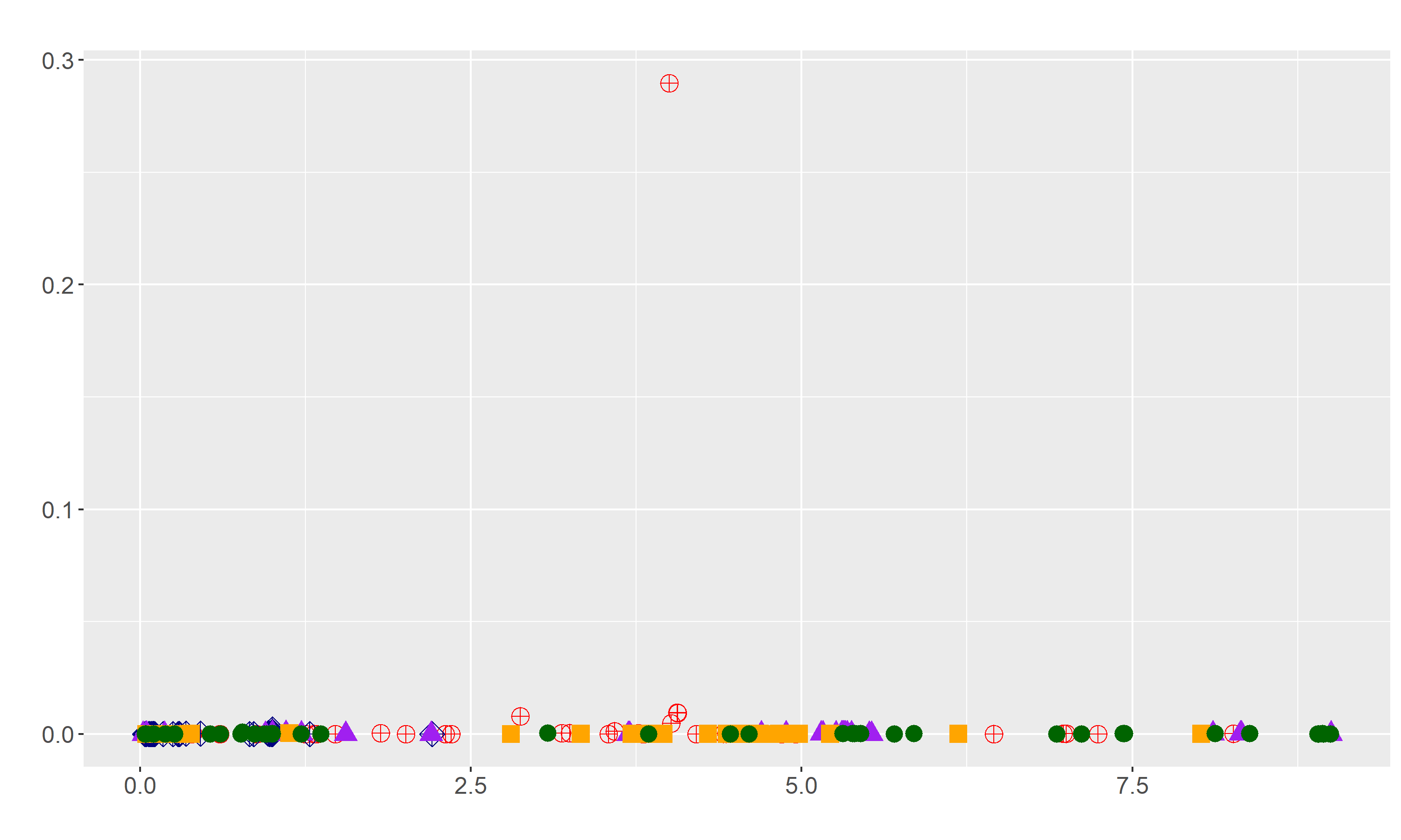}};
			\node(s2) at (2,0){\includegraphics[width=1cm]{figs/leyend_local.png}};
			\node[rotate=90](y) at (-3.4,0){{\color{darkgray}\tiny{$\log(\text{error}+1)$}}};
			\node(x) at (0,-2){{\color{darkgray}\tiny{distance}}};
		\end{tikzpicture}
		\caption{Local solvers: Distance to ref. solution vs error.}
\end{subfigure}

\begin{subfigure}{0.495\textwidth}
		\centering
		\begin{tikzpicture}
			\node(s1) at (0,0){\includegraphics[width=0.9\textwidth]{figs/hivP_distancia_erro_global_status0.png}};
			\node(s2) at (2,0){\includegraphics[width=1cm]{figs/leyend_global.png}};
			\node[rotate=90](y) at (-3.4,0){{\color{darkgray}\tiny{$\log(\text{error}+1)$}}};
			\node(x) at (0,-2){{\color{darkgray}\tiny{distance}}};
		\end{tikzpicture}
		\caption{Global solvers: Distance to ref. solution vs error when ``solved''.}
\end{subfigure}

\caption{Results for problem \hivP.}
\label{fig:hivP}
\end{figure}

\newpage 
\subsection{\lvF}
\begin{figure}[!htbp]
\centering
\begin{subfigure}{0.495\textwidth}
		\centering
		\begin{tikzpicture}
			\node(s1) at (0,0){\includegraphics[width=0.9\textwidth]{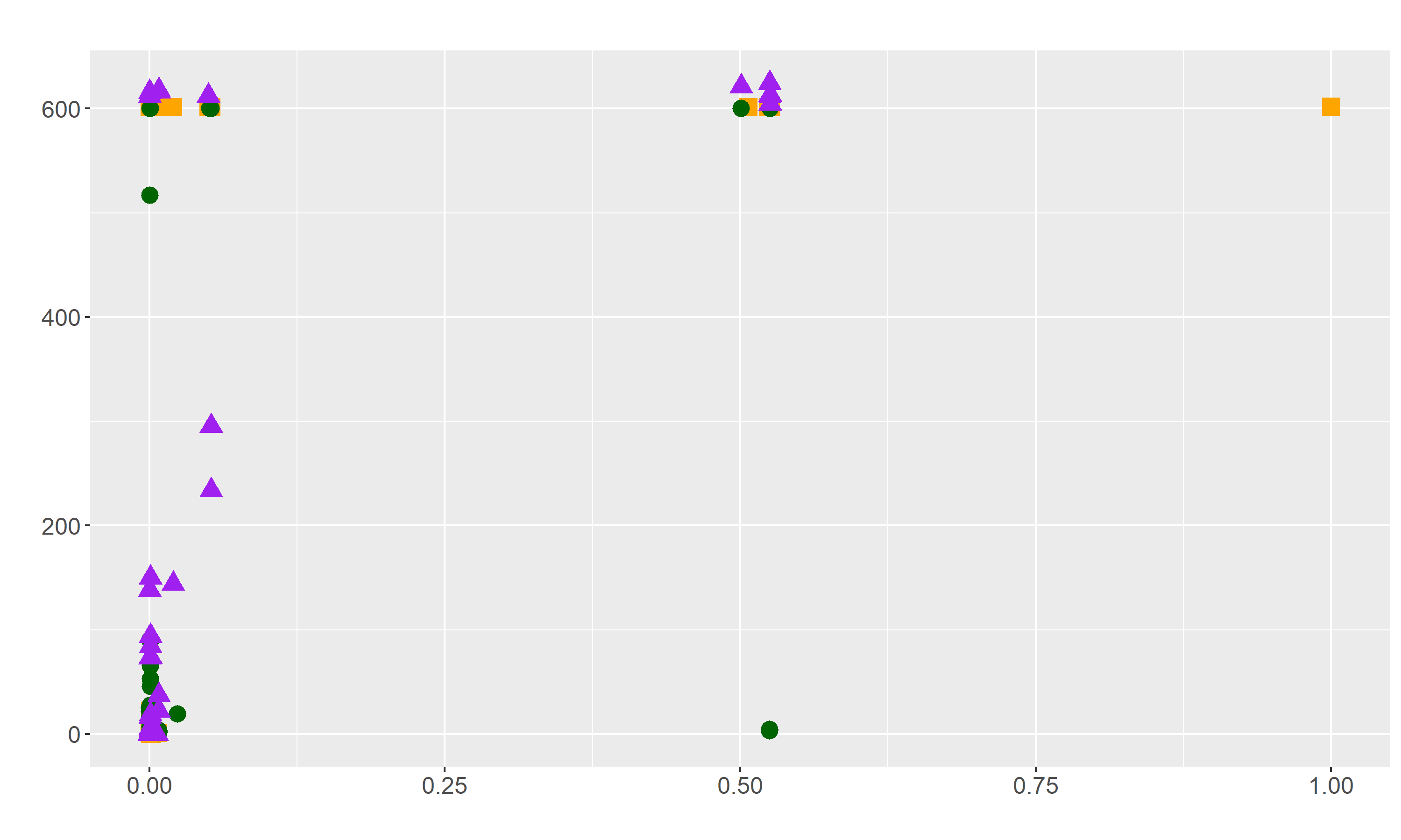}};
			\node[rotate=90](y) at (-3.35,0){{\color{darkgray}\tiny{time}}};
			\node(x) at (0,-2){{\color{darkgray}\tiny{distance}}};
		\end{tikzpicture}
		\caption{Global solvers: Distance to ref. solution vs time.}
\end{subfigure}
\begin{subfigure}{0.495\textwidth}
		\centering
		\begin{tikzpicture}
			\node(s1) at (0,0){\includegraphics[width=0.9\textwidth]{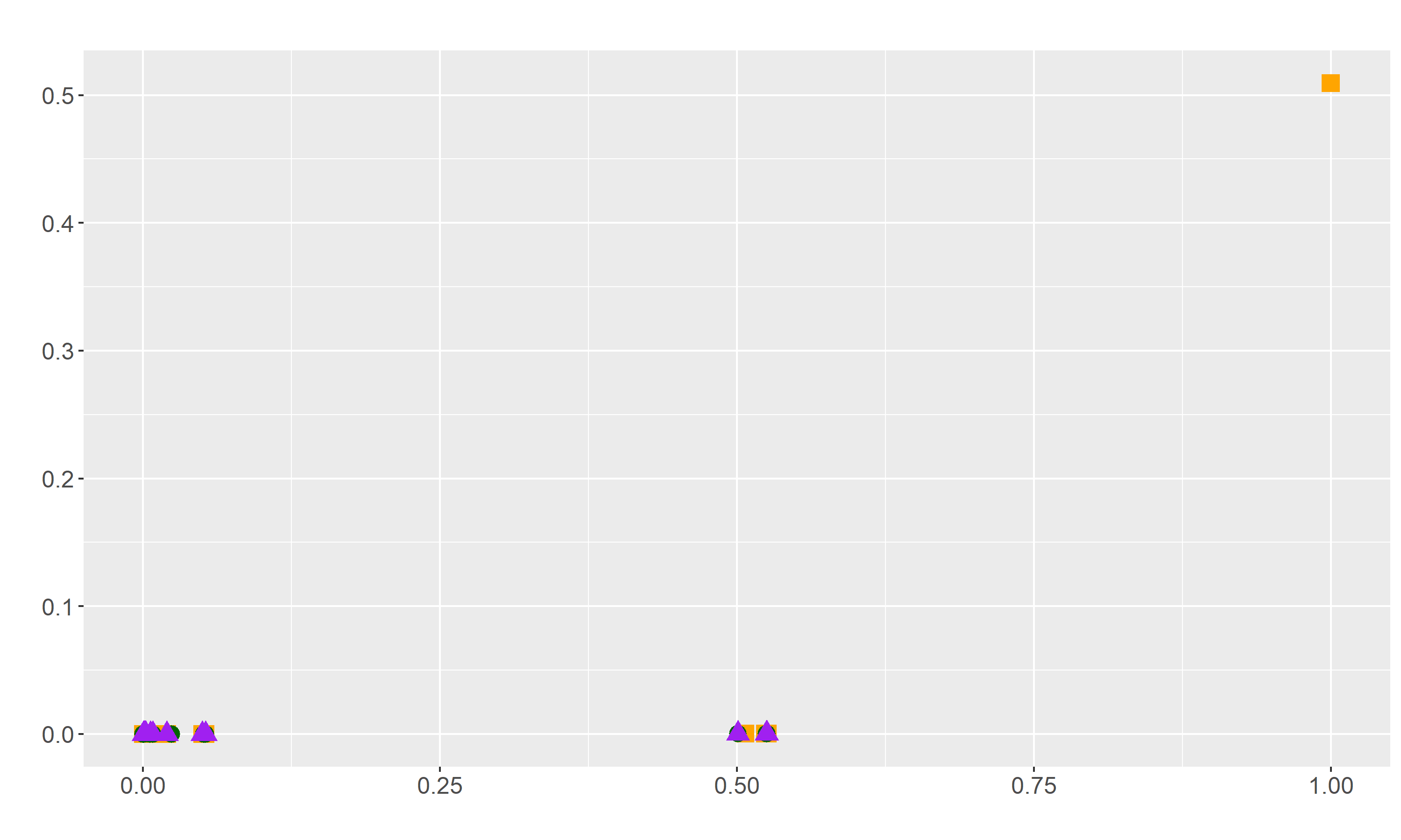}};
		\node(s2) at (2,0){\includegraphics[width=1cm]{figs/leyend_global.png}};
		\node[rotate=90](y) at (-3.4,0){{\color{darkgray}\tiny{$\log(\text{error}+1)$}}};
		\node(x) at (0,-2){{\color{darkgray}\tiny{distance}}};
		\end{tikzpicture}
		\caption{Global solvers: Distance to ref. solution vs error.}
\end{subfigure}

\vspace{0.2cm}

\begin{subfigure}{0.495\textwidth}
		\centering
		\begin{tikzpicture}
			\node(s1) at (0,0){\includegraphics[width=0.9\textwidth]{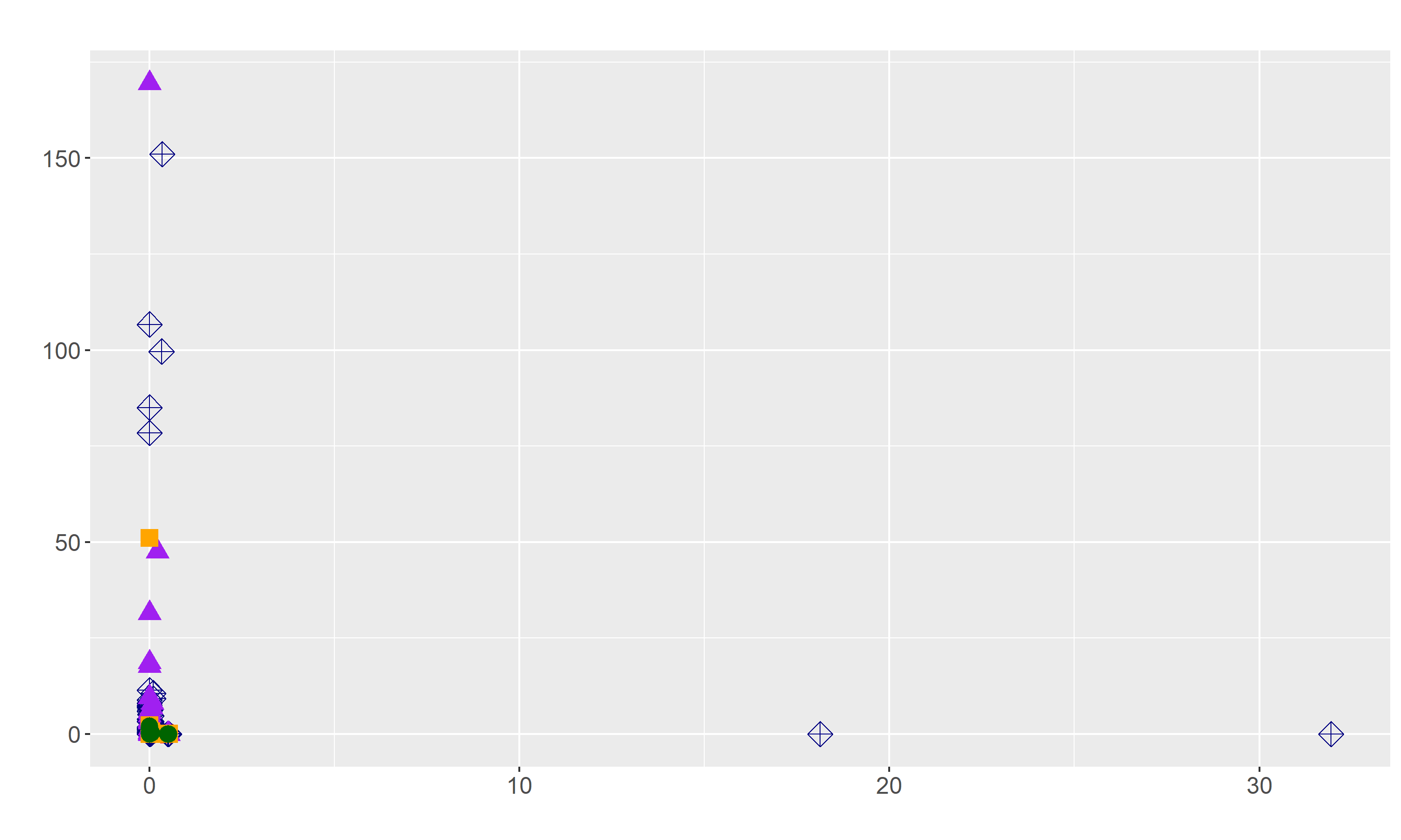}};
			\node[rotate=90](y) at (-3.35,0){{\color{darkgray}\tiny{time}}};
			\node(x) at (0,-2){{\color{darkgray}\tiny{distance}}};
		\end{tikzpicture}
		\caption{Local solvers: Distance to ref. solution vs time.}
\end{subfigure}
\begin{subfigure}{0.495\textwidth}
		\centering
		\begin{tikzpicture}
			\node(s1) at (0,0){\includegraphics[width=0.9\textwidth]{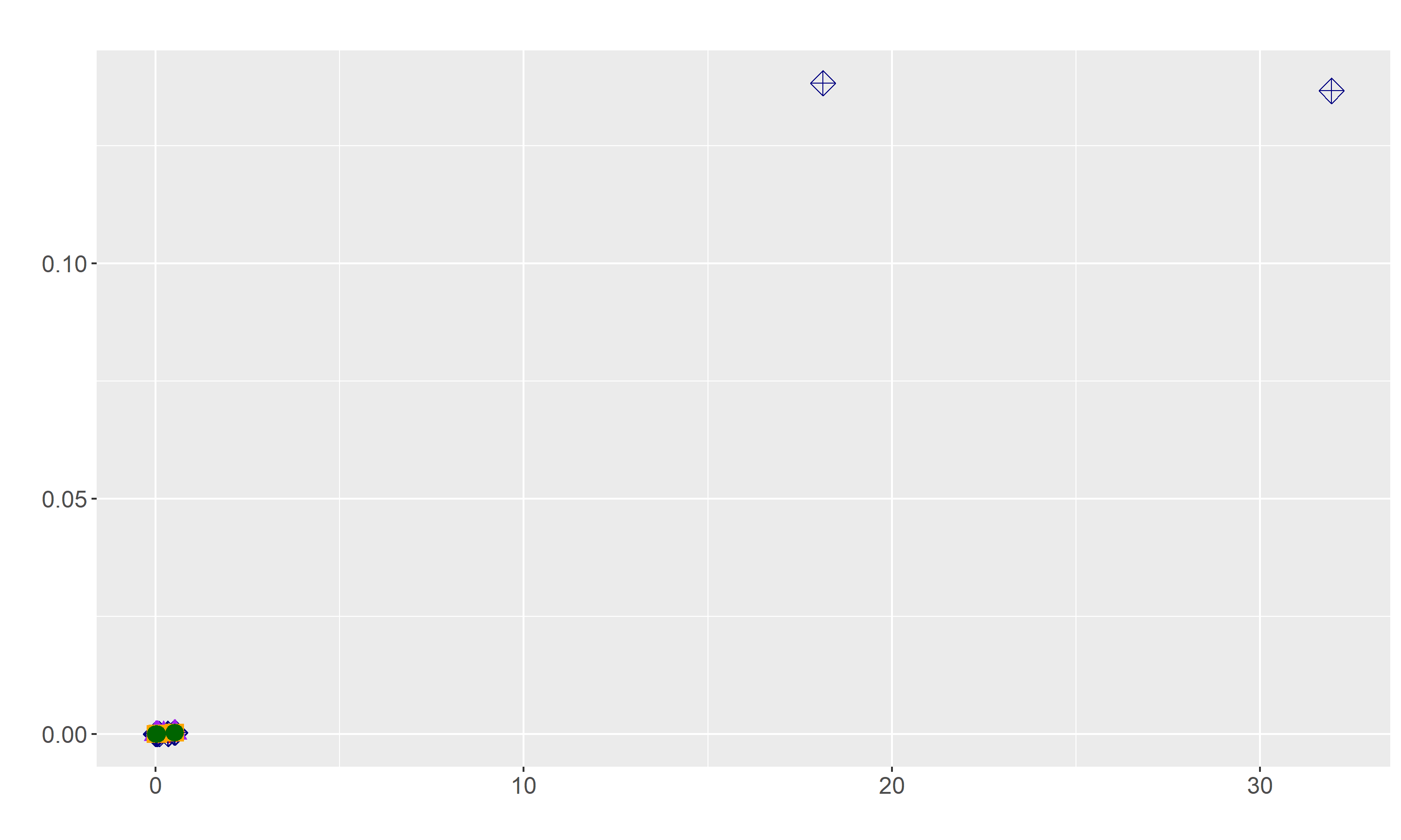}};
			\node(s2) at (2,0){\includegraphics[width=1cm]{figs/leyend_local.png}};
			\node[rotate=90](y) at (-3.4,0){{\color{darkgray}\tiny{$\log(\text{error}+1)$}}};
			\node(x) at (0,-2){{\color{darkgray}\tiny{distance}}};
		\end{tikzpicture}
		\caption{Local solvers: Distance to ref. solution vs error.}
\end{subfigure}

\begin{subfigure}{0.495\textwidth}
		\centering
		\begin{tikzpicture}
			\node(s1) at (0,0){\includegraphics[width=0.9\textwidth]{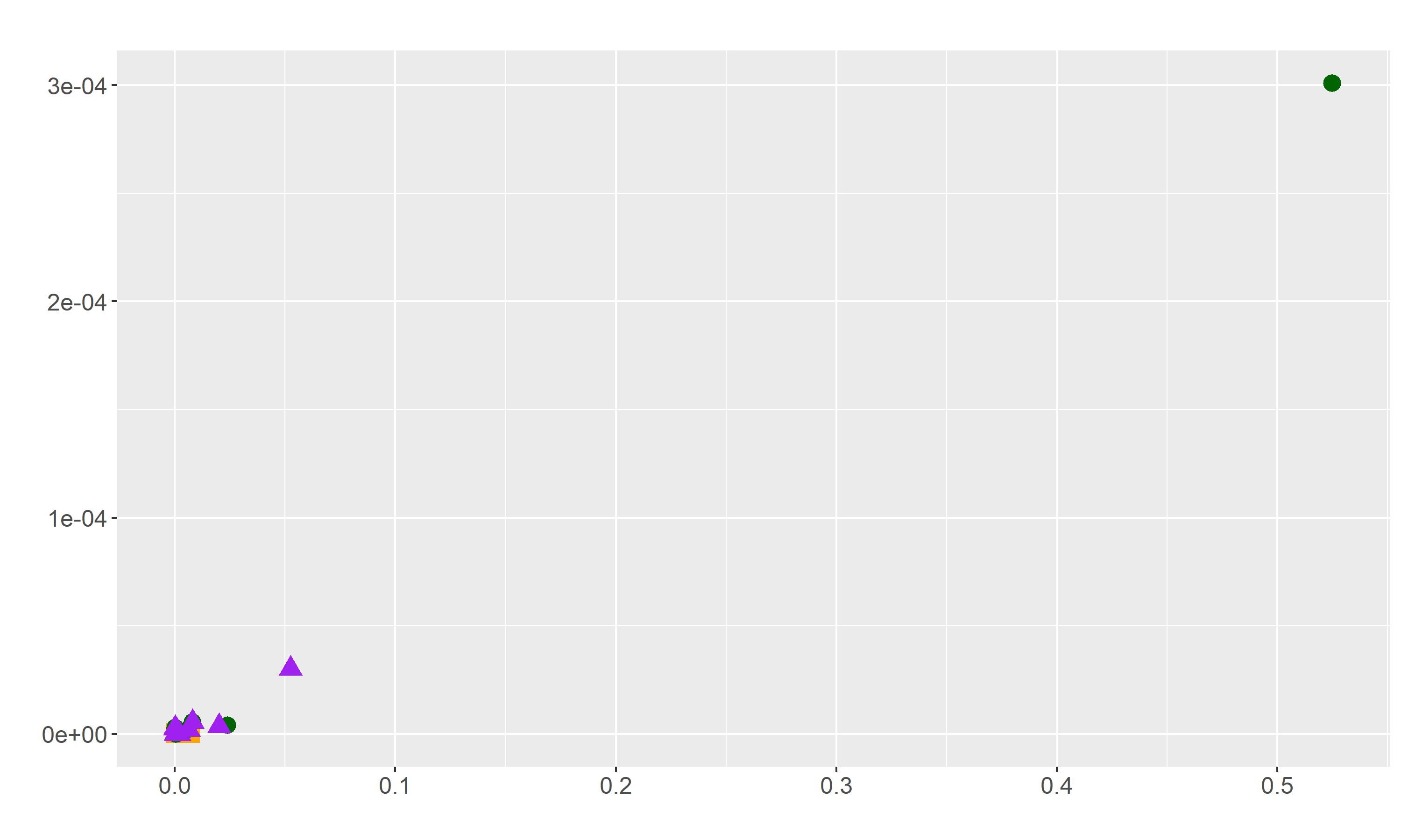}};
			\node(s2) at (2,0){\includegraphics[width=1cm]{figs/leyend_global.png}};
			\node[rotate=90](y) at (-3.4,0){{\color{darkgray}\tiny{$\log(\text{error}+1)$}}};
			\node(x) at (0,-2){{\color{darkgray}\tiny{distance}}};
		\end{tikzpicture}
		\caption{Global solvers: Distance to ref. solution vs error when ``solved''.}
\end{subfigure}

\caption{Results for problem \lvF.}
\label{fig:lvF}
\end{figure}

\newpage 
\subsection{\lvP}
\begin{figure}[!htbp]
\centering
\begin{subfigure}{0.495\textwidth}
		\centering
		\begin{tikzpicture}
			\node(s1) at (0,0){\includegraphics[width=0.9\textwidth]{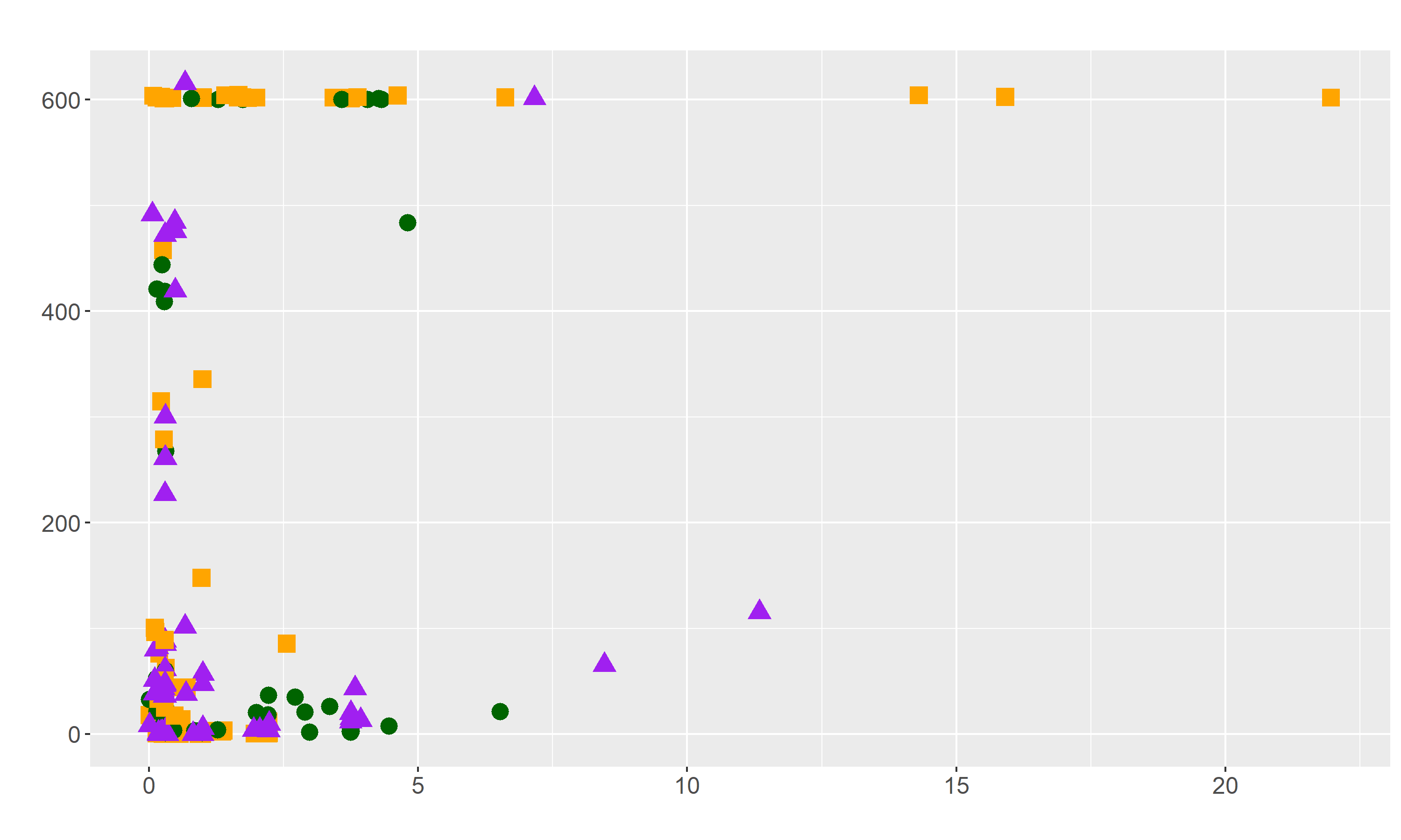}};
			\node[rotate=90](y) at (-3.35,0){{\color{darkgray}\tiny{time}}};
			\node(x) at (0,-2){{\color{darkgray}\tiny{distance}}};
		\end{tikzpicture}
		\caption{Global solvers: Distance to ref. solution vs time.}
\end{subfigure}
\begin{subfigure}{0.495\textwidth}
		\centering
		\begin{tikzpicture}
			\node(s1) at (0,0){\includegraphics[width=0.9\textwidth]{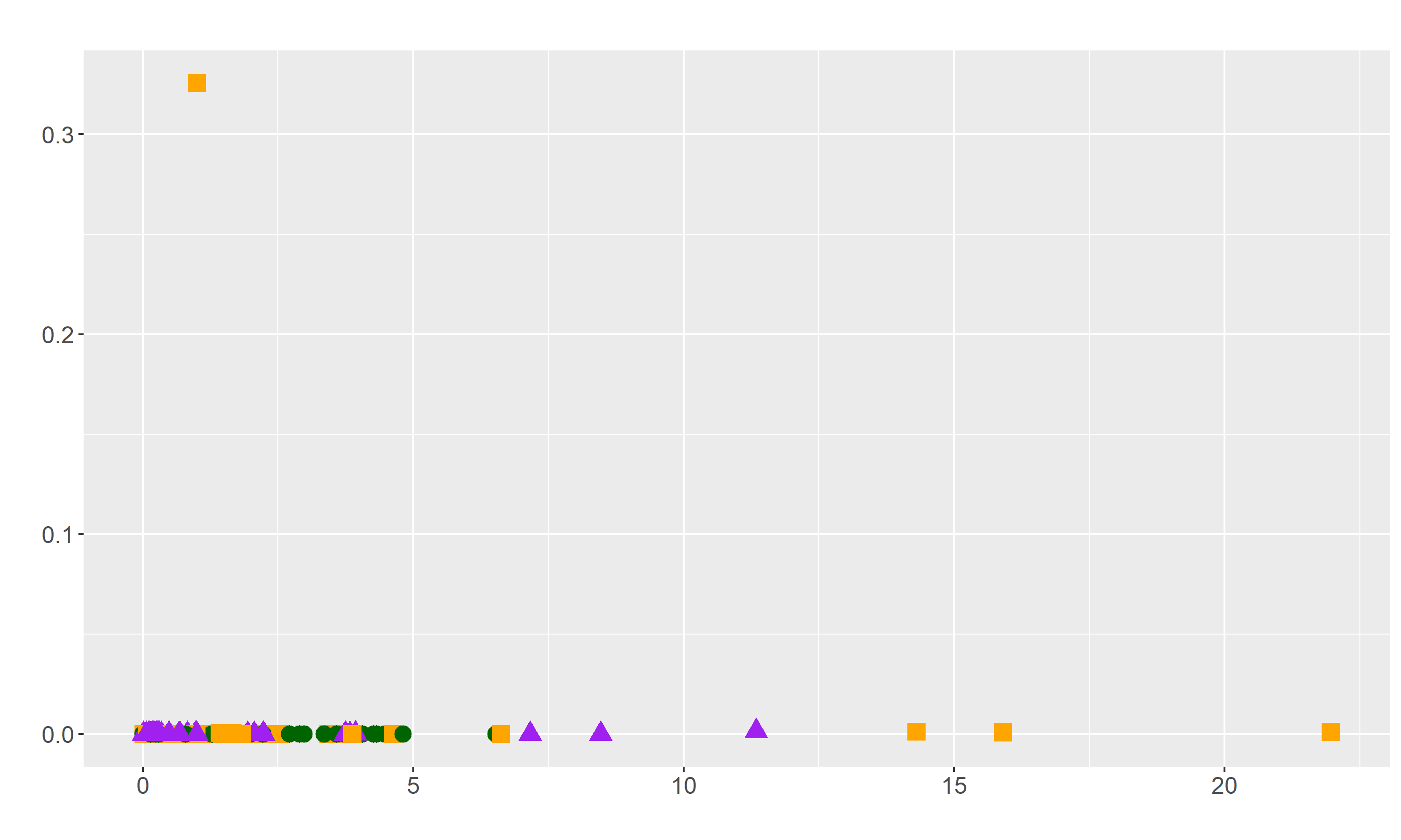}};
		\node(s2) at (2,0){\includegraphics[width=1cm]{figs/leyend_global.png}};
		\node[rotate=90](y) at (-3.4,0){{\color{darkgray}\tiny{$\log(\text{error}+1)$}}};
		\node(x) at (0,-2){{\color{darkgray}\tiny{distance}}};
		\end{tikzpicture}
		\caption{Global solvers: Distance to ref. solution vs error.}
\end{subfigure}

\vspace{0.2cm}

\begin{subfigure}{0.495\textwidth}
		\centering
		\begin{tikzpicture}
			\node(s1) at (0,0){\includegraphics[width=0.9\textwidth]{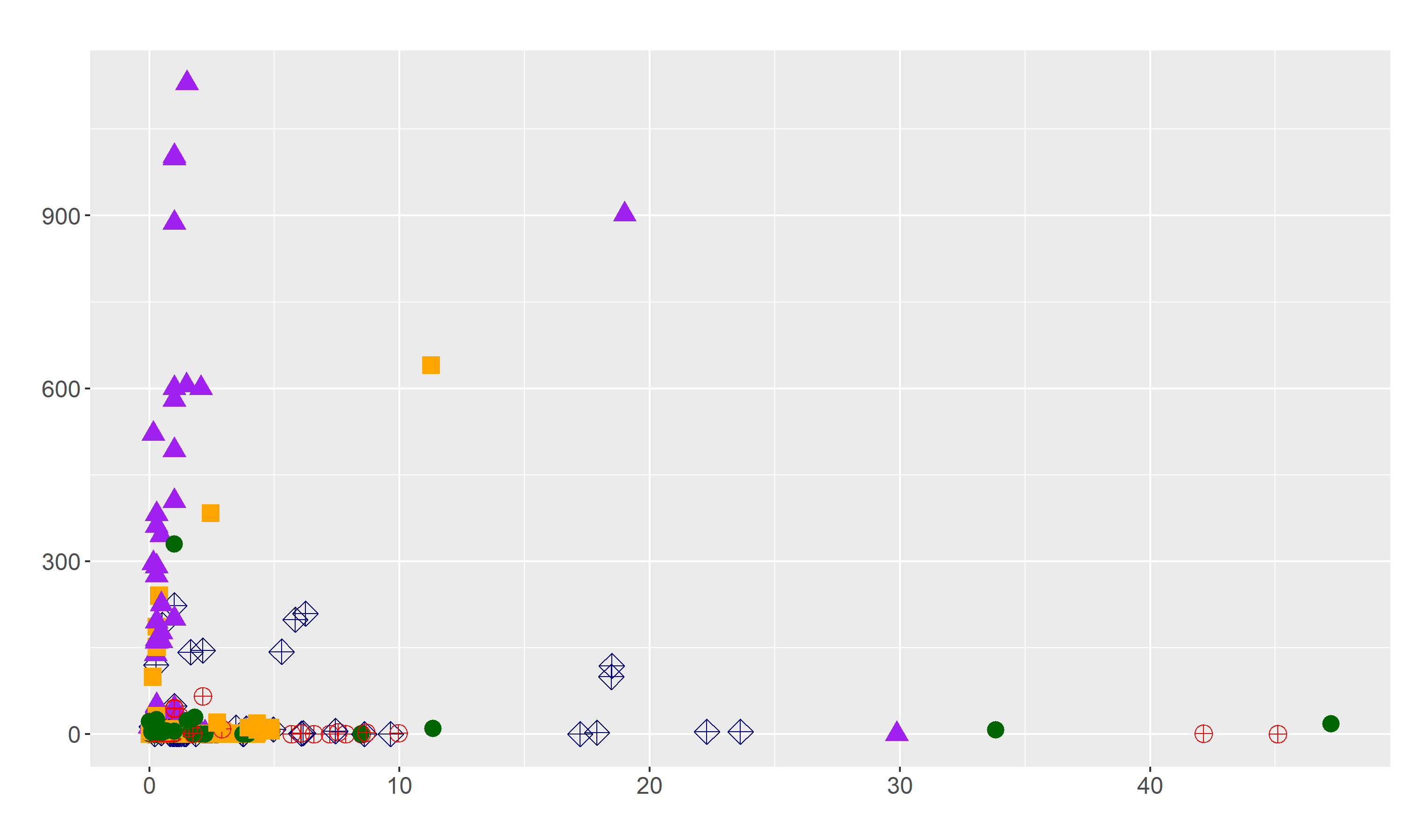}};
			\node[rotate=90](y) at (-3.35,0){{\color{darkgray}\tiny{time}}};
			\node(x) at (0,-2){{\color{darkgray}\tiny{distance}}};
		\end{tikzpicture}
		\caption{Local solvers: Distance to ref. solution vs time.}
\end{subfigure}
\begin{subfigure}{0.495\textwidth}
		\centering
		\begin{tikzpicture}
			\node(s1) at (0,0){\includegraphics[width=0.9\textwidth]{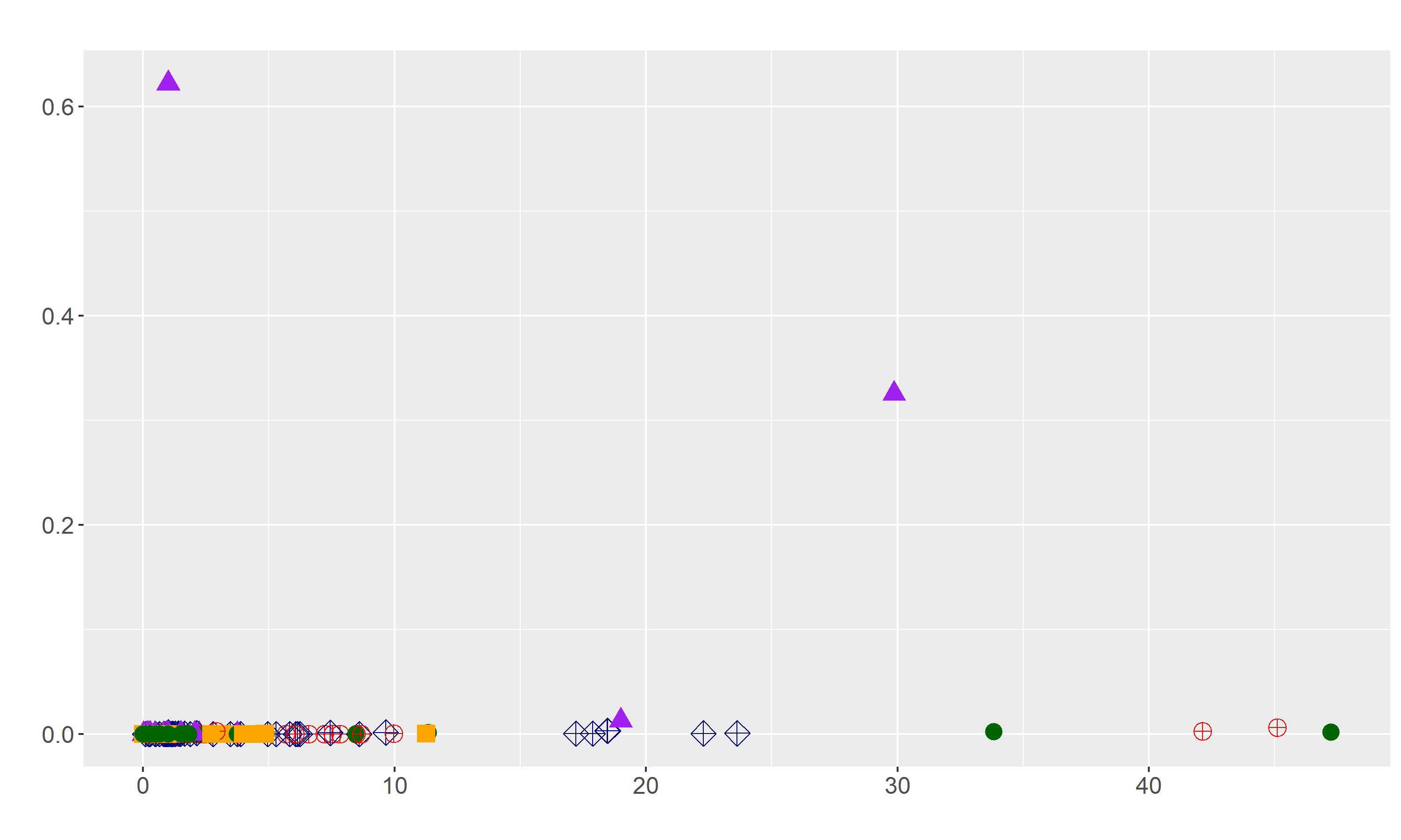}};
			\node(s2) at (2,0){\includegraphics[width=1cm]{figs/leyend_local.png}};
			\node[rotate=90](y) at (-3.4,0){{\color{darkgray}\tiny{$\log(\text{error}+1)$}}};
			\node(x) at (0,-2){{\color{darkgray}\tiny{distance}}};
		\end{tikzpicture}
		\caption{Local solvers: Distance to ref. solution vs error.}
\end{subfigure}

\begin{subfigure}{0.495\textwidth}
		\centering
		\begin{tikzpicture}
			\node(s1) at (0,0){\includegraphics[width=0.9\textwidth]{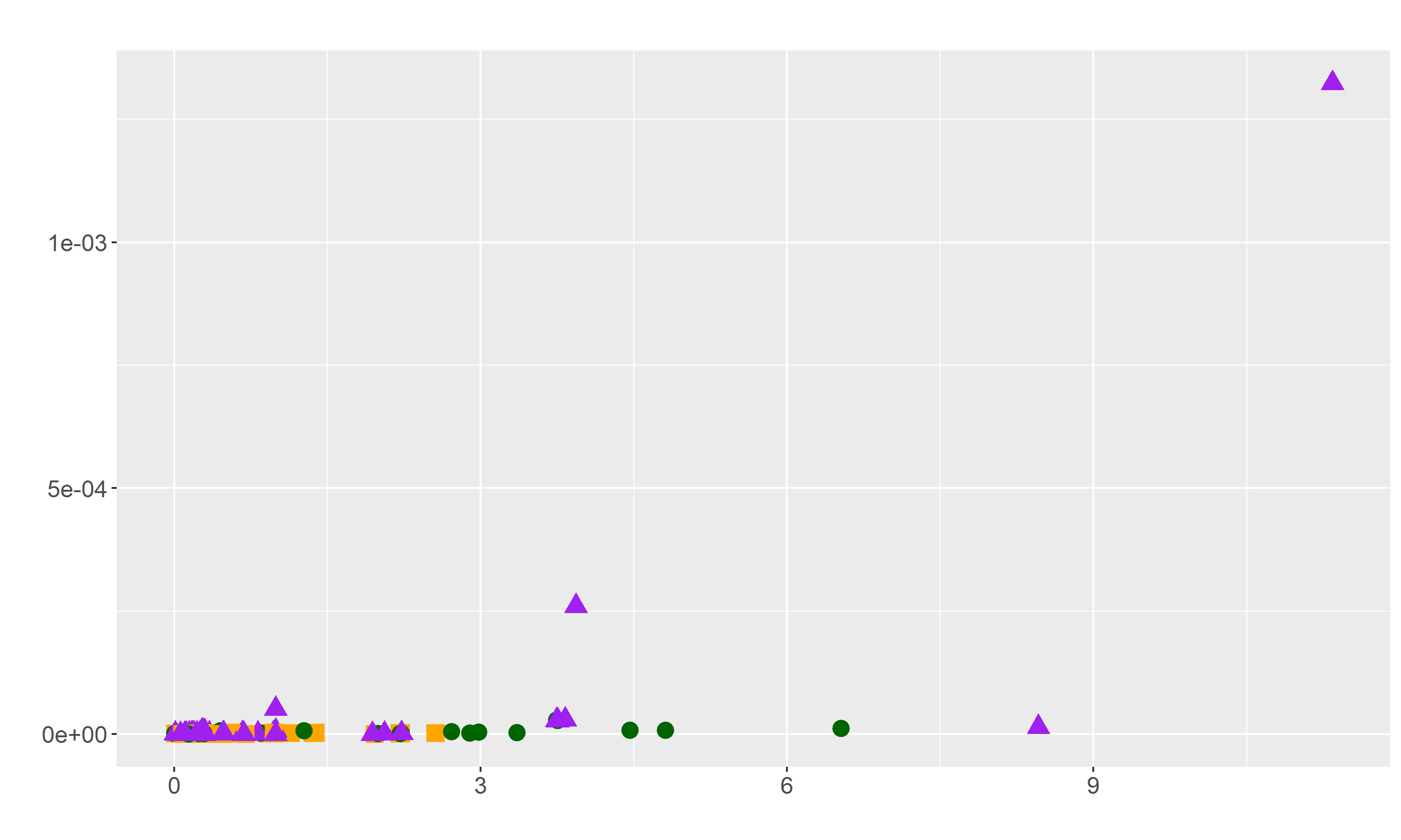}};
			\node(s2) at (2,0){\includegraphics[width=1cm]{figs/leyend_global.png}};
			\node[rotate=90](y) at (-3.4,0){{\color{darkgray}\tiny{$\log(\text{error}+1)$}}};
			\node(x) at (0,-2){{\color{darkgray}\tiny{distance}}};
		\end{tikzpicture}
		\caption{Global solvers: Distance to ref. solution vs error when ``solved''.}
\end{subfigure}

\caption{Results for problem \lvP.}
\label{fig:lvP}
\end{figure}

\newpage 
\subsection{\FHN}
\begin{figure}[!htbp]
\centering
\begin{subfigure}{0.495\textwidth}
		\centering
		\begin{tikzpicture}
			\node(s1) at (0,0){\includegraphics[width=0.9\textwidth]{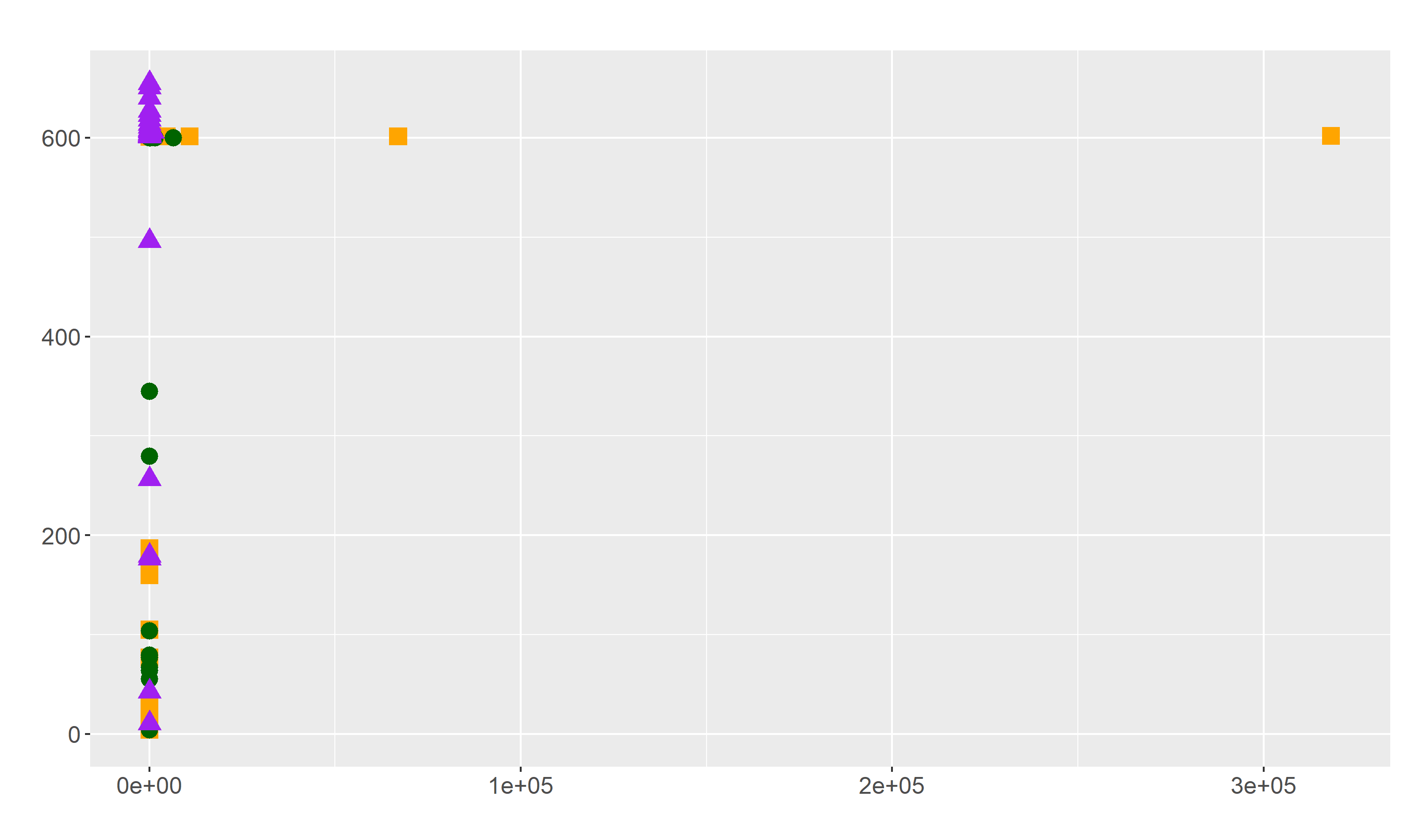}};
			\node[rotate=90](y) at (-3.35,0){{\color{darkgray}\tiny{time}}};
			\node(x) at (0,-2){{\color{darkgray}\tiny{distance}}};
		\end{tikzpicture}
		\caption{Global solvers: Distance to ref. solution vs time.}
\end{subfigure}
\begin{subfigure}{0.495\textwidth}
		\centering
		\begin{tikzpicture}
			\node(s1) at (0,0){\includegraphics[width=0.9\textwidth]{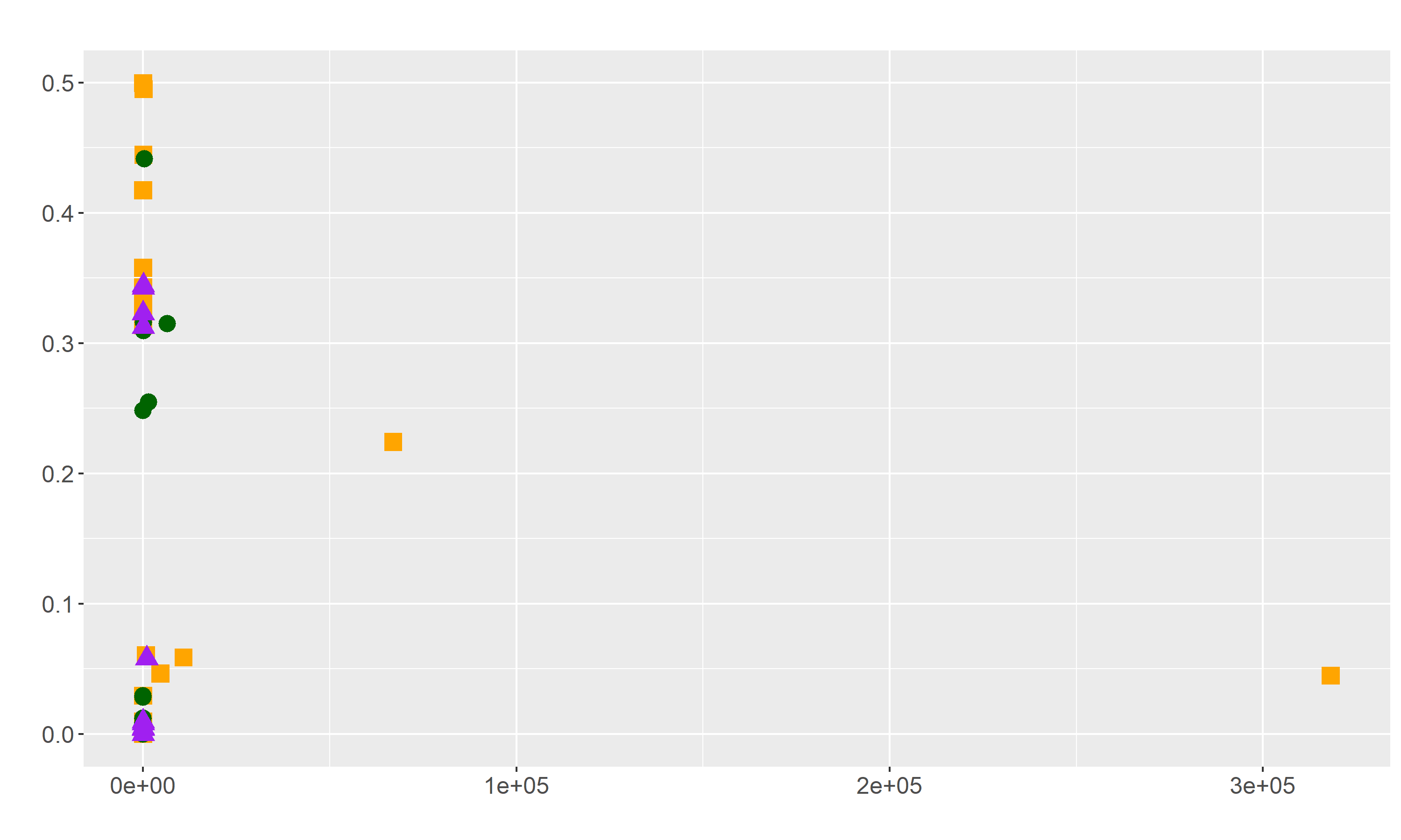}};
		\node(s2) at (2,0){\includegraphics[width=1cm]{figs/leyend_global.png}};
		\node[rotate=90](y) at (-3.4,0){{\color{darkgray}\tiny{$\log(\text{error}+1)$}}};
		\node(x) at (0,-2){{\color{darkgray}\tiny{distance}}};
		\end{tikzpicture}
		\caption{Global solvers: Distance to ref. solution vs error.}
\end{subfigure}

\vspace{0.2cm}

\begin{subfigure}{0.495\textwidth}
		\centering
		\begin{tikzpicture}
			\node(s1) at (0,0){\includegraphics[width=0.9\textwidth]{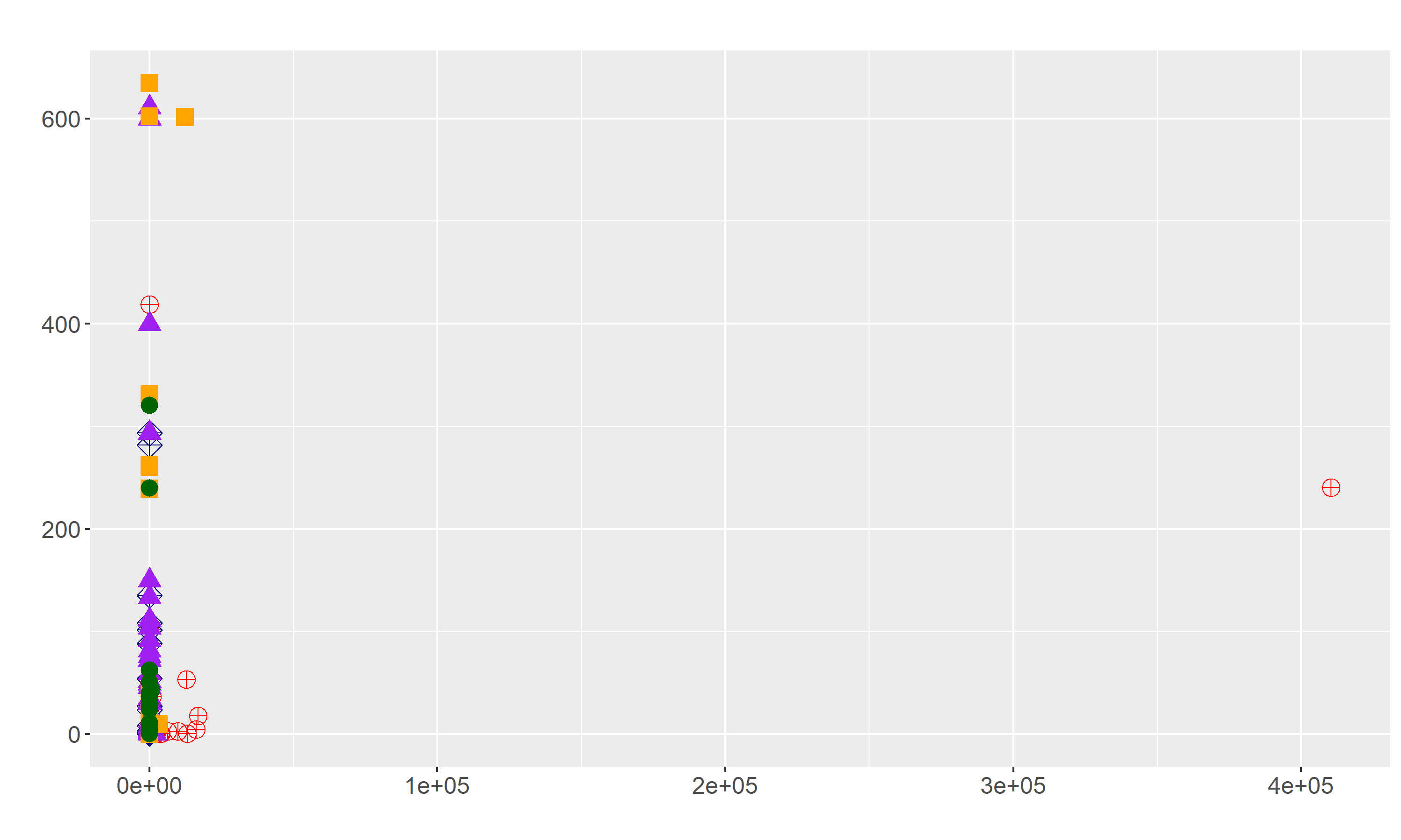}};
			\node[rotate=90](y) at (-3.35,0){{\color{darkgray}\tiny{time}}};
			\node(x) at (0,-2){{\color{darkgray}\tiny{distance}}};
		\end{tikzpicture}
		\caption{Local solvers: Distance to ref. solution vs time.}
\end{subfigure}
\begin{subfigure}{0.495\textwidth}
		\centering
		\begin{tikzpicture}
			\node(s1) at (0,0){\includegraphics[width=0.9\textwidth]{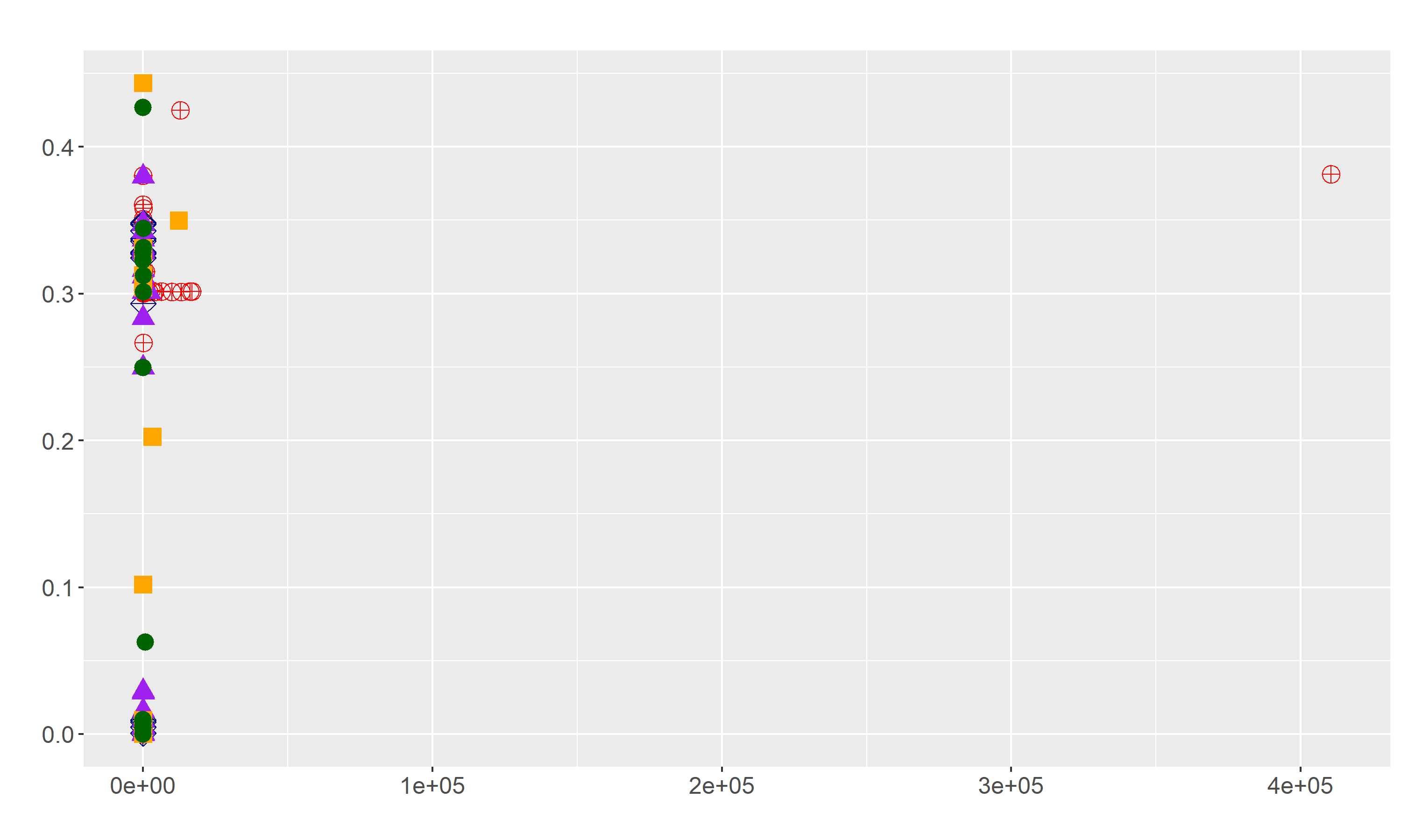}};
			\node(s2) at (2,0){\includegraphics[width=1cm]{figs/leyend_local.png}};
			\node[rotate=90](y) at (-3.4,0){{\color{darkgray}\tiny{$\log(\text{error}+1)$}}};
			\node(x) at (0,-2){{\color{darkgray}\tiny{distance}}};
		\end{tikzpicture}
		\caption{Local solvers: Distance to ref. solution vs error.}
\end{subfigure}

\begin{subfigure}{0.495\textwidth}
		\centering
		\begin{tikzpicture}
			\node(s1) at (0,0){\includegraphics[width=0.9\textwidth]{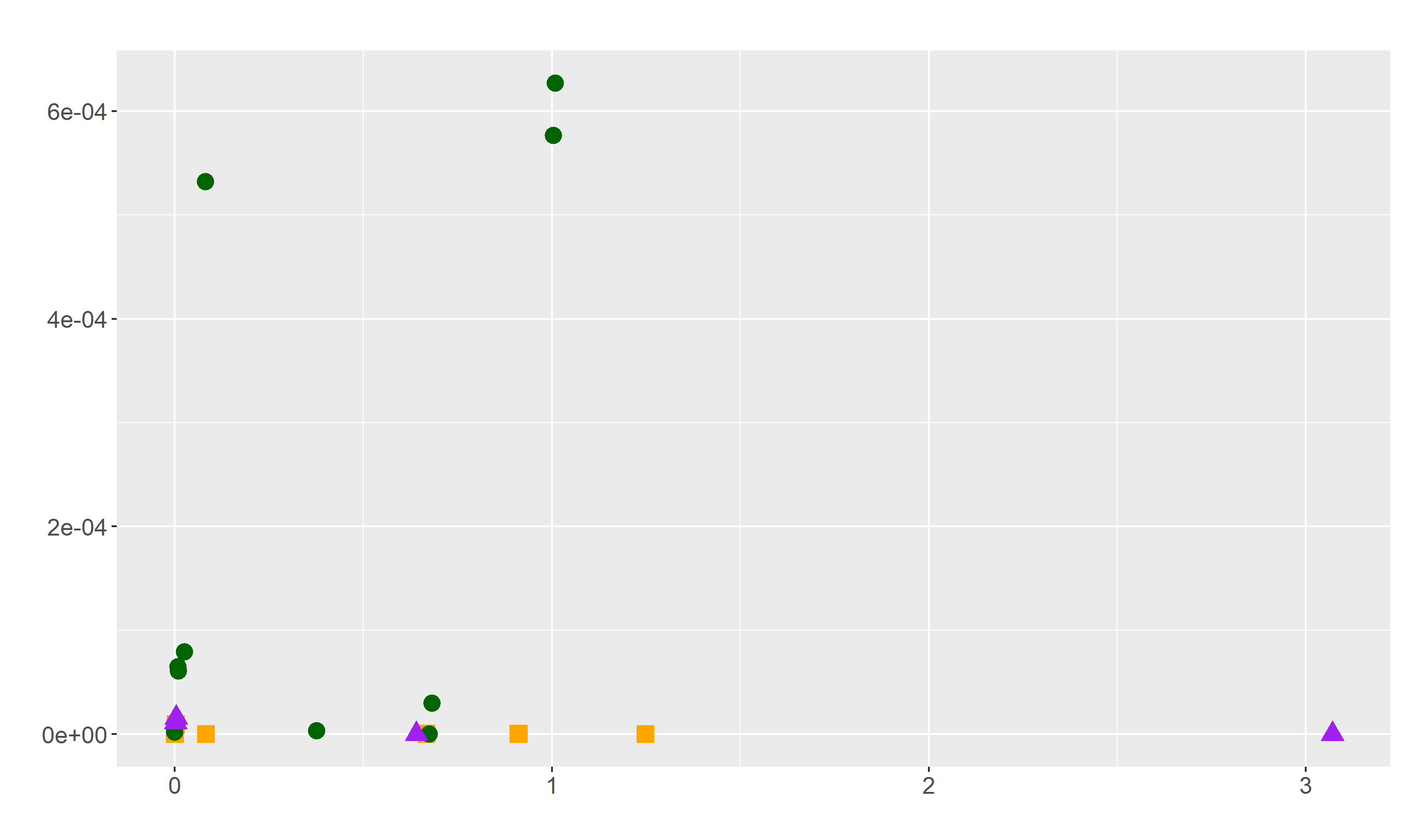}};
			\node(s2) at (2,0){\includegraphics[width=1cm]{figs/leyend_global.png}};
			\node[rotate=90](y) at (-3.4,0){{\color{darkgray}\tiny{$\log(\text{error}+1)$}}};
			\node(x) at (0,-2){{\color{darkgray}\tiny{distance}}};
		\end{tikzpicture}
		\caption{Global solvers: Distance to ref. solution vs error when ``solved''.}
\end{subfigure}

\caption{Results for problem \FHN.}
\label{fig:FHN}
\end{figure}

\newpage 
\subsection{\crausteF}
\begin{figure}[!htbp]
\centering
\begin{subfigure}{0.495\textwidth}
		\centering
		\begin{tikzpicture}
			\node(s1) at (0,0){\includegraphics[width=0.9\textwidth]{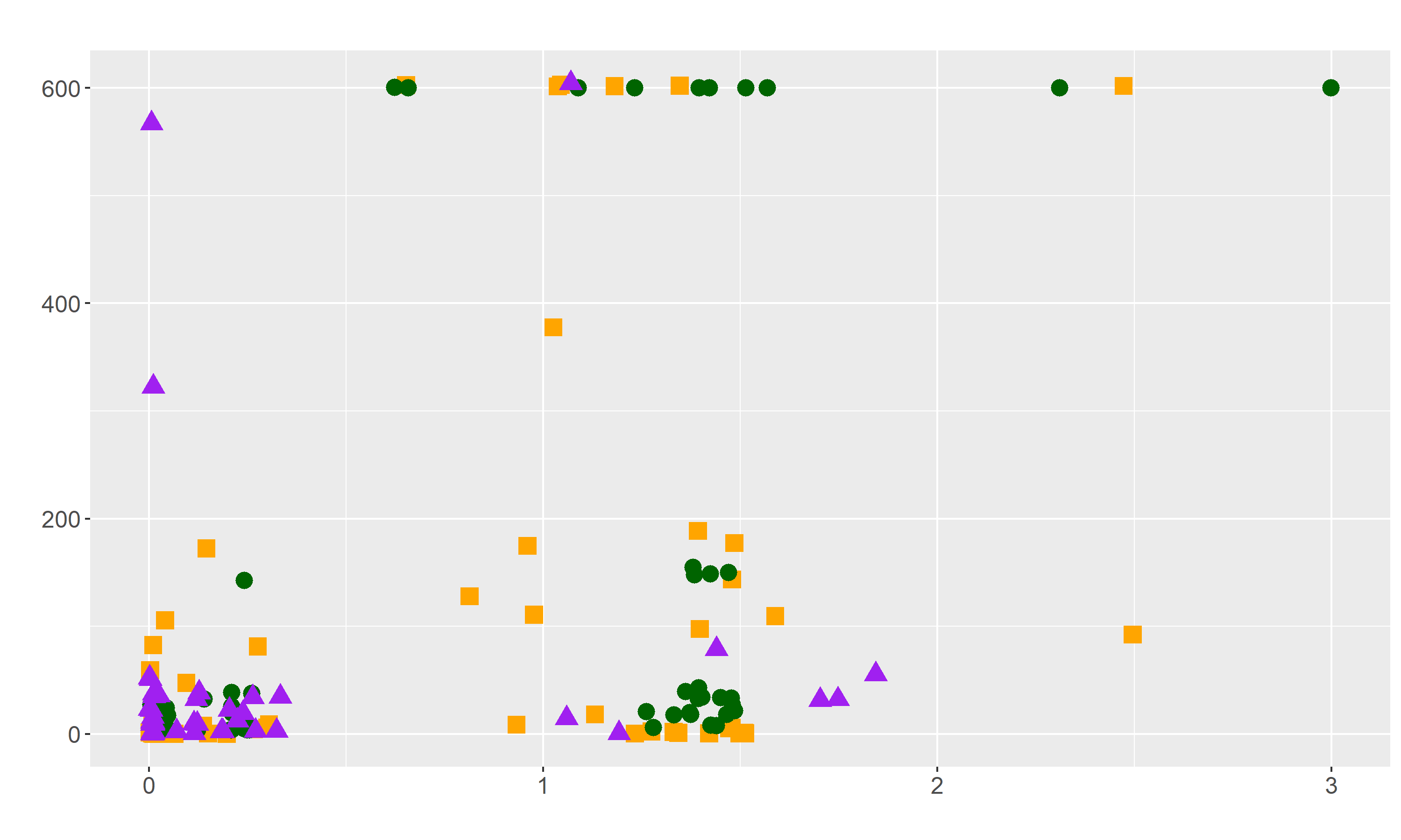}};
			\node[rotate=90](y) at (-3.35,0){{\color{darkgray}\tiny{time}}};
			\node(x) at (0,-2){{\color{darkgray}\tiny{distance}}};
		\end{tikzpicture}
		\caption{Global solvers: Distance to ref. solution vs time.}
\end{subfigure}
\begin{subfigure}{0.495\textwidth}
		\centering
		\begin{tikzpicture}
			\node(s1) at (0,0){\includegraphics[width=0.9\textwidth]{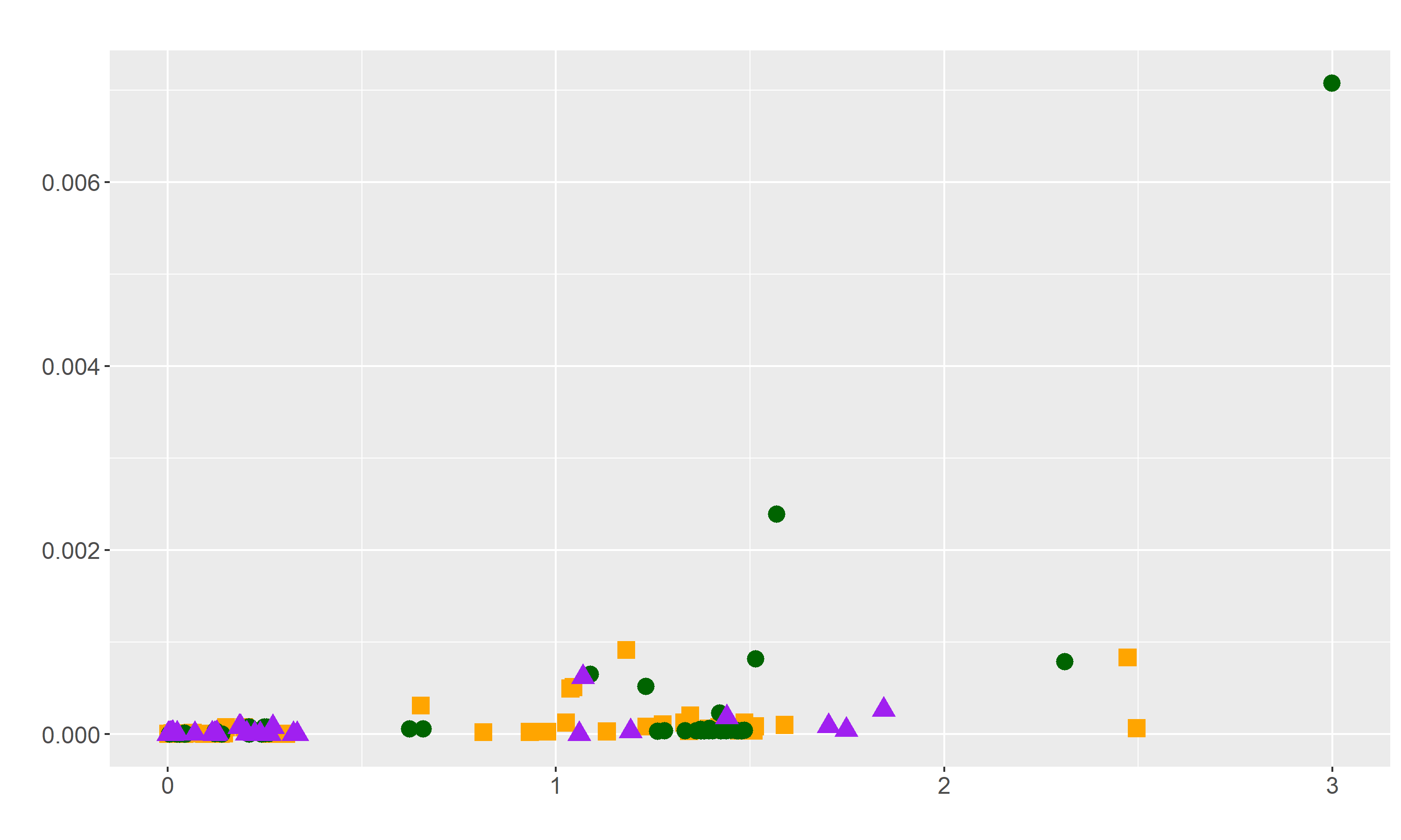}};
		\node(s2) at (2,0){\includegraphics[width=1cm]{figs/leyend_global.png}};
		\node[rotate=90](y) at (-3.4,0){{\color{darkgray}\tiny{$\log(\text{error}+1)$}}};
		\node(x) at (0,-2){{\color{darkgray}\tiny{distance}}};
		\end{tikzpicture}
		\caption{Global solvers: Distance to ref. solution vs error.}
\end{subfigure}

\vspace{0.2cm}

\begin{subfigure}{0.495\textwidth}
		\centering
		\begin{tikzpicture}
			\node(s1) at (0,0){\includegraphics[width=0.9\textwidth]{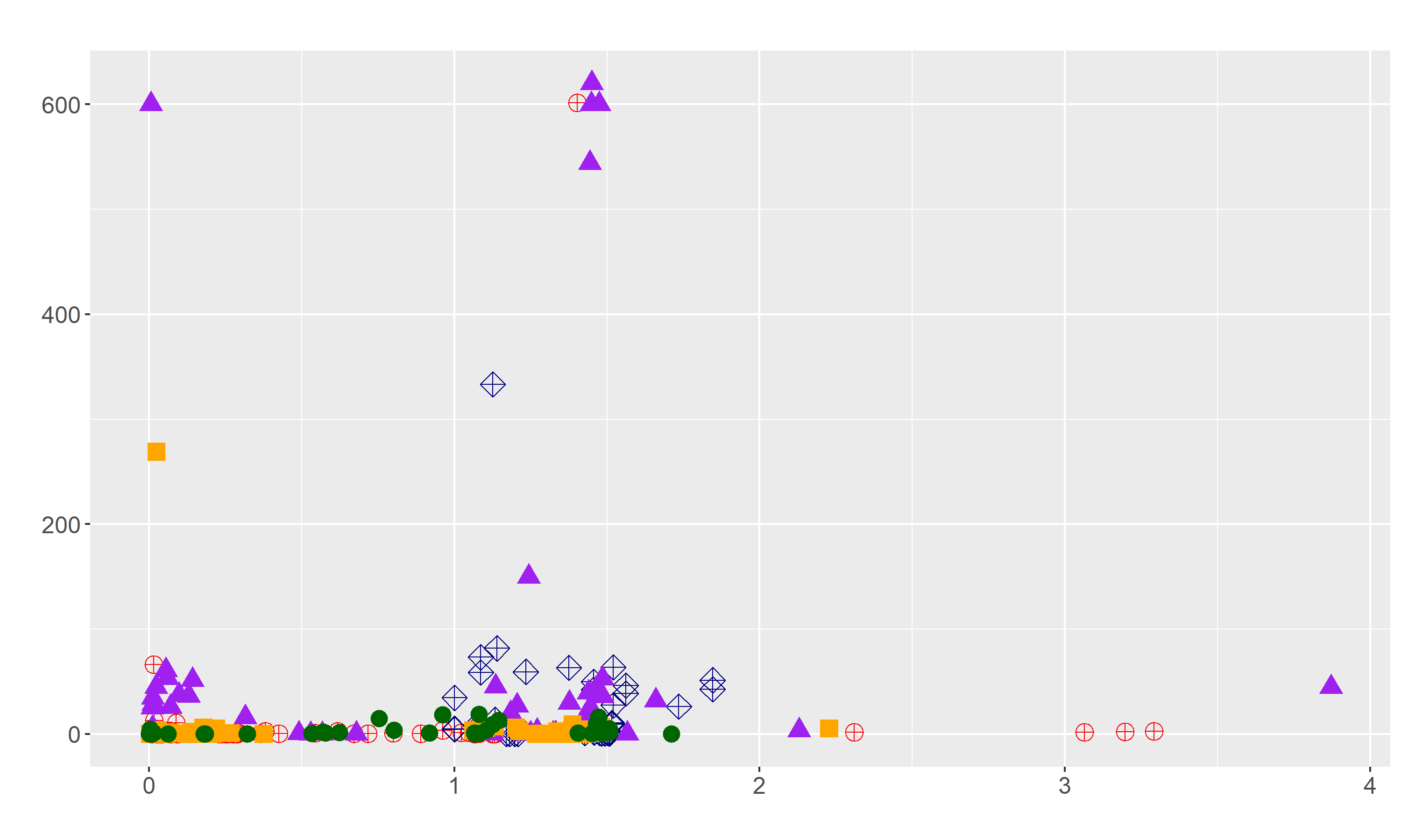}};
			\node[rotate=90](y) at (-3.35,0){{\color{darkgray}\tiny{time}}};
			\node(x) at (0,-2){{\color{darkgray}\tiny{distance}}};
		\end{tikzpicture}
		\caption{Local solvers: Distance to ref. solution vs time.}
\end{subfigure}
\begin{subfigure}{0.495\textwidth}
		\centering
		\begin{tikzpicture}
			\node(s1) at (0,0){\includegraphics[width=0.9\textwidth]{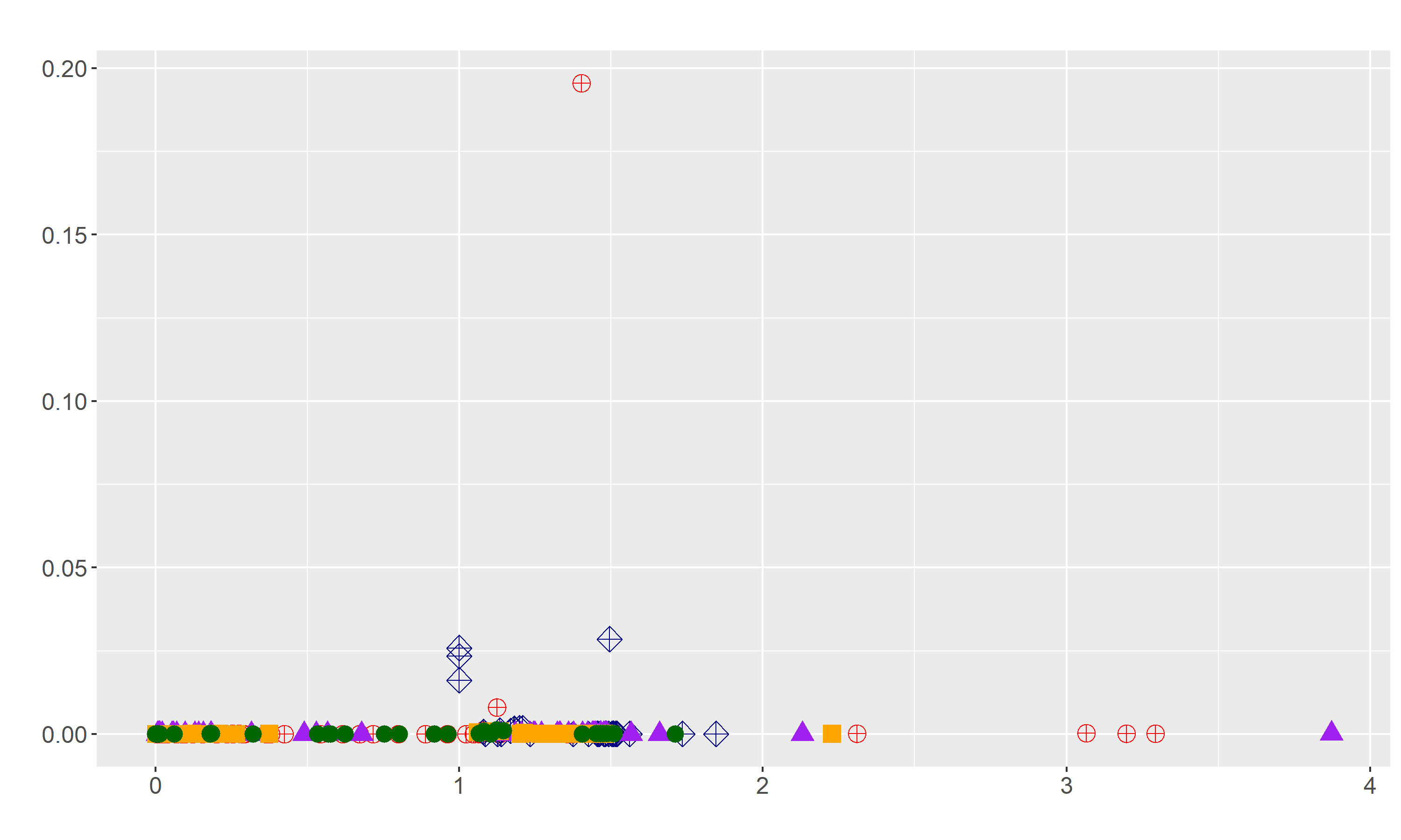}};
			\node(s2) at (2,0){\includegraphics[width=1cm]{figs/leyend_local.png}};
			\node[rotate=90](y) at (-3.4,0){{\color{darkgray}\tiny{$\log(\text{error}+1)$}}};
			\node(x) at (0,-2){{\color{darkgray}\tiny{distance}}};
		\end{tikzpicture}
		\caption{Local solvers: Distance to ref. solution vs error.}
\end{subfigure}

\begin{subfigure}{0.495\textwidth}
		\centering
		\begin{tikzpicture}
			\node(s1) at (0,0){\includegraphics[width=0.9\textwidth]{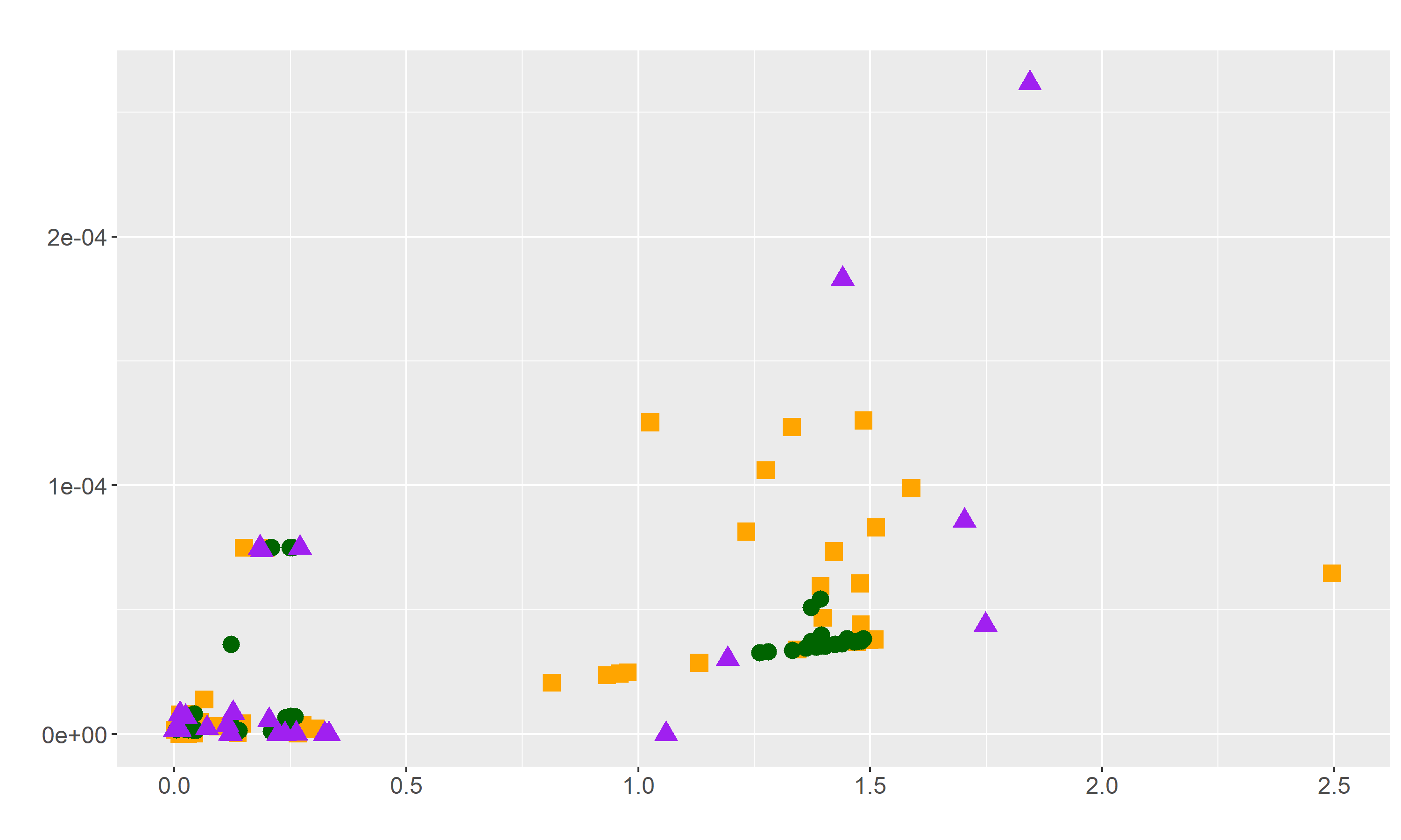}};
			\node(s2) at (2,0){\includegraphics[width=1cm]{figs/leyend_global.png}};
			\node[rotate=90](y) at (-3.4,0){{\color{darkgray}\tiny{$\log(\text{error}+1)$}}};
			\node(x) at (0,-2){{\color{darkgray}\tiny{distance}}};
		\end{tikzpicture}
		\caption{Global solvers: Distance to ref. solution vs error when ``solved''.}
\end{subfigure}

\caption{Results for problem \crausteF.}
\label{fig:crausteF}
\end{figure}

\newpage 
\subsection{\crausteP}
\begin{figure}[!htbp]
\centering
\begin{subfigure}{0.495\textwidth}
		\centering
		\begin{tikzpicture}
			\node(s1) at (0,0){\includegraphics[width=0.9\textwidth]{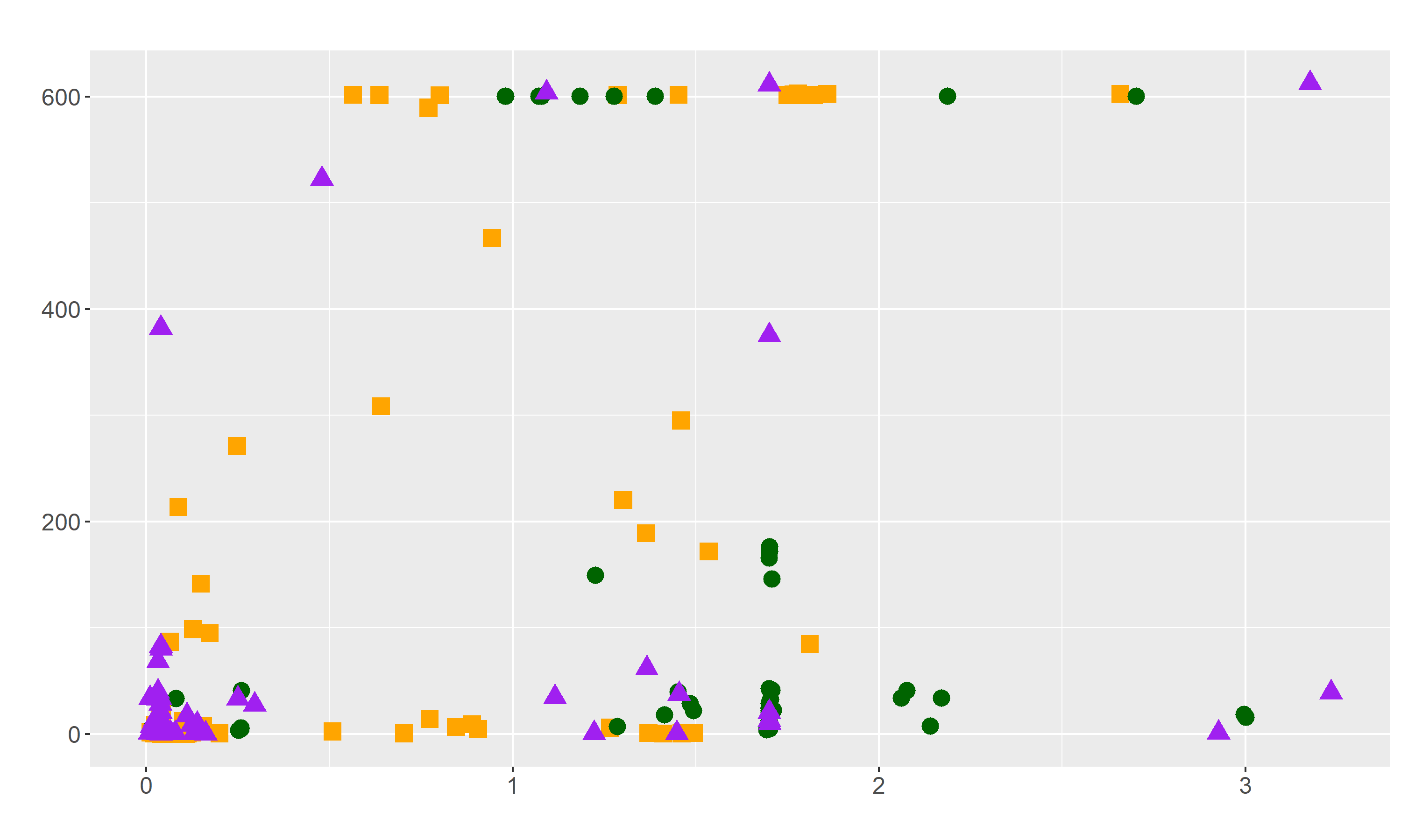}};
			\node[rotate=90](y) at (-3.35,0){{\color{darkgray}\tiny{time}}};
			\node(x) at (0,-2){{\color{darkgray}\tiny{distance}}};
		\end{tikzpicture}
		\caption{Global solvers: Distance to ref. solution vs time.}
\end{subfigure}
\begin{subfigure}{0.495\textwidth}
		\centering
		\begin{tikzpicture}
			\node(s1) at (0,0){\includegraphics[width=0.9\textwidth]{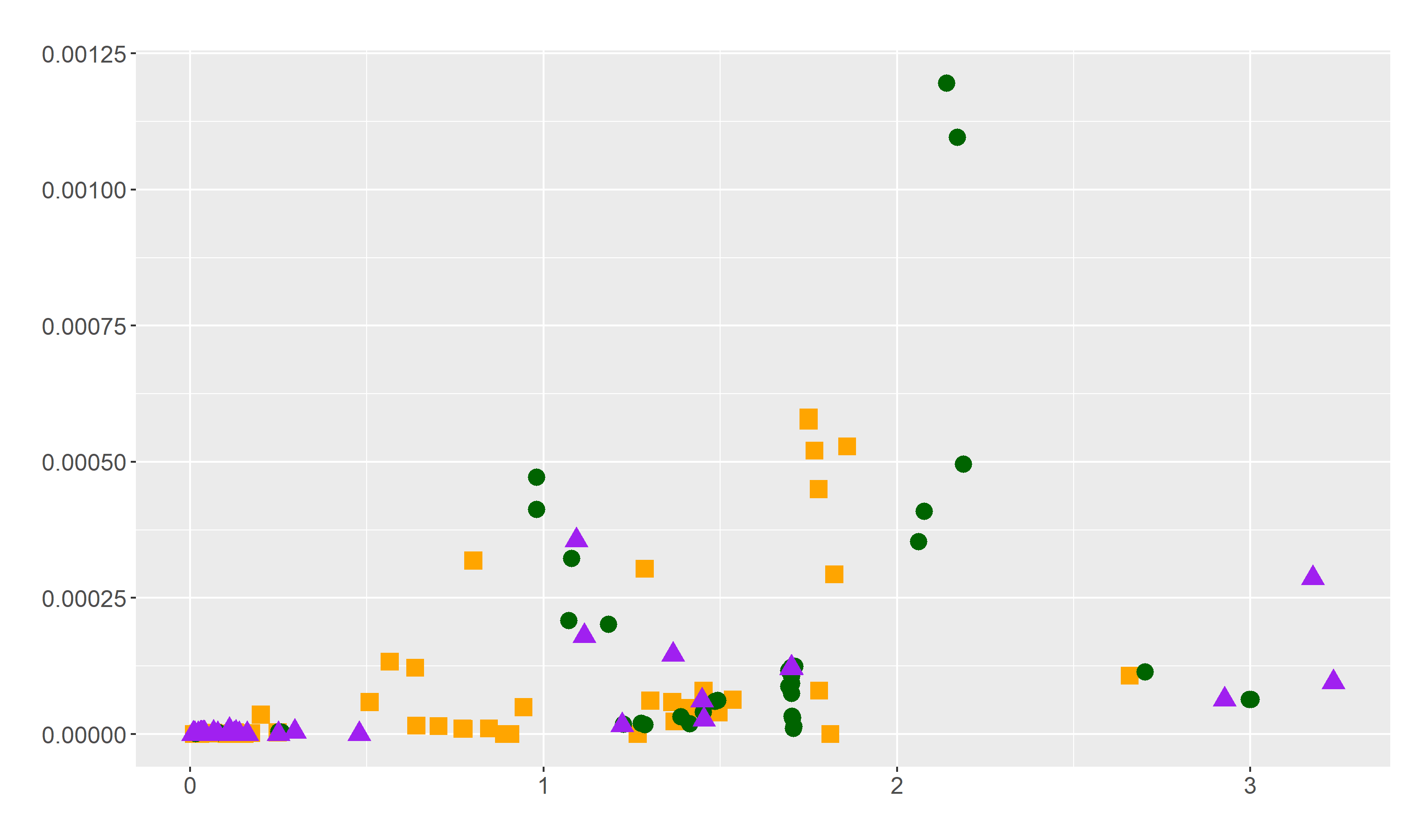}};
		\node(s2) at (2,0){\includegraphics[width=1cm]{figs/leyend_global.png}};
		\node[rotate=90](y) at (-3.4,0){{\color{darkgray}\tiny{$\log(\text{error}+1)$}}};
		\node(x) at (0,-2){{\color{darkgray}\tiny{distance}}};
		\end{tikzpicture}
		\caption{Global solvers: Distance to ref. solution vs error.}
\end{subfigure}

\vspace{0.2cm}

\begin{subfigure}{0.495\textwidth}
		\centering
		\begin{tikzpicture}
			\node(s1) at (0,0){\includegraphics[width=0.9\textwidth]{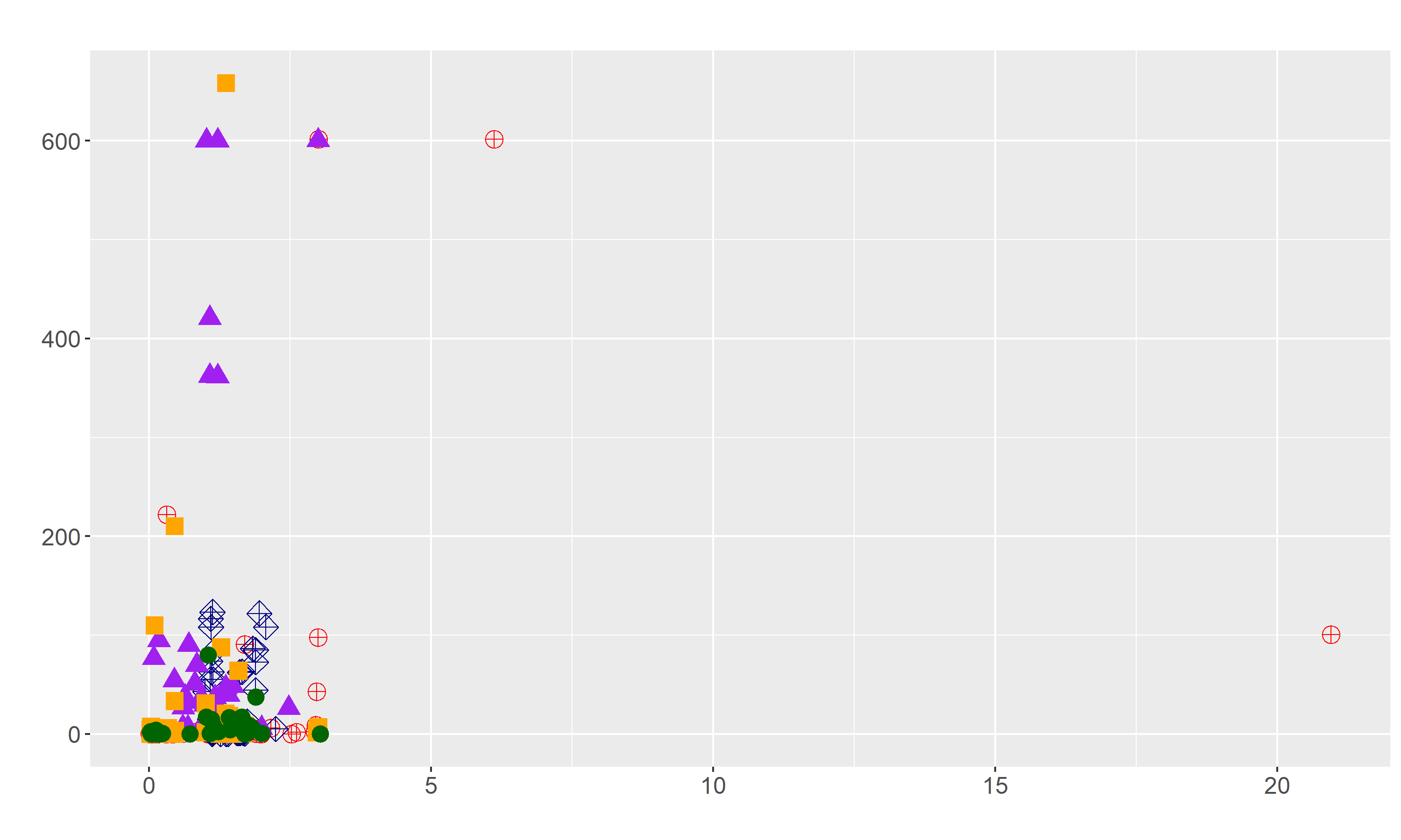}};
			\node[rotate=90](y) at (-3.35,0){{\color{darkgray}\tiny{time}}};
			\node(x) at (0,-2){{\color{darkgray}\tiny{distance}}};
		\end{tikzpicture}
		\caption{Local solvers: Distance to ref. solution vs time.}
\end{subfigure}
\begin{subfigure}{0.495\textwidth}
		\centering
		\begin{tikzpicture}
			\node(s1) at (0,0){\includegraphics[width=0.9\textwidth]{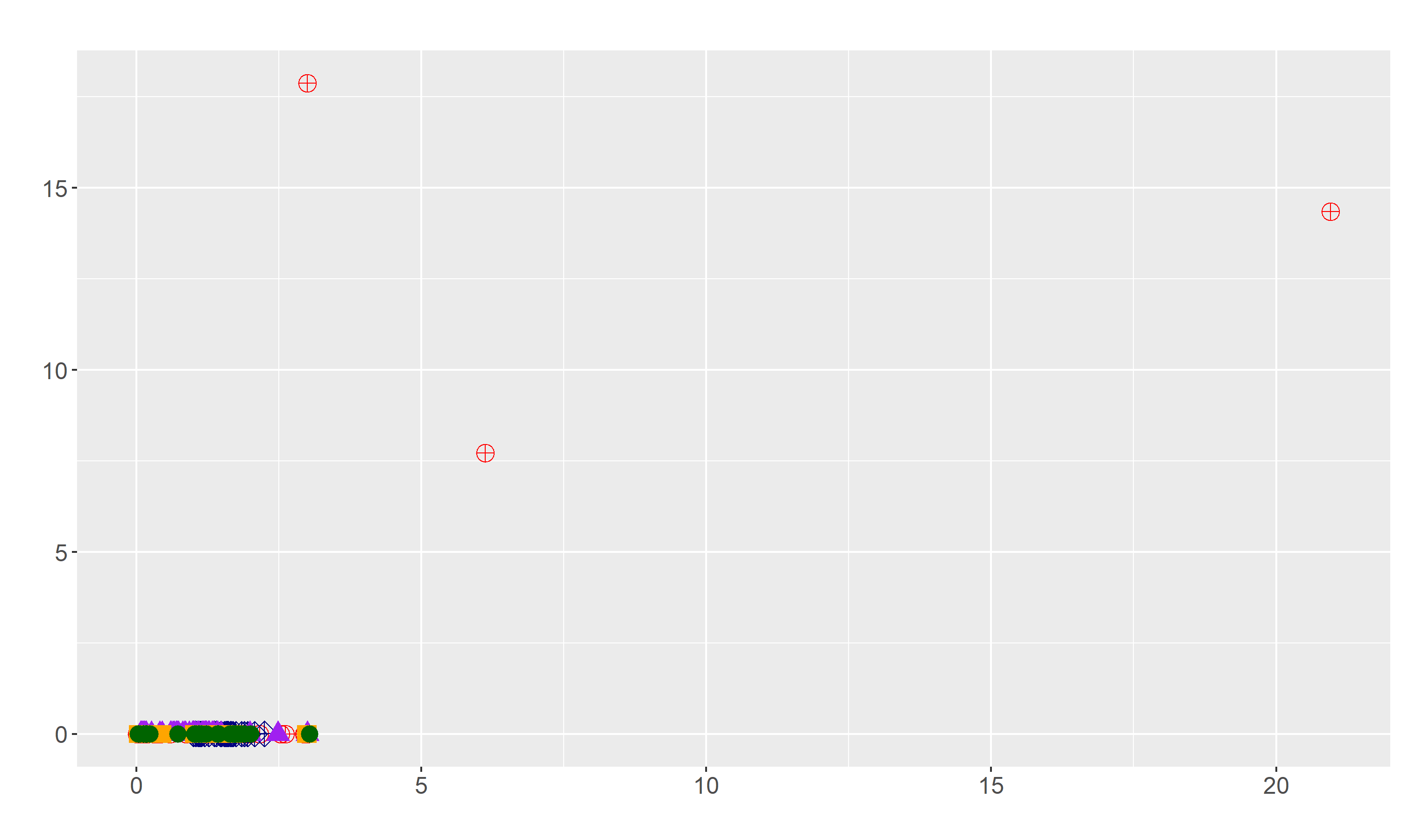}};
			\node(s2) at (2,0){\includegraphics[width=1cm]{figs/leyend_local.png}};
			\node[rotate=90](y) at (-3.4,0){{\color{darkgray}\tiny{$\log(\text{error}+1)$}}};
			\node(x) at (0,-2){{\color{darkgray}\tiny{distance}}};
		\end{tikzpicture}
		\caption{Local solvers: Distance to ref. solution vs error.}
\end{subfigure}

\begin{subfigure}{0.495\textwidth}
		\centering
		\begin{tikzpicture}
			\node(s1) at (0,0){\includegraphics[width=0.9\textwidth]{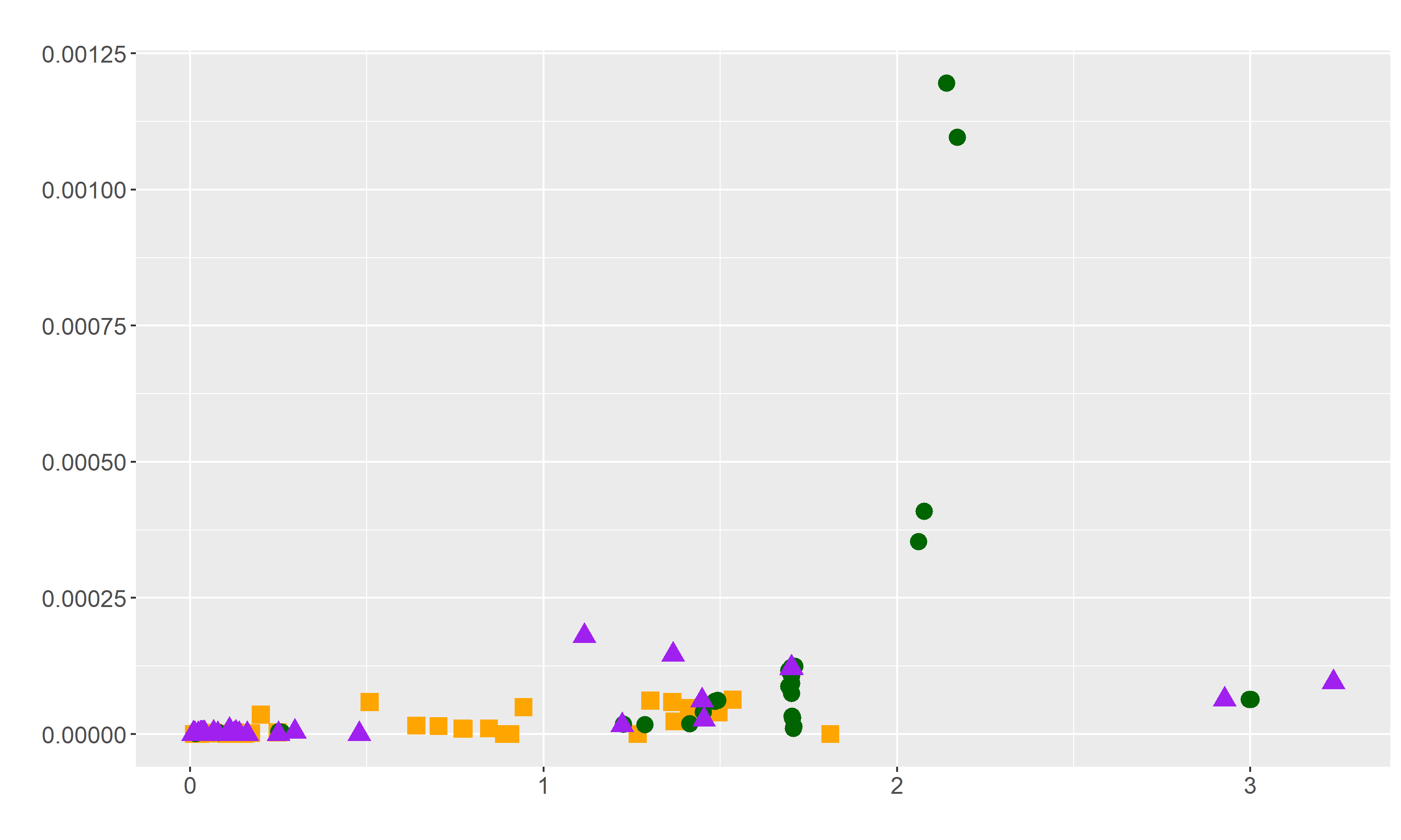}};
			\node(s2) at (2,0){\includegraphics[width=1cm]{figs/leyend_global.png}};
			\node[rotate=90](y) at (-3.4,0){{\color{darkgray}\tiny{$\log(\text{error}+1)$}}};
			\node(x) at (0,-2){{\color{darkgray}\tiny{distance}}};
		\end{tikzpicture}
		\caption{Global solvers: Distance to ref. solution vs error when ``solved''.}
\end{subfigure}

\caption{Results for problem \crausteP.}
\label{fig:crausteP}
\end{figure}

\newpage 
\subsection{\BBG}
\begin{figure}[!htbp]
\centering
\begin{subfigure}{0.495\textwidth}
		\centering
		\begin{tikzpicture}
			\node(s1) at (0,0){\includegraphics[width=0.9\textwidth]{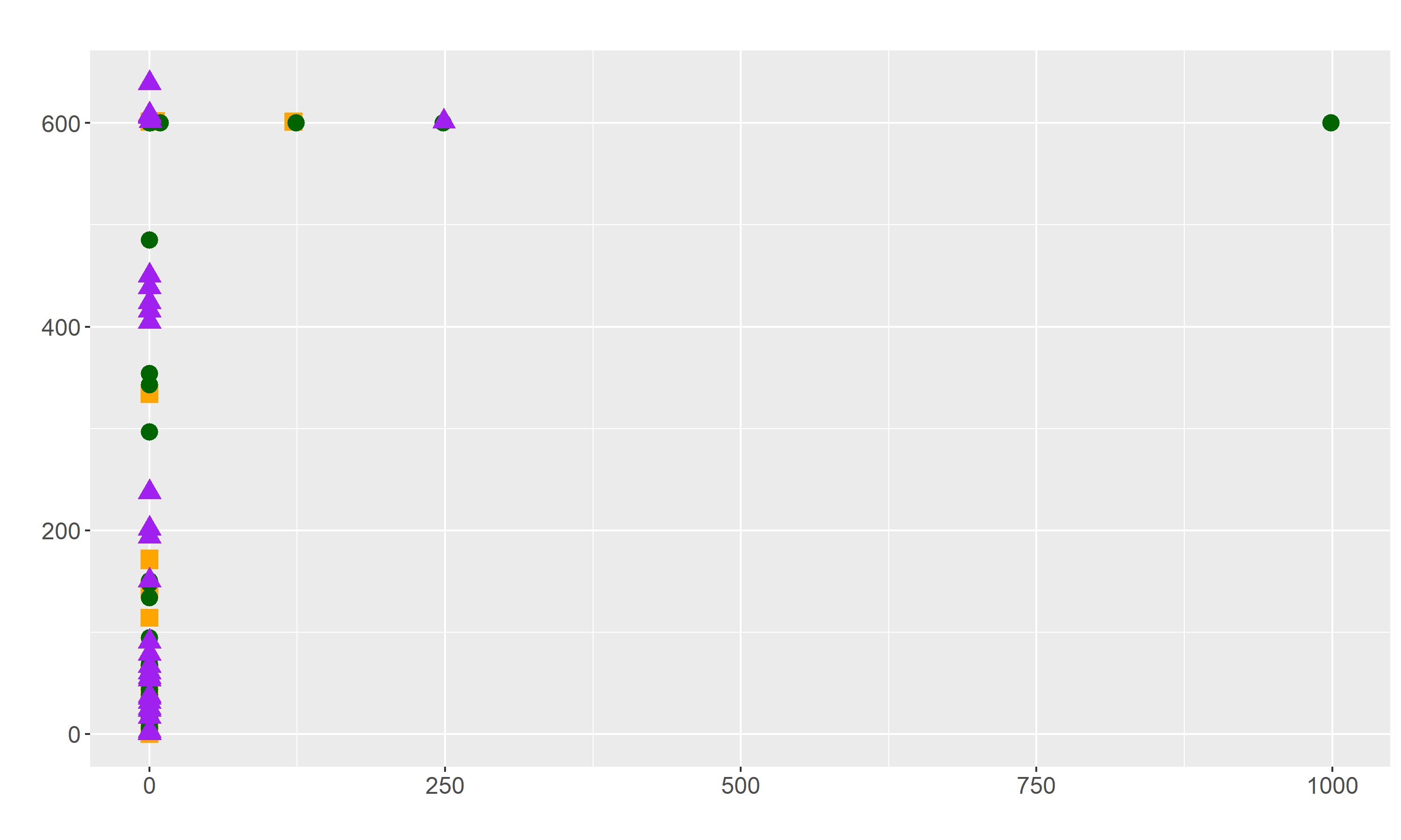}};
			\node[rotate=90](y) at (-3.35,0){{\color{darkgray}\tiny{time}}};
			\node(x) at (0,-2){{\color{darkgray}\tiny{distance}}};
		\end{tikzpicture}
		\caption{Global solvers: Distance to ref. solution vs time.}
\end{subfigure}
\begin{subfigure}{0.495\textwidth}
		\centering
		\begin{tikzpicture}
			\node(s1) at (0,0){\includegraphics[width=0.9\textwidth]{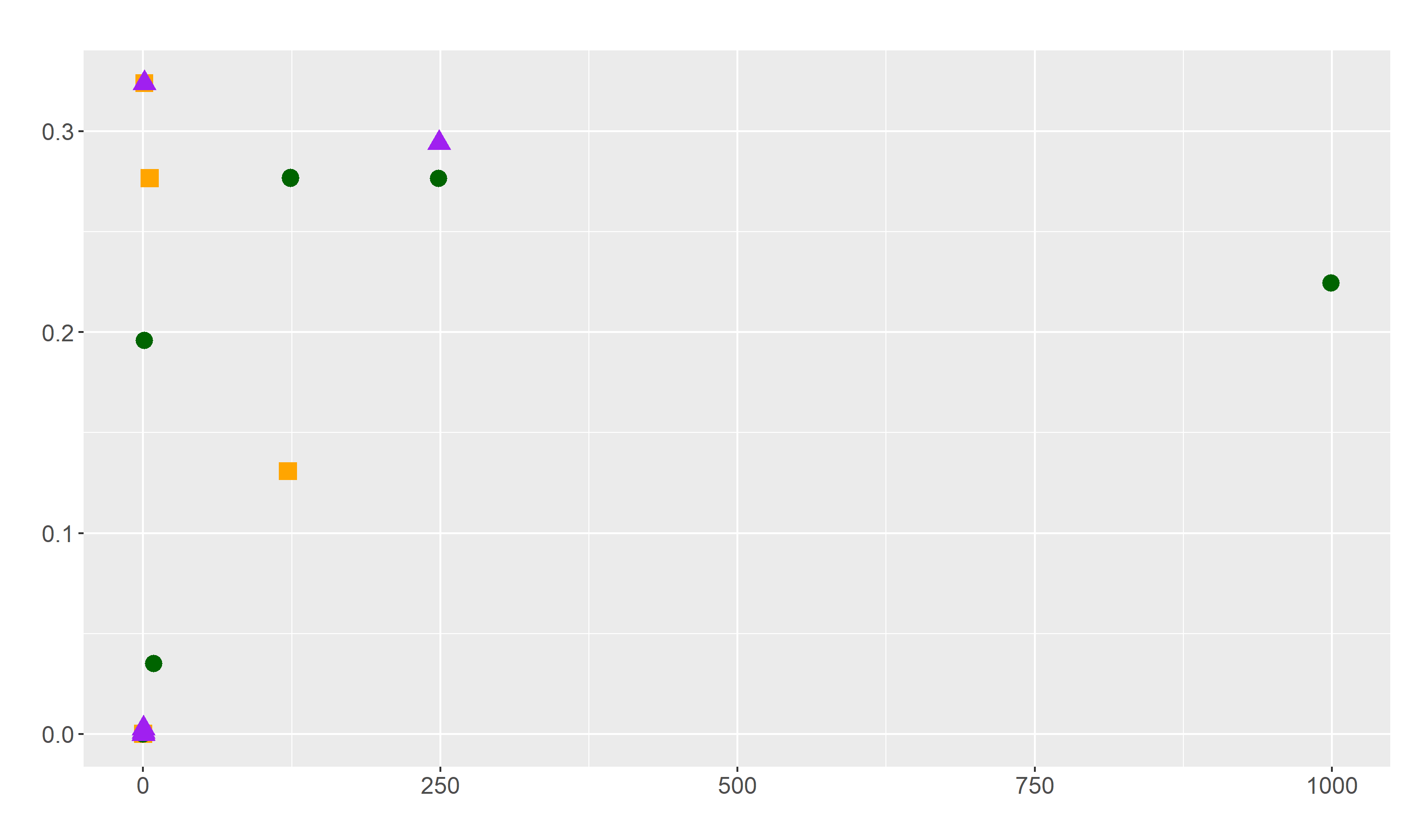}};
		\node(s2) at (2,0){\includegraphics[width=1cm]{figs/leyend_global.png}};
		\node[rotate=90](y) at (-3.4,0){{\color{darkgray}\tiny{$\log(\text{error}+1)$}}};
		\node(x) at (0,-2){{\color{darkgray}\tiny{distance}}};
		\end{tikzpicture}
		\caption{Global solvers: Distance to ref. solution vs error.}
\end{subfigure}

\vspace{0.2cm}

\begin{subfigure}{0.495\textwidth}
		\centering
		\begin{tikzpicture}
			\node(s1) at (0,0){\includegraphics[width=0.9\textwidth]{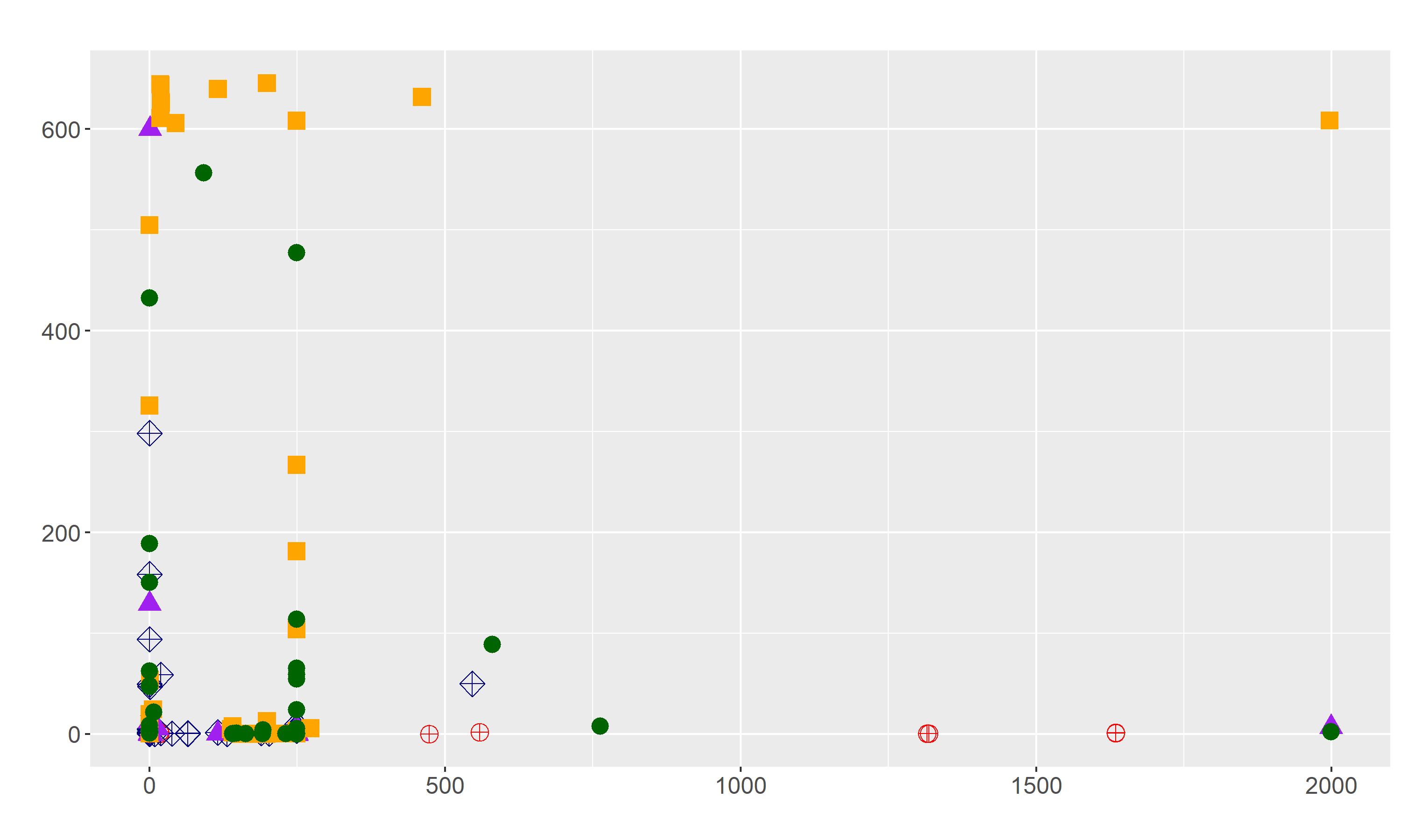}};
			\node[rotate=90](y) at (-3.35,0){{\color{darkgray}\tiny{time}}};
			\node(x) at (0,-2){{\color{darkgray}\tiny{distance}}};
		\end{tikzpicture}
		\caption{Local solvers: Distance to ref. solution vs time.}
\end{subfigure}
\begin{subfigure}{0.495\textwidth}
		\centering
		\begin{tikzpicture}
			\node(s1) at (0,0){\includegraphics[width=0.9\textwidth]{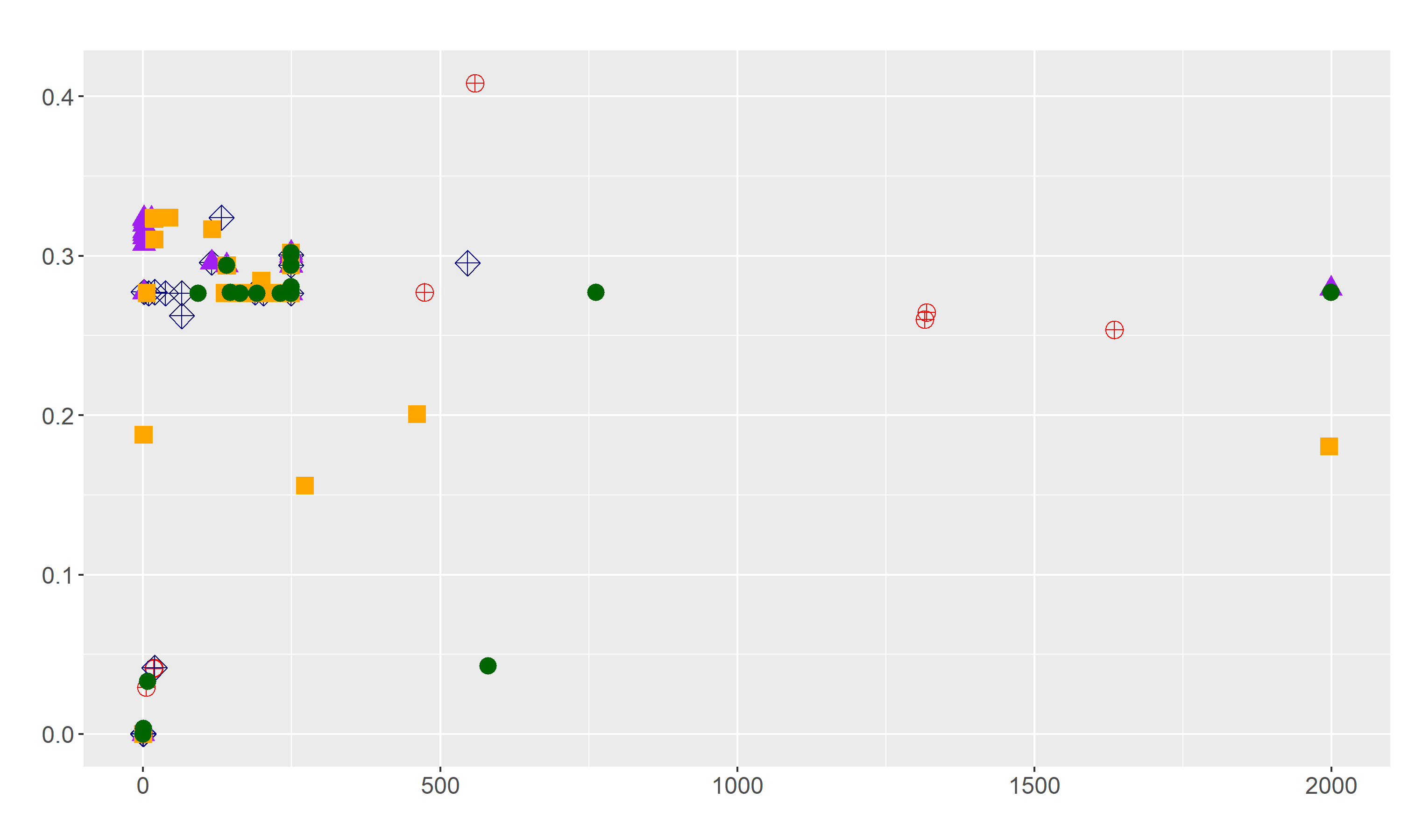}};
			\node(s2) at (2,0){\includegraphics[width=1cm]{figs/leyend_local.png}};
			\node[rotate=90](y) at (-3.4,0){{\color{darkgray}\tiny{$\log(\text{error}+1)$}}};
			\node(x) at (0,-2){{\color{darkgray}\tiny{distance}}};
		\end{tikzpicture}
		\caption{Local solvers: Distance to ref. solution vs error.}
\end{subfigure}

\begin{subfigure}{0.495\textwidth}
		\centering
		\begin{tikzpicture}
			\node(s1) at (0,0){\includegraphics[width=0.9\textwidth]{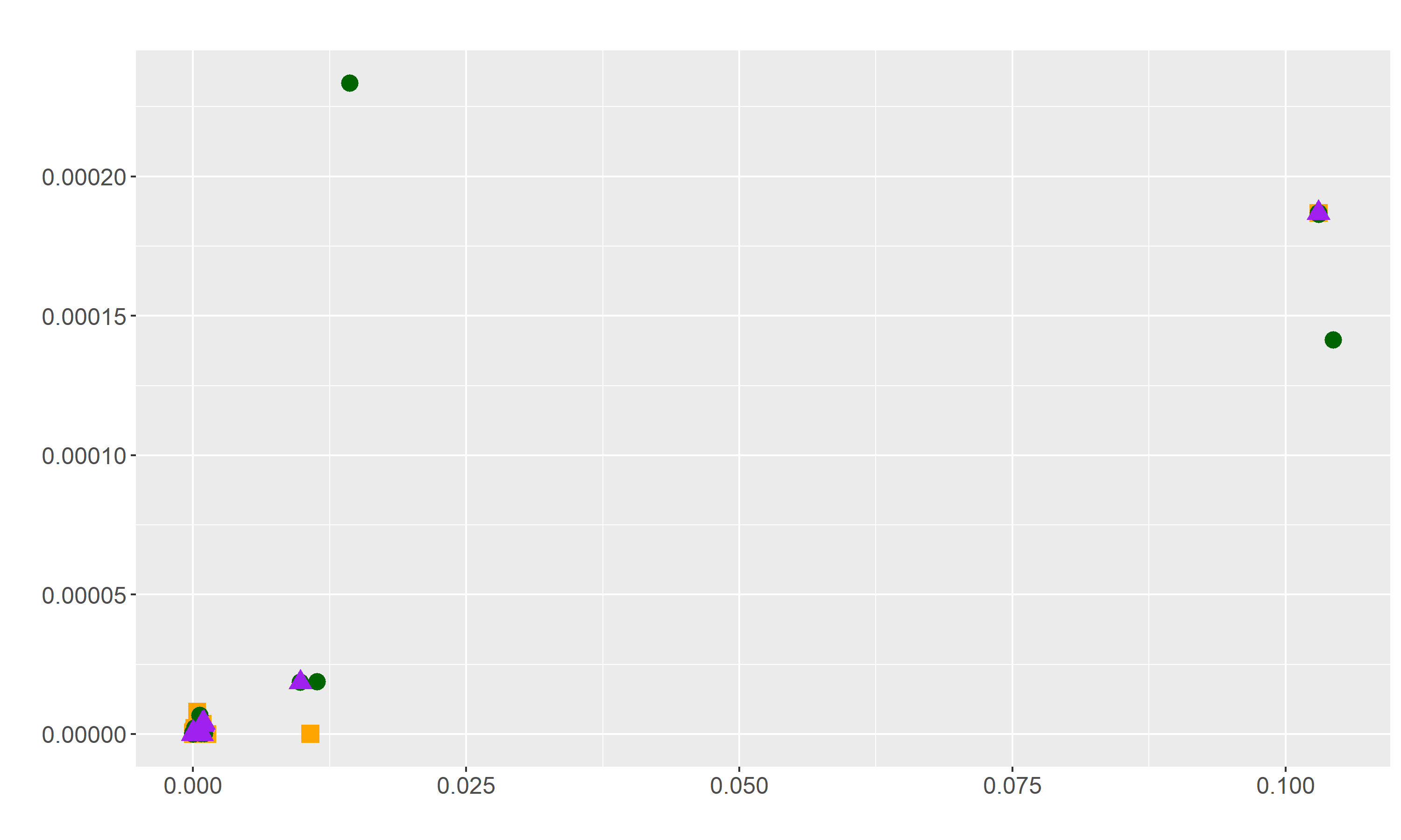}};
			\node(s2) at (2,0){\includegraphics[width=1cm]{figs/leyend_global.png}};
			\node[rotate=90](y) at (-3.4,0){{\color{darkgray}\tiny{$\log(\text{error}+1)$}}};
			\node(x) at (0,-2){{\color{darkgray}\tiny{distance}}};
		\end{tikzpicture}
		\caption{Global solvers: Distance to ref. solution vs error when ``solved''.}
\end{subfigure}

\caption{Results for problem \BBG.}
\label{fig:BBG}
\end{figure}
\end{document}